\documentclass[notitlepage,10pt, twoside]{article}

\usepackage{amssymb}
\usepackage{amsmath}
\usepackage{amscd}
\usepackage{latexsym}
\usepackage{epic}
\usepackage{eepic}
\usepackage{mathrsfs}
\usepackage[mathcal]{euscript}
\usepackage{eufrak}
\usepackage[all]{xy}

\newcommand{\vq}{(  V,   q)}

\newcommand{\vqpr}{(  V,   q')}

\newcommand{\covq}{  C_0(  V,  q)}
\newcommand{\covqpr}{  C_0(  V,  q')}

\newcommand{\simvqvqpr}{\hbox{\rm Sim}[\vq,\vqpr]}
\newcommand{\isomvqvqpr}{\hbox{\rm Iso}[\vq,\vqpr]}
\newcommand{\sisomvqvqpr}{\hbox{\rm S-Iso}[\vq,\vqpr]}

\newcommand{\isomcovqcovqpr}{\hbox{\rm Iso}[\covq,\covqpr]}
\newcommand{\isomcovqcovqprprime}{\hbox{\rm Iso}'[\covq,\covqpr]}
\newcommand{\sisomcovqcovqpr}{\hbox{\rm S-Iso}[\covq,\covqpr]}

\newcommand{\ovq}{{\hbox{\rm O}\vq}}

\newcommand{\vqi}{(  V,   q, I)}
\newcommand{\vprqpripr}{(  V',   q', I')}
\newcommand{\vprqpri}{(  V',   q', I)}
\newcommand{\vqpri}{(  V,   q', I)}

\newcommand{\covqi}{  C_0(  V,  q, I)}
\newcommand{\covqpri}{  C_0(  V,  q', I)}
\newcommand{\covprqpripr}{  C_0(  V',  q', I')}

\newcommand{\simvqivprqpripr}{\hbox{\rm Sim}[\vqi,\vprqpripr]}

\newcommand{\isomvqvprqpri}{\hbox{\rm Iso}[\vqi,\vprqpri]}
\newcommand{\simvqvqpri}{\hbox{\rm Sim}[\vqi,\vqpri]}
\newcommand{\isomvqvqpri}{\hbox{\rm Iso}[\vqi,\vqpri]}
\newcommand{\sisomvqvqpri}{\hbox{\rm S-Iso}[\vqi,\vqpri]}

\newcommand{\simvqvqi}{\hbox{\rm Sim}[\vqi,\vqi]}
\newcommand{\isomvqvqi}{\hbox{\rm Iso}[\vqi,\vqi]}
\newcommand{\sisomvqvqi}{\hbox{\rm S-Iso}[\vqi,\vqi]}

\newcommand{\isomcovqicovprqpripr}{\hbox{\rm Iso}[\covqi,\covprqpripr]}

\newcommand{\isomcovqcovqpri}{\hbox{\rm Iso}[\covqi,\covqpri]}
\newcommand{\isomcovqcovqprprimei}{\hbox{\rm Iso}'[\covqi,\covqpri]}
\newcommand{\sisomcovqcovqpri}{\hbox{\rm S-Iso}[\covqi,\covqpri]}

\newcommand{\sovqi}{{\hbox{\rm SO}}\vqi} 
\newcommand{\ovqi}{{\hbox{\rm O}\vqi}} 
\newcommand{\govqi}{{\hbox{\rm GO}\vqi}}

\newcommand{\autcovqi}{{\hbox{\rm Aut}}(\covqi)}   

\newcommand{\sautcovqi}{{\hbox{\rm S-Aut}}(\covqi)}   
\newcommand{\autcovqprimei}{{\hbox{\rm Aut}'}(\covqi)}   

\newcommand{\isolambdatwovivpripr}{{\hbox{\rm Iso}}[{\Lambda}^2(V)\otimes I^{-1},{\Lambda}^2(V')\otimes {(I')}^{-1}]}

\newcommand{\autdetlambdatwovi}{\hbox{\rm Aut}[\hbox{\rm det}(\Lambda^2(V)\otimes I^{-1})]}

\newcommand{\falgW}{\hbox{{\sf Alg}}_{{W}}}
\newcommand{\algW}{\hbox{\rm Alg}_{{W}}}

\newcommand{\faaWw}{{\hbox{\sf Assoc}_{{{W,w}}}}}
\newcommand{\aaWw}{{\hbox{\rm Assoc}_{{{W,w}}}}}
\newcommand{\stabw}{\hbox{\rm Stab}_w}
\newcommand{\gl}[1]{\hbox{\rm GL}_{#1}}
\newcommand{\fazuWw}{\hbox{\sf Azu}_{{W,w}}}
\newcommand{\azuWw}{\hbox{\rm Azu}_{{W,w}}}
\newcommand{\azuWiwi}{\hbox{\rm Azu}_{{W_i,w_i}}}
\newcommand{\azuWprwpr}{\hbox{\rm Azu}_{{W',w'}}}

\newcommand{\spazuWw}{\hbox{\rm SpAzu}_{{W,w}}}
\newcommand{\spazuWiwi}{\hbox{\rm SpAzu}_{{W_i,w_i}}}
\newcommand{\spazuWprwpr}{\hbox{\rm SpAzu}_{{W',w'}}}

\newcommand{\falgw}{\hbox{{\sf Alg}}_{W}}

\newcommand{\spec}[1]{\hbox{Spec\,}\left({#1}\right)}

\newcommand{\sym}[2]{\hbox{\sf Sym}_{#1}\left[\left(
{{#2}_{#1}}^{\vee}\otimes_{#1}{{#2}_{#1}}^{\vee}\otimes_{#1}{{#2}_{#1}}\right)^{\vee}\right]}

\newcommand{\closure}[1]{\langle{#1}\rangle}

\newcommand{\ncs}{{\cal N}_C^s(4,0)}
\newcommand{\nc}{{\cal N}_C(4,0)}

\newcommand{\ucss}{{\cal U}_C^{ss}(2,0)}
\newcommand{\ucs}{{\cal U}_C^s(2,0)}

\newtheorem{theorem}{Theorem}[section] 
\newtheorem{remark}[theorem]{Remark}
\newtheorem{definition}[theorem]{Definition} 
\newtheorem{lemma}[theorem]{Lemma} 
\newtheorem{proposition}[theorem]{Proposition} 
\newtheorem{corollary}[theorem]{Corollary} 

\begin{document}                                                                                   

\title{Line-Bundle-Valued Ternary Quadratic Bundles Over Schemes \\
{\small\bf Dedicated to the Memory of Professor Martin Kneser}}

\author{Venkata Balaji THIRUVALLOOR EESANAIPAADI\\
Mathematisches Institut, Georg-August-Universit\"at\\ 
Bunsenstrasse 3-5, D-37073 G\"ottingen, Germany\\
e-mail: {\tt tevbal@uni-math.gwdg.de}\\}

\date{July, 2005}

\maketitle 

\begin{abstract}
We study degenerations of rank 3 quadratic forms using those of rank 4
Azumaya algebras, and extend what is known for good forms and Azumaya
algebras. By considering line-bundle-valued forms, we extend the theorem
of Max-Albert Knus that the Witt-invariant---the even-Clifford algebra
of a form---suffices for classification. The general, usual and
special orthogonal groups of a form are determined in terms of
automorphism groups of its Witt-invariant. Martin Kneser's
characteristic-free notion of semiregular form is used. Examples of
non-existence of good forms and Azumaya structures are given. 
\end{abstract}




{\tableofcontents}

\section{Introduction}
To explain the results of this work, in this introduction we shall
 restrict ourselves to the case of affine schemes for the sake of
 simplicity, leaving the general formulation to \S  \ref{sec1}. So
 let $R$ be a commutative ring with $1.$ The objects of our study are
 ternary quadratic forms $q:V\rightarrow R$ on projective $R$-modules $V$ of
 rank 3. Zariski-locally these are quadratic forms in three
 variables. We shall denote such forms by pairs $(V,q).$ 
A similarity from $(V,q)$ to
 $(V',q')$ is a pair $(g,m)$, where $g:V\cong V'$ and 
$m:R\cong R$ are linear isomorphisms such that $q'\circ g=m\circ q.$ 
When $m$ is the identity, we call the similarity an isometry. 

Our broad aim is to study the sets of similarities between
 ternary forms (when they exist) in terms of their even-Clifford
 algebras and indicate applications. We are especially interested in
 bad quadratic forms, and would like to develop a theory
 for such forms. We obtain results extending (though not deduced from,
 but rather
 motivated by) what is known for the good ones. 
 
What are good forms? These are the same as the usual regular (or
nonsingular) forms, provided none of the residue fields of $R$ is of
 characteristic two. The correct technical notion is that of {\em
 semiregular} form, due to Martin Kneser, on which we will elaborate
 in \S\ref{semireg-bil-forms}, page~\pageref{semireg-bil-forms}.
 At this point we only remark that
 this notion works universally, so includes the case of regular forms
 as well. 

For a semiregular form $q:V\longrightarrow R$, its even-Clifford algebra
$A=C_0(V,q)$ is an Azumaya $R$-algebra with underlying module of rank
4.
 Azumaya means that the natural $R$-algebra homomorphism 
$$A\otimes_R A^{\textrm{op}}\rightarrow\textrm{End}_{R-\textrm{Mod}}(A): 
a\otimes b^{\textrm{op}}\mapsto (x\mapsto axb^{\textrm{op}})$$
is an isomorphism, where the superscript refers to the opposite
algebra.
 Matrix algebras, for example algebras of $(2\times 2)$ 
matrices in our case, are simple examples of Azumaya algebras, and
they are all one gets if say $R$ were an algebraically closed
field. 
Any similarity $(g,m):(V,q)\cong (V',q')$ induces an
isomorphism of $R$-algebras $C_0(g,m):C_0(V,q)\cong C_0(V',q').$ So if
we denote the set of isometry classes of semiregular 
 ternary quadratic forms by
${\cal Q}_3^{sr}(R)$, and the set of $R$-algebra isomorphism classes of 
rank 4 Azumaya $R$-algebras by ${\cal AZU}_4(R)$, then we get a natural
map, which is in fact functorial in $R$: 
$$(\dag)\hspace*{0.25in} {\cal Q}_3^{sr}(R)\longrightarrow {\cal AZU}_4(R) \textrm{ via } (V,q)\mapsto C_0(V,q).$$
It is well-known that the above map is surjective but may not be injective. 
In fact, given an Azumaya algebra $A$ of rank 4, there is a
well-defined unique standard involution $\sigma_A$ on $A$, which
defines a regular quadratic form: the norm $n_A:A\longrightarrow R$
and a linear form: the trace $tr_A:
A\longrightarrow R.$ Then the restriction of $n_A$ to the kernel $A'$ of
$tr_A$, which is a rank 3 $R$-module (due to the surjectivity of the
trace),
 remains good i.e., is a semiregular quadratic form, and one
has an isomorphism of $R$-algebras 
$$C_0(A', n_A|A')\cong A.$$
Thus the map $(\dag)$ is surjective. Let us explain why it is
not injective. Let $I$ be a rank 1 projective $R$-module (i.e., an
invertible module) equipped with an $R$-linear isomorphism  
$$h:I^2:=I\otimes I\cong R.$$ The pair $(I,h)$ is called a {\em
discriminant module} and can be thought of as a good symmetric
bilinear form on $I$ as well as a prescribed square root of $R.$
Now given a ternary quadratic module $(V,q)$, we may define a new
quadratic module: 
$$(V,q)\otimes (I, h):=(V\otimes I, q\otimes h): (q\otimes h)(v\otimes
l):=q(v)h(l\otimes l).$$ Then it can be verified that we have an
isomorphism of $R$-algebras 
$$C_0(V,q)\cong C_0(V\otimes I, q\otimes
h).$$
 Thus the map in $(\dag)$ may  not be injective. Now the set
of isomorphism classes of $(I,h)$ forms a commutative group of
exponent 2, the group operation induced by the usual tensor
product. It is denoted by $\textrm{Disc}(R)$ and acts on the left side
of $(\dag).$ It is a Theorem of Max-Albert Knus that 
 we have a bijection 
$$(\dag\dag)\hspace*{0.25in} {\cal Q}_3^{sr}(R)/\textrm{Disc(R)} \cong 
{\cal AZU}_4(R),$$ which may also be thought of as a statement in
cohomology (cf.~discussion following Theorem~\ref{Max-Knus-Theorem},
page~\pageref{coho-interpretation}).
 A natural question is to
ask if there is a statement analogous to 
$(\dag\dag)$  if we consider  
non-semiregular forms on the left side as well. Our central result,
Theorem~\ref{bijectivity} (page~\pageref{bijectivity}),
 shows that the answer is yes, thus providing a
limiting version of the cohomological statement
$(\dag\dag).$ However, when considering non-semiregular forms, 
 we would need to make some changes to our present formulation; 
these we set out to explain next. 

Firstly, since we wanted to include degenerate forms in the left side
of $(\dag)$, and since such forms are limits, locally in the
Zariski-topology, of semiregular forms, it is natural to consider
similar limits on the right side. In other words, we consider rank 4
$R$-algebras which are schematic limits of rank 4 Azumaya
algebras. They were introduced in \cite{tevb-paper1}, and are recalled
in \S\ref{subsec2.7} (page~\pageref{subsec2.7}).
 Instead of going now into the definition of such
limiting algebras, we note that it follows from {\em op.~cit.,} that
they are precisely those rank 4 $R$-algebras which are
Zariski-locally isomorphic to even-Clifford algebras of ternary
quadratic forms. Let us denote the set of isomorphism classes of such
algebras by ${\cal SPAZU}_4(R).$ This will be the replacement for
the right side of $(\dag\dag)$. 

Given an algebra $A$ representing an element of ${\cal
SPAZU}_4(R)$, all we know by definition is that there are elements $f_i$
that generate $R$, and ternary quadratic $R_{f_i}$-forms $(V_i, q_i)$
such that there are isomorphisms $A\otimes_R{R_{f_i}}\cong C_0(V_i,
q_i).$ It is not clear {\em a priori} that there is a choice for which
the $V_i$ would glue to give a rank 3 $R$-module $V$, and even if this were
so, that the $q_i$ would glue to give a quadratic form on $V.$ In
fact, it may not. We show that there is a
choice for which the $V_i$ glue to give a $V$,  but that the $q_i$ glue to
give a quadratic form $q$ {\em with values in the invertible
module } $I:=\textrm{det}^{-1}(A)$, such that 
 $$C_0(V,q,I)\cong A,$$ (cf.~Theorem
\ref{structure-of-specialisation}, part (a),
page~\pageref{structure-of-specialisation}) and further that that we
could get a
quadratic form with values in $R$ iff $\textrm{det}(A)\in
2.\textrm{Pic}(R)$ (cf.~Theorem \ref{structure-of-specialisation}, part
(b)). At this point we note that the even Clifford algebra
$C_0(V,q,I)$ is the one defined by Bichsel-Knus in
\cite{Bichsel-Knus}, briefly recalled in \S\ref{subsec2.3},
page~\pageref{subsec2.3}, which reduces to the familiar even-Clifford
algebra when the invertible module is $R$ itself. 
Thus, going back to our idea of replacing the left side of
$(\dag)$ with degenerate forms, we see that we should work with
ternary quadratic forms with values in invertible
$R$-modules. Let us denote the set of isometry classes of such
quadratic modules by ${\cal Q}_3(R).$
  Further, since we had to divide out by the action of
discriminant modules to get the bijection $(\dag\dag)$,
and since such discriminant modules are square roots of $R$, it is
natural that we should consider square roots of invertible $R$-modules
as well. Such objects could be called {\em twisted} discriminant
modules, and they do form a group under the tensor product which we
denote by $\textrm{T-Disc}(R).$ Any such object is given by a triple 
$(L, h, J)$ where $h:L\otimes L\cong J$ is an $R$-linear isomorphism
with $L, J$ being invertible $R$-modules ($h$ makes $L$ a specified
square-root of $J$). Given a quadratic form $q:V\rightarrow I$ with
values in the invertible module $I$, one defines in the obvious way: 
$$(V,q,I)\otimes(L, h, J):=(V\otimes L, q\otimes h, I\otimes
J):v\otimes l\mapsto q(v)\otimes h(l\otimes l),$$ and there is a 
canonical isomorphism of $R$-algebras 
$$C_0(V\otimes I, q\otimes h, L\otimes J)\cong C_0(V, q, I).$$
 It happens that this generalisation to
values in invertible modules does not disturb the bijection
$(\dag\dag)$ i.e., if we considered semiregular
quadratic modules with values in invertible modules on the left side
modulo the action of $\textrm{T-Disc}(R)$, we still obtain a bijection
to the right side. Theorem~\ref{bijectivity}
(page~\pageref{bijectivity}) shows that we do have a bijection 
$$(\ddag)\hspace*{0.25in}{\cal Q}_3(R)/\textrm{T-Disc}(R)\cong {\cal
SPAZU}_4(R)$$
which is a limiting version of $(\dag\dag).$ We next indicate the
ingredients that go into the proof of $(\ddag)$, which are of interest
by themselves, since they involve the study of the sets of
similarities between ternary quadratic modules. 

Given ternary quadratic forms $q,q':V\rightarrow I$, we may consider
the set $\textrm{Sim}(q,q')$
of similarities from $q$ to $q'$: these are pairs $(g,l)$
where $g\in\textrm{GL}(V)$ and $l\in R^*$ that satisfy
$q'g=\lambda_lq$ where $\lambda_l$ is the element of $\textrm{Aut}(I)$
given by multiplication by $l.$ The second member $l$ of $(g,l)$ 
 is called the {\em multiplier} of the similarity. The set 
$\textrm{Sim}(q,q')$ is a subset of $\textrm{GL}(V)\times R^*.$ The subset 
$\textrm{Iso}(q,q')$ of
{\em orthogonal transformations} or {\em isometries} are those
similarities which have trivial multipliers. A smaller subset,
$\textrm{S-Iso}(q,q')$, consists 
of {\em special} orthogonal transformations, and corresponds to those for
which $\textrm{det}(g)=1.$ The latter two subsets may be
considered as subsets
of $\textrm{GL}(V)$ and $\textrm{SL}(V)$ 
respectively. 

Since the even-Clifford algebra is functorial in $q$,
any similarity $q\rightarrow q'$ gives rise to an isomorphism of
$R$-algebras $C_0(V,q,I)\cong C_0(V,q',I).$
Theorem~\ref{lifting-of-isomorphisms},
page~\pageref{lifting-of-isomorphisms}, identifies the images of the  
subsets of similarities (defined in the previous paragraph)
  in the set of isomorphisms
$\textrm{Iso}(C_0(q), C_0(q'))$ 
from $C_0(q)$ to $C_0(q').$ It shows firstly, that each isomorphism
of even-Clifford algebras lifts to a similarity; further, given such
an isomorphism $\phi$, Prop.~\ref{transfer-to-lambda2},
page~\pageref{transfer-to-lambda2}, shows how to define an $R$-linear 
automorphism $\phi_{\Lambda^2}$ of $\Lambda^2(V)\otimes I^{-1}.$ We
may thus consider those isomorphisms $\phi$ for which
$\textrm{det}(\phi_{\Lambda^2})\in R^*$ is a square, (respectively, is 
1), and denote it by $\textrm{Iso}'(C_0(q), C_0(q'))$ (respectively,
by $\textrm{S-Iso}(C_0(q), C_0(q')).$) Then
Theorem~\ref{lifting-of-isomorphisms} shows that the image of 
$\textrm{Iso}(q,q')$ is preciesly $\textrm{Iso}'(C_0(q),C_0(q'))$ and
that of $\textrm{S-Iso}(q,q')$ is precisely $\textrm{S-Iso}(C_0(q),
C_0(q')).$ 

These surjectivities are accomplished by explicitly constructing 
 families of sections indexed by the integers. The proof of
the existence of $\phi_{\Lambda^2}$ involves the use of 
{\em Bourbaki's Tensor Operations} for forms with values
in invertible modules (cf.~\S\ref{subsec2.4}). In order to construct the
sections, we make computations analysing  
 the correspondence from \cite{tevb-paper1} 
between locally-even-Clifford algebras
and specialisations of Azumaya algebras with the same underlying
module. The further uses of
Theorem~\ref{lifting-of-isomorphisms} are principally threefold: it is
used to prove the injectivity part of $(\ddag)$
(Theorem~\ref{bijectivity}) and later on in the proof
of the surjectivity part (Theorem~\ref{structure-of-specialisation},
(a)). It is also used 
 in Theorem~\ref{lifting-of-automorphisms},
 page~\pageref{lifting-of-automorphisms}, in the 
 description of the
general, usual and special orthogonal groups
$\textrm{GO}(q):=\textrm{Sim}(q,q)$,
$\textrm{O}(q):=\textrm{Iso}(q,q)$ and
$\textrm{SO}(q):=\textrm{S-Iso}(q,q)$ respectively, as fitting into a
neat commutative diagram. Further, $\textrm{GO}(q)$ is shown to be a
semi-direct product, and when $R$ is an integral domain, and $q$ is
semiregular even at one point of $\textrm{Spec}(R)$, then all the
subgroups of automorphisms of $C_0(q)$ described above are seen to
coincide, so that $\textrm{O}(q)$ is the semi-direct product of 
the multiplicative subgroup $\mu_2(R)$ consisting of the square roots
of 1 and of $\textrm{SO}(q).$  

The injectivity of Theorem~\ref{bijectivity} is the statement that two
ternary quadratic modules with isomorphic even-Clifford algebras
differ by an isometry after one of them is tensored with a suitable
twisted discriminant module. The existence of such a module is reduced
to the lifting of a given isomorphism of even-Cliffords to a similarity in
the free case (cf.\S\ref{sec3}). 

A suitable member of the family of sections constructed to show the
surjectivities in the proof of Theorem~\ref{lifting-of-isomorphisms}
allows one to get a cocycle defining a ternary quadratic module from a
suitably chosen cocycle for a given specialised algebra
(cf.\S\ref{sec5}). 

Since Bourbaki's Tensor Operations form a key tool in our
computations, we explain how they are of use to us. Given two ternary
quadratic forms $q,q':V\rightarrow I$, one may ask how one would
compare, if at all, the underlying modules of their even-Clifford
algebras. Now a bilinear form $b:V\times V\rightarrow I$ defines a
quadratic form $q_b:V\rightarrow I$ by $q_b(v):=b(v,v).$ If there
exists a $b$ so that $q_b=q-q'$, then the Tensor Operations 
provide us with an explicitly computable  $R$-linear isomorphism 
$$\psi_b: C_0(V,q, I)\cong C_0(V,q', I).$$ 
Since we can always find a $b$ such that $q_b=q$, by taking $q'=0$ we
see that we could, via $\psi_b$, transport the algebra structure on
$C_0(V,q, I)$ to one on the underlying module $W$ of 
$$C_0(V, 0,I)=R\oplus \Lambda^2(V)\otimes I^{-1}$$
with multiplicative identity $w=1_R.$ In this way, we may map the
space of $I$-valued bilinear forms on $V$ into the space of
$R$-algebra structures on the {\em same} underlying module $W$ with
identity $w.$ And we could do this functorially in $R$-algebras $S.$
Theorem~\ref{bilinear-forms-as-specialisations},
page~\pageref{bilinear-forms-as-specialisations}, shows that in this
way we obtain bijectively for any  $S$, all the locally-even-Clifford
$S$-algebra structures on $W\otimes S$ with identity $w\otimes S.$ In
other words, bilinear forms {\em correspond} to such algebras. 
This generalizes a similar result due to S.~Ramanan for the case $R=k$
an algebraically closed field (cf.~\S2.3 of \cite{phdthesis}). When
$\textrm{char.}(k)\neq 2$, it was given another interpretation by
C.~S.~Seshadri in \cite{css-desing} (cf.~Theorem~\ref{Thetaisom},
page~\pageref{Thetaisom}), which will also be our main key to
computations. 

Not only in the proof of
Theorem~\ref{bilinear-forms-as-specialisations}, but throughout this
work, we use the description from \cite{tevb-paper1} of the
locally-even-Clifford algebras as schematic specialisations of Azumaya
algebras since it is this description that allows us to compute
without losing track of global information. 

For the rest of this introduction, we depart to the language of
schemes. Suppose $W$ is a vector bundle of rank 4 over a scheme
$X.$ If there exists an Azumaya algebra structure on $W$, then it can
be seen that $W$ has to be self-dual. One may ask for a converse: if
$W$ is self dual, does there exist a global Azumaya algebra structure
on $W$? We answer this question more generally, and 
show in Theorem~\ref{relations-with-self-duality},  
page~\pageref{relations-with-self-duality}, 
that this is indeed the case if we on the one hand let $X$ 
to be a proper scheme of finite type over a base scheme with connected fibres,
 while on the other hand we assume that $W$ has
square rank $n^2$ with $n\geq 2$, and moreover that there exists some
associative unital algebra structure $A$  on $W$ which is Azumaya at
some point of each fiber over the base.
 We then show in fact that $A$ has to be Azumaya,
giving a ``punctual to global'' result.  

By an application of Theorem~\ref{bijectivity} to this result, we
deduce  
 that the
hypothesis of self-duality on the underlying bundle of the
even-Clifford algebra of a ternary quadratic bundle implies that the
quadratic bundle is semiregular everywhere if it is semiregular even
at a single point of each fiber over the base scheme. 

Such results can be used to give examples of
rank 4 vector bundles $W$ over which there do not exist any Azumaya
algebra structures but which nevertheless admit algebra structures
that are Azumaya on a nonempty open subscheme. They can also be used 
to give examples of pairs $(V,I)$, 
consisting of a rank 3
vector bundle $V$ and a line bundle $I$, such that there exists no semiregular
quadratic form on $V$ with values in $I$, though there do exist forms
that are semiregular on a nonempty open subscheme.  
 For any given irreducible smooth
projective curve $C$ of genus $g\geq 2$ over an algebraically closed field,  
 such examples  
naturally occur on certain smooth projective varieties $\nc$
 of dimension $4g-3$, where $g$ is the genus of $C$. In fact these $\nc$  
 appear as fine moduli spaces of certain parabolic stable vector
bundles of rank 4 and degree zero over $C$ and arise as Seshadri
desingularisations of the moduli space of semistable vector bundles of
rank 2 and degree 0 over $C.$ The whole situation can be generalised
to the case of a curve relative to a base scheme and the results are
given in Theorem~\ref{non-existence-examples},
page~\pageref{non-existence-examples}. A detailed account with proofs will
appear in \cite{tevb-paper3}. 

Suppose $X$ is a scheme and $W$ a rank 4
vector bundle on $X$ with a nowhere-vanishing global section $w.$ In
Theorem~\ref{surjectivity-affine-limited},
page~\pageref{surjectivity-affine-limited}, we delve a little into the
topological and geometric properties of the scheme of specialised
algebra structures on $W$ with identity $w$ (functorially in
$X$-schemes). If $X$ is locally-factorial, we show that the natural
homomorphism from the Picard group of
$X$ to that of the scheme of specialised algebra structures on $W$
with unit $w$ is an isomorphism. 
We also study the specialised algebras that are
nowhere Azumaya. When $X$ is locally-factorial and $W$ is
self-dual, these algebras define a Weil divisor $D_X$ such
that some positive integer multiple $n.D_X$ is principal. Therefore
the natural surjective homomorphism from the Picard group of the
scheme of specialisations to that of the open subscheme of Azumaya
algebra structures is an isomorphism iff $n=1$ and in this situation, 
 every specialised
algebra structure on $W\otimes T$ arises from a quadratic form with
values in the trivial line bundle, for any $X$-scheme $T.$   
For example, this happens if $X$ is 
locally-factorial and if there exists an Azumaya algebra structure on
$W.$
 The proof uses the
stratification of the variety of specialisations of $(2\times
2)$-matrix algebras over an algebraically closed field
(Theorem~\ref{stratification}, page~\pageref{stratification}). 

\paragraph*{\bf Plan of the Paper.}
 Background material is to be found in
 \S\ref{sec2}, which also fixes some notation, explains definitions
  and recalls results for future use.
The formulation of the statements of the main results follow in
 \S\ref{sec1}.
 The proofs will occupy  
 \S\ref{sec3} through \S\ref{sec7}.
\paragraph*{\bf Note on Numbering.}
Various items such as theorems, propositions, definitions etc are all 
consecutively numbered within a section, irrespective of the
subsection that they may appear in.

\section{Notations and Preliminaries}
\label{sec2}
This section collects together definitions and results necessary for
future use.
We omit proofs for statements involving quadratic and bilinear 
 forms, since these can be reduced to the corresponding 
 results for the affine case which are treated in Knus' book \cite{Knus}. 
 For the systematic treatment of the generalised 
 Clifford algebra and its properties, we refer the reader to the paper of 
Bichsel-Knus \cite{Bichsel-Knus}. We also collect some preliminaries 
 on the notion of schematic image and recall some 
 background material from Part A of \cite{tevb-paper1}. 
\subsection{Quadratic and Bilinear Forms with Values in a Line Bundle} 
\label{subsec2.1}
 Let $V$ be a vector bundle (of constant positive rank)
 and  $I$ a line bundle on a scheme $X.$ Let $\cal V$ and $\cal I$
 respectively
 denote the coherent locally-free sheaves corresponding to $V$ and
 $I.$ We define the coherent locally-free sheaves of bilinear and 
 alternating forms on $V$ with values in $I$ respectively as follows: 
$$\hbox{\rm Bil}_{({\cal V}, {\cal I})}
 :=(T^2_{{\cal O}_X}({\cal V}))^\vee\otimes{\cal
 I}\textrm{\hspace*{3mm}and\hspace*{3mm}}
\hbox{\rm Alt}^2_{({\cal V}, {\cal I})}
 :=(\Lambda^2_{{\cal O}_X}({\cal V}))^\vee\otimes{\cal I}.$$
We let $\hbox{\rm Bil}_{(V,I)}$ and $\hbox{\rm Alt}^2_{(V,I)}$ 
 denote the corresponding vector bundles. 
Now define the (coherent, locally-free) sheaf of $I$-valued quadratic 
 forms on $V$  by the exactness of the following sequence: 
$$0\longrightarrow \hbox{\rm Alt}^2_{({\cal V},{\cal I})}
 \longrightarrow \hbox{\rm Bil}_{({\cal V},{\cal I})} \longrightarrow
 \hbox{\rm Quad}_{({\cal V},{\cal I})}\longrightarrow 0.$$ Let
 the corresponding bundle  of $I$-valued 
 quadratic forms 
 on $V$ be denoted by $\hbox{\rm Quad}_{(V,I)}.$
 Thus a
 bilinear form (resp. alternating form, resp. quadratic form) with values in
 $I$ on $V$ over an open set $U\hookrightarrow X$  is by 
 definition a section over $U$ of the vector bundle 
$$\hbox{\rm Bil}_{(V, I)}\textrm{ (resp. of }\hbox{\rm Alt}^2_{(V,I)},
\textrm{ resp. of }\hbox{\rm Quad}_{(V,I)}),$$ or equivalently, an 
 element of 
$$\Gamma(U, \hbox{\rm Bil}_{({\cal V}, {\cal I})})\textrm{ (resp. 
 of }\Gamma(U, \hbox{\rm Alt}^2_{({\cal V}, {\cal I})}),\textrm{ resp. 
 of }\Gamma(U, \hbox{\rm Quad}_{({\cal V}, {\cal I})})).$$ 
 In terms of the corresponding (geometric) vector bundles over $X$,
 the 
 exact sequence above 
 translates into  the following 
 sequence of morphisms of vector bundles, with the first one a closed immersion and the second one a Zariski locally-trivial principal 
$\hbox{\rm Alt}^2_{(V, I)}$-bundle: 
$$\hbox{\rm Alt}^2_{(V,I)} \hookrightarrow \hbox{\rm Bil}_{(V,I)} \twoheadrightarrow \hbox{\rm Quad}_{(V,I)}.$$ 
Given a quadratic form $q\in\Gamma(U, \hbox{\rm Quad}_{({\cal V},{\cal I})})$,
 recall 
 that  the usual 
 `associated' symmetric 
 bilinear form $b_q\in\Gamma(U, \hbox{\rm Bil}_{({\cal V},{\cal I})})$ 
 is given on sections (over open subsets of $U$) by 
$$v\otimes v'\mapsto q(v+v')-q(v)-q(v').$$
 Given a (not-necessarily symmetric!) bilinear form
 $b$, we also have the induced quadratic form $q_b$ given on sections by 
 $v\mapsto b(v\otimes v).$
 A global quadratic form may not be induced 
 from a global bilinear form, unless we assume something more, for e.g., 
 that the scheme is affine, or more 
 generally  that the sheaf cohomology group
 $\hbox{\rm H}^1(X, \hbox{\rm Alt}^2_{({\cal V}, {\cal I})})=0$.
\subsection{Sets of Similarities of Quadratic Bundles}
\label{newsubsec2.2}
Let $\vqi$ and $\vprqpripr$ be quadratic bundles on the scheme $X.$
 We denote by $$\simvqivprqpripr$$ the set of generalised similarities 
 from  $\vqi$ to $\vprqpripr.$ These consist of pairs  
  $(g,m)$ such that 
 $g:V\cong V'$ and $m:I\cong I'$ are  linear isomorphisms and 
 the following diagram commutes (where $q$ and $q'$ are considered 
 as morphisms of sheaves of {\em sets}): $$\begin{CD}
{V} @>{g}>{\cong}> {V'}\\
@V{q}VV @VV{q'}V\\
I @>{\cong}>m> I'
\end{CD}$$
 When $I=I'$, since an $m\in\hbox{\rm Aut}(I)$ may be thought of as 
 multiplication by a scalar
 $l\in\Gamma(X, {\cal O}_X^*)\equiv\hbox{\rm Aut}(I)$, we may call the 
 isomorphism $(g,m)$ as an $I$-similarity with multiplier $l.$ In such 
 a case we may as well  
 denote $(g,m)$ by the pair $(g, l)$ and we often write 
 $$g:(V, q, I)\cong_{l}(V',q',I).$$
 Let  
 $$\isomvqvprqpri$$ be the subset of isometries (i.e., those pairs $(g,m)$ 
 with $m=\hbox{\rm Identity}$ or $I$-similarities with trivial 
 multipliers). 
  When $ V= V'$, the 
 subset of isometries with trivial determinant is denoted 
 $$\sisomvqvqpri.$$ On taking $ q= q'$
  these sets 
 naturally become subgroups of 
$$\hbox{\rm Aut}(V)\times \Gamma(X, {\cal O}_X^*)
=\hbox{\rm GL}(V)\times \Gamma(X, {\cal O}_X^*)$$
 and we define
\begin{center}\begin{tabular}{c}
$\simvqvqi=:\govqi\supset\isomvqvqi=:\ovqi\supset$\\
 $\supset\sisomvqvqi=:\sovqi.$ \end{tabular}\end{center}  
Of course, $\ovqi$ and $\sovqi$ may be thought of as subgroups of 
 $\hbox{\rm GL}(V)\equiv\hbox{\rm GL}(V)\times\{1\}$ and 
$\hbox{\rm SL}(V)\equiv\hbox{\rm SL}(V)\times\{1\}$ respectively.   
\subsection{Semiregular Bilinear forms}\label{semireg-bil-forms}
\label{subsec2.2}
 Fundamental problems in dealing with quadratic forms over
 arbitrary base schemes arise essentially from two abnormalities 
 in characteristic two: firstly,  
  the mapping that associates a quadratic form to its symmetric 
 bilinear form is not bijective and secondly, 
 there do not exist regular quadratic forms on any odd-rank bundle.
 The remedy for this is to consider  
 semiregular quadratic 
 forms, a concept due to M.Kneser \cite{Kneser} and 
 elaborated upon 
 by Knus in \cite{Knus}, which in fact works over an arbitrary base 
 scheme (and hence in a characteristic-free way) and further reduces to 
 the usual notion of regular form in characteristics $\neq 2.$

 Let $\hbox{\rm Spec}(R)=U\hookrightarrow X$ be an open affine subscheme 
 of $X$ such that $V|U$ is trivial. Let  $\hbox{\rm Quad}_V:=\hbox{\rm
   Quad}_{(V, {\cal O}_X)}.$
 Consider a  quadratic form $q\in\Gamma(U, \hbox{\rm Quad}_V)$ on $V|U$ 
 and  its associated
 symmetric bilinear form $b_q$. The matrix of this bilinear form relative
 to any fixed basis is a symmetric matrix of odd rank and in particular, 
 if $U$ is of pure characteristic two (i.e., the residue field of
 any point of $U$ is 
 of characteristic two), then this matrix is also alternating
 and is  hence singular, immediately implying that $q$ cannot be regular.
 However, computing the 
 the determinant of such a matrix  in {\em formal variables} $\{\zeta_i,
\zeta_{ij}\}$
  shows that it is twice the 
 following  polynomial 
$$P_{3}(\zeta_{i},\zeta_{ij})
=4\zeta_{1}\zeta_{2}\zeta_{3}+\zeta_{12}\zeta_{13}\zeta_{23}
-(\zeta_{1}\zeta_{23}^{2}+\zeta_{2}\zeta_{13}^{2}+\zeta_{3}\zeta_{12}^{2}).$$
 The value $P_3(q(e_i),b_q(e_i,e_j))$ corresponding to a basis 
 $\{e_1,e_2,e_3\}$ is called the 
{\em half-discriminant} of $q$ relative to that basis,  
 and $q$ is said to be semiregular if its half-discriminant is 
 a unit. It turns out that this definition is independent of the 
 basis chosen (\S 3, Chap.IV, \cite{Knus}).

 Even if $V|U$ were 
 only projective (i.e., locally-free but not free), the semiregularity 
 of $q$ may be defined as the semiregularity of $q\otimes_R R_{\mathfrak m}$ 
 for each maximal ideal ${\mathfrak m}\subset R$, and it turns out that with 
this definition, the notion of a quadratic form
 being semiregular is local and is 
 well-behaved under base-change (Prop.3.1.5, Chap.IV, \cite{Knus}).

 We may thus define the subfunctor of $\hbox{\rm Quad}_V:=\hbox{\rm
   Quad}_{(V, {\cal O}_X)}$
 of semiregular quadratic forms. This subfunctor is 
 represented by a $\hbox{\rm GL}_V$-invariant
  open subscheme 
$$i:\hbox{\rm Quad}_V^{sr}\hookrightarrow 
 \hbox{\rm Quad}_V$$ because, 
 over each affine open subscheme $U\hookrightarrow X$
 which trivialises $V$, it  
 corresponds to localisation by the non-zerodivisor $P_3$. Note that   
 this canonical open immersion 
 is affine and 
 schematically dominant as well.
 We next turn to semiregular bilinear forms. Recall that 
 we had defined a bilinear form $b$ to be semiregular iff its 
 induced quadratic form $q_b$ is semiregular.  
 Thus by definition, $\hbox{\rm Bil}_V^{sr}$ is the fiber product:
$$\begin{CD}
\hbox{\rm Bil}_V @>{p}>> \hbox{\rm Quad}_V\\
@A{i'}AA @AA{i}A\\
\hbox{\rm Bil}_V^{sr} @>>{p'}> \hbox{\rm Quad}_V^{sr}
\end{CD}$$
Since $p$ is a Zariski-locally-trivial principal $\hbox{\rm Alt}^2_V$-bundle,
 it is smooth and surjective (in particular faithfully flat). 
 It therefore follows that the affineness and 
 schematic dominance of $i$ imply those of 
 $i'.$ We record these facts below. 
\begin{proposition}\label{bilvsr-def}
The open immersion 
$$\hbox{\rm Bil}_V^{sr}\hookrightarrow \hbox{\rm Bil}_V$$
 is a $\hbox{\rm GL}_V$-equivariant schematically dominant affine morphism. 
 Further this open immersion behaves well under base-change (relative 
 to $X$).
\end{proposition}
\subsection{The Generalised Clifford Algebra of Bichsel-Knus}
\label{subsec2.3}
Let $R$ be a commutative ring (with 1), $I$ an invertible $R$-module and
 $V$ a projective $R$-module. Consider the Laurent-Rees algebra of 
 $I$ defined by 
$$L[I]:=R\oplus
\left(\bigoplus_{n> 0}(T^n(I) \oplus T^n(I^{-1}))\right)$$
 and define the $\mathbb Z$-gradation on the 
 tensor product of algebras $TV\otimes L[I]$ 
 by requiring elements of $V$ (resp. of $I$) 
 to be of degree one (resp. of degree two).
 Let $q:V\longrightarrow I$ be an $I$-valued quadratic form on $V.$
Following the definition of Bichsel-Knus \cite{Bichsel-Knus}, let
 $J(q, I)$ be the two-sided ideal of $TV\otimes L[I]$ generated by the set 
 $$\{(x\otimes_{TV} x)\otimes 1_{L[I]}-1_{TV}\otimes q(x)\thinspace
|\thinspace x\in V\}$$ and let the generalised Clifford algebra of $q$ be 
 defined by  $$\widetilde{C}(V,q,I):=TV\otimes L[I]/J(q,I).$$ This is an 
 $\mathbb Z$-graded algebra by definition. Let $C_n$ be the submodule of 
 elements of degree $n.$ Then $C_0$ is a subalgebra, playing the role of 
 the even Clifford algebra in the classical situation (i.e., $I=R$) and 
 $C_1$ is a $C_0$-bimodule. Bichsel and Knus baptize $C_0$ and $C_1$ 
 respectively as the {\em even Clifford algebra} and the {\em Clifford
  module} associated to the triple $(V,q,I).$
  The generalised Clifford algebra satisfies 
 an appropriate universal property which 
 ensures it behaves well functorially. Since $V$ is projective,  
 the canonical maps 
$$V\longrightarrow \widetilde{C}(V,q,I)\textrm{ and }
 L[I]\longrightarrow \widetilde{C}(V,q,I)$$ are injective. For proofs of 
 these facts, see Sec.3 of \cite{Bichsel-Knus}. 
 If $(V,q,I)$ is an $I$-valued quadratic form on the vector bundle $V$
 over a scheme $X$, with $I$ a line bundle, then the above construction 
 may be carried out to define the generalised Clifford algebra bundle 
 $\widetilde{C}(V,q,I)$ which is an $\mathbb Z$-graded algebra bundle on 
 $X.$ Its degree zero subalgebra bundle is denoted $\covqi$ and is
 called the even Clifford algebra bundle of $\vqi.$ 
\subsection{Bourbaki's Tensor Operations with Values in a Line Bundle}
\label{subsec2.4}
 Let $R$ be a commutative ring and $L[I]:=R\oplus
\left(\bigoplus_{n> 0}(T^n(I) \oplus T^n(I^{-1}))\right)$ as above. 
 We denote by $\otimes_T$ (resp. by $\otimes_L$)
 the tensor product and by $1_T$ (resp. $1_L$) 
 the unit element in the algebra $TV$ (resp. in $L[I]$).
\begin{theorem}[{\rm with the above notations}]\label{Bourbaki}\   
\begin{description}
\item[(1)]
Let $q:V\longrightarrow I$ be an $I$-valued quadratic form on $V$ and 
$f\in \hbox{\rm Hom}_R(V, I).$ Then there exists an $R$-linear
 endomorphism $t_f$ of 
 the algebra $TV\otimes L[I]$ which is unique with respect to the first three
 of the following properties it satisfies:
\begin{description}
\item[(a)] for each $\lambda\in L[I]$ we have  $t_f(1_T\otimes \lambda)
=0;$
 \item[(b)] for any $x\in V$, $y\in TV$, and $\lambda\in L[I]$ we have 
$$t_f((x\otimes_T y)\otimes \lambda)
=y\otimes(f(x)\otimes_L\lambda)-(x\otimes_T1_T)t_f(y\otimes\lambda);$$ 
\item[(c)] if $J(q, I)$ is
 the two-sided ideal of $TV\otimes L[I]$ as defined in
 \S\ref{subsec2.3} above, then we have 
$$t_f(J(q, I))\subset J(q, I);$$
\item[(d)] $t_f$ is homogeneous of degree $+1$ (except for elements which 
 it does not annihilate);
\item[(e)] if $\widetilde{C}\vqi$ is the generalised 
 Clifford algebra of $\vqi$ as defined in \S\ref{subsec2.3} above,
 then due to assertion (c) above, $t_f$ 
induces a $\mathbb Z$-graded endomorphism of degree $+1$ denoted by 
 $$d_f^q:\widetilde{C}(V,q, I)\longrightarrow \widetilde{C}(V,q, I);$$
\item[(f)]$t_f\circ t_f=0;$
\item[(g)]if $g\in \hbox{\rm Hom}_R(V,I)$, 
then $t_f\circ t_g + t_g\circ t_f=0;$
\item[(h)]if $\alpha\in \hbox{\rm End}_R(V)$ and 
 $\alpha^*f\in \hbox{\rm Hom}_R(V, I)$ is
 defined by  $x\mapsto f(\alpha(x))$ then 
 $$t_f\circ (T(\alpha)\otimes \hbox{\rm Id}_{L[I]})
=(T(\alpha)\otimes\hbox{\rm Id}_{L[I]})\circ t_{\alpha^*f};$$ 
\item[(i)]$t_f\equiv 0$ on the $R$-subalgebra of $TV\otimes L[I]$ 
generated by 
$$\hbox{\rm kernel}(f)\otimes L[I].$$ In fact, atleast when  
 $V$ is projective, the smallest $R$-subalgebra of
 $TV\otimes L[I]$ containing 
$$\hbox{\rm kernel}(f)\otimes R.1_L$$ is given by  
 $$\hbox{\rm Image}(T(\hbox{\rm kernel}(f))\otimes R.1_L$$ and $t_f$ vanishes 
 on this $R$-subalgebra. 
\end{description} 
\item[(2)] Let $q, q':V\longrightarrow I$ be two $I$-valued 
 quadratic forms whose
 difference is the quadratic form $q_b$ induced by 
 an $I$-valued 
  bilinear form 
$$b\in\hbox{\rm Bil}_R(V, I)
:=\hbox{\rm Hom}_R(V\otimes_R V, I).$$ This means that   
$$q'(x)-q(x)=q_b(x):=b(x,x)\thinspace \forall x\in V.$$
 Further, for any $x\in V$ denote by $b_x$ the element of 
 $\hbox{\rm Hom}_R(V, I)$  given by 
 $y\mapsto b(x,y).$ Then there exists an $R$-linear automorphism 
 $\Psi_b$ of 
 $TV\otimes L[I]$
 which is unique with respect to the first three of the following 
 properties it satisfies: 
\begin{description}
\item[(a)] for any $\lambda\in L[I]$ we have 
 $\Psi_b(1_T\otimes \lambda)=(1_T\otimes\lambda);$
\item[(b)] for any $x\in V$, $y\in TV$ and $\lambda\in L[I]$ we have  
$$\Psi_b((x\otimes_T y)\otimes\lambda)
=(x\otimes 1_L).\Psi_b(y\otimes\lambda)+t_{b_x}(\Psi_b(y\otimes\lambda));$$ 
\item[(c)]$\Psi_b(J(q', I))\subset J(q, I);$ 
\item[(d)]by the previous property, $\Psi_b$ induces an
 isomorphism of $\mathbb Z$-graded 
 $R$-modules $$\psi_b:\widetilde{C}(V,q', I)\cong \widetilde{C}(V,q, I);$$
 in particular, given a 
 quadratic form $q_1:V\longrightarrow I$, since there always 
 exists an $I$-valued  bilinear form $b_1$ that induces $q_1$ (i.e.,
 such that $q_1=q_{b_1}$), setting $q'=q_1$, $q=0$ and $b=b_1$ in the above 
 gives an   
 $\mathbb Z$-graded linear
 isomorphism $$\psi_{b_1}:\widetilde{C}(V,q_1, I)\cong 
 \widetilde{C}(V,0, I)=\Lambda(V)\otimes L[I];$$ 
\item[(e1)] $\Psi_b(T^{2n}V\otimes L[I])\subset
 \oplus_{(i\leq n)}(T^{2i}V\otimes L[I]);$
\item[(e2)] $\Psi_b(T^{2n+1}V\otimes L[I])\subset
 \oplus_{(\hbox{\rm\tiny odd }i\leq 2n+1)}(T^{i}V\otimes L[I]);$
\item[(e3)] $\Psi_b(T^{2n}V\otimes I^{-n})\subset
 \oplus_{(i\leq n)}(T^{2i}V\otimes I^{-i});$ 
\item[(e4)] $\Psi_b(T^{2n+1}V\otimes I^{-n})\subset
 \oplus_{(\hbox{\rm\tiny odd }i\leq 2n+1)}(T^{i}V\otimes I^{\frac{1-i}{2}});$
\item[(f)]in particular, for $x,x'\in V$, 
$$\Psi_b((x\otimes_T x')\otimes 1_L)
 =(x\otimes_T x')\otimes 1_L + 1_T\otimes b(x,x')$$ 
 so that for $$\psi_b:C_0(V,q_b, I)\cong C_0(V,0, I)
=\oplus_{n\geq 0}(\Lambda^{2n}(V)\otimes I^{-n})$$ we have 
 $$\psi_b(((x\otimes_T x')\otimes \zeta)\textrm{ mod }J(q_b, I))=(x\wedge x')
\otimes\zeta+ \zeta(b(x,x')).1$$ for any $x,x'\in V$ and 
$\zeta\in I^{-1}\equiv\hbox{\rm Hom}_R(I, R);$  
\item[(g)]if $f\in \hbox{\rm Hom}_R(V, I)$, and $t_f$ is given by (1) above,
 then $\Psi_b\circ t_f=t_f\circ \Psi_b;$
\item[(h)]for $I$-valued bilinear forms $b_i$ on $V$, 
 $$\Psi_{b_1+b_2}=\Psi_{b_1}\circ 
 \Psi_{b_2}\textrm{ and }\Psi_0=\hbox{\rm Identity on } {TV\otimes L[I]};$$
\item[(i)]for any $\alpha\in\hbox{\rm End}_R(V)$, we have 
 $$\Psi_b\circ
 (T(\alpha)\otimes \hbox{\rm Id}_{L[I]})=(T(\alpha)\otimes\hbox{\rm Id}_{L[I]})
\circ \Psi_{(b.\alpha)}$$ where $(b.\alpha)(x,x'):=b(\alpha(x),\alpha(x'))\thinspace\forall x,x'\in V;$ 
\item[(j)]by property (h), 
one has a homomorphism of groups $$(\hbox{\rm Bil}_R(V, I), +)
\longrightarrow (\hbox{\rm Aut}_R(TV\otimes L[I]), \circ):b\mapsto \Psi_b;$$
 the associative
 unital monoid $(\hbox{\rm End}_R(V), \circ)$ acts on $\hbox{\rm Bil}_R(V, I)$ 
 on the right by $b'\leadsto b'.\alpha$ and acts on the 
 left (resp. on the right) of $\hbox{\rm End}_R(TV\otimes L[I])$ by
 $$\alpha.\Phi:=(T(\alpha)\otimes\hbox{\rm
 Id}_{L[I]})\circ\Phi\textrm{ (resp. by }\Phi.\alpha:=\Phi\circ
 (T(\alpha)\otimes\hbox{\rm Id}_{L[I]})\textrm{ )},$$
 and the  homomorphism
 $b\mapsto \Psi_b$ satisfies 
$$\alpha.\Psi_{(b.\alpha)}=\Psi_b.\alpha;$$ 
 the group $\hbox{\rm Aut}_R(V)=\hbox{\rm GL}_R(V)$ acts on
  the left of $\hbox{\rm Bil}_R(V, I)$ by 
 $$g.b: (x,x')\mapsto b(g^{-1}(x), g^{-1}(x'))$$ and on the left of 
 $\hbox{\rm Aut}_R(TV\otimes L[I])$ by conjugation via the natural 
 group homomorphism $$\hbox{\rm GL}_R(V)\longrightarrow
 \hbox{\rm Aut}_R(TV\otimes L[I]): 
 g\mapsto T(g)\otimes\hbox{\rm Id}_{L[I]}$$
 namely we have  
$$g.\Phi:=(T(g)\otimes\hbox{\rm Id}_{L[I]})\circ \Phi\circ 
 (T(g^{-1})\otimes\hbox{\rm Id}_{L[I]}),$$ and the 
 homomorphism $b\mapsto \Psi_b$ is $\hbox{\rm GL}_R(V)$-equivariant: 
 $\Psi_{g.b}=g.\Psi_b.$ 
\end{description} 
\item[(3)] For a commutative $R$-algebra $S$ (with 1), let 
 $$(q\otimes_R S), (q'\otimes_R S):
 (V\otimes_R S=:V_S)\longrightarrow (I\otimes_R S=:I_S)$$ 
 be the $I_S$-valued quadratic forms induced from 
 the quadratic forms $q,q'$ of (2) above
 and $$(b\otimes_R S)\in\hbox{\rm Bil}_S(V_S, I_S)$$
 the $I_S$-valued $S$-bilinear form induced 
 from the bilinear form $b$ of (2) above. 
 Then as a result  of the uniqueness properties (2a)--(2c) satisfied 
 by $\Psi_b$ and $\Psi_{(b\otimes_R S)}$, the $S$-linear automorphisms
 $(\Psi_b\otimes_R S)$ and $\Psi_{(b\otimes_R S)}$ 
may be canonically  identified. In particular,
 the $\mathbb Z$-graded $S$-linear isomorphism 
 $$(\psi_{b}\otimes_R S):\widetilde{C}(V_S,(q'\otimes_R\thinspace S), I_S)\cong\widetilde{C}(V_S,(q\otimes_R S), I_S)$$ 
 induced from $\psi_b$ of
 (2d) above may be canonically identified with $\psi_{(b\otimes_R S)}.$ 
\end{description}
\end{theorem}
\begin{remark}\rm 
As mentioned in \cite{Bichsel-Knus}, F. van Oystaeyen has observed that 
 $L[I]$ is a faithfully-flat splitting for $I$, and the generalised 
 Clifford algebra $\widetilde{C}(V,q,I)$ is nothing but the ``classical'' 
 Clifford algebra of the triple 
$$(V\otimes_R L[I], q\otimes_R\thinspace L[I], 
 I\otimes_R L[I]\equiv L[I])$$
 over $L[I].$  In the same vein, $I$-valued forms 
 (both multilinear and quadratic) on an $R$-module $V$ can be treated as 
 the usual ($L[I]$-valued) forms on $V\otimes_R L[I].$ 
 With this in mind, the proof of Theorem \ref{Bourbaki} follows from the 
 usual Bourbaki tensor operations with respect to $V\otimes_R L[I]$ on 
 $L[I].$  (See \S 9, Chap.9, \cite{Bourbaki} or 
 para.1.7, Chap.IV, \cite{Knus} 
 for the ``classical'' Bourbaki operations). 
 However one needs to remember that the $\mathbb Z$-gradation on 
 $TV\otimes_R L[I]$ as defined above is different from the usual 
 ${\mathbb Z}_{\geq 0}$-gradation on $T(V\otimes_R L[I]).$  
\end{remark}
\subsection{Tensoring by Bilinear- and 
 Twisted Discriminant Bundles}
\label{subsec2.5}
Let $V,M$ be vector bundles on a scheme $X$ and let $I,J$ be line bundles on 
 $X.$ Let $q$ be a quadratic form on $V$ with values in $I$ and let 
 $b$ be a symmetric bilinear form on $M$ with values in $J.$ By abuse
 of notation, we also use $b$ to denote the  corresponding 
$J$-valued linear form on $M\otimes M.$ 
\begin{proposition}[{\rm with the above notations}]
\label{tensoring-with-sym-bil-module}\
\begin{description}
\item[(1)] The tensor product of $\vqi$ with $(M,b, J)$ gives a unique 
 quadratic bundle 
$$(V\otimes M, 
 q\otimes b, I\otimes J).$$ The quadratic form $q\otimes b$ 
 on $V\otimes M$ is given on 
 sections by $$v\otimes m\mapsto q(v)\otimes b(m\otimes m)$$ 
 and has associated bilinear 
 form $b_{q\otimes b}=b_q\otimes b.$
\item[(2)] When $M$ is a line bundle, $(M,b,J)$ is regular (=nonsingular)
  iff $(M, q_b, J)$  is 
 semiregular iff 
$$b:M\otimes M\cong J$$ is an 
 isomorphism.
\item[(3)] Let $V$ be of odd rank and $M$ a
 line bundle such that $(M,b,J)$ is regular.
 Then $\vqi$ is semiregular iff
 $$\vqi\otimes(M,b,J)=(V\otimes M, q\otimes b, I\otimes J)$$ is semiregular. 
\end{description}
\end{proposition} 
\begin{definition}  A triple $(L,h,J)$, consisting of a linear
  isomorphism $h:L\otimes L\cong J$, with $L$ and $J$ both line
 bundles, is called a {\em twisted discriminant bundle} on $X.$
\end{definition}
A twisted discriminant bundle $(L,h,J)$ specifies $L$ as a 
square root of the line bundle $J$ via $h$.
 The terminology is
motivated by the following: when $J$ is the trivial
line bundle, such a datum is referred to as a {\em discriminant bundle}
in \S 3, Chap.III, of Knus' book \cite{Knus}.
 By part (2) of the preceding
Proposition, $h$ is a regular bilinear form on $L$
(necessarily symmetric) with values in $J$, so that we may speak of an
isometry between two twisted discriminant bundles $(L,h,J)$ and 
$(L',h',J')$: it is a pair $(\zeta, \eta)$ consisting of linear
isomorphisms $\zeta:L\cong L'$ and $\eta:J\cong J'$ such that 
$\eta h=h'(\zeta\otimes\zeta).$ 
\begin{lemma}\label{t-disc-defn}
 On the set $\textrm{T-Disc}(X)$ of isometry
classes of twisted discriminant bundles on $X$, we have a natural 
group structure  induced
by the tensor product. $\textrm{T-Disc}(X)$ is functorial in $X.$
If we consider only isometry classes of discriminant bundles, i.e., of 
triples $(L, h, {\cal O}_X)$, then we obtain a subgroup 
$\textrm{Disc}(X)\subset\textrm{T-Disc}(X)$ which is of exponent 2. 
\end{lemma}
\begin{proposition}\label{tensoring-with-disc-module}
Let $V$ and $V'$ be vector bundles of the same rank on the scheme $X$, 
 $(L,h, J)$ a twisted discriminant line bundle on $X$ and 
 $$\alpha:V'\cong V\otimes L$$ an isomorphism of bundles. 
\begin{description}
\item[(1)] Over any open subset
 $U\hookrightarrow X$, given a bilinear form 
 $$b'\in\Gamma\left(U, \hbox{\rm Bil}_{(V',I)}\right),$$
 we can define a 
 bilinear form 
$$b\in\Gamma\left(U, \hbox{\rm Bil}_{(V, I\otimes J^{-1})} \right)$$ 
 using $\alpha$ and $h$ as follows: we let 
$$b:=(b'\otimes J^{-1})\circ{(\zeta_{(\alpha, h)})}^{-1}$$
 where $$\zeta_{(\alpha,h)}:V'\otimes V'\otimes J^{-1}\cong
 V\otimes V$$ is the linear isomorphism given
 by the composition of the following natural morphisms:
\begin{center}\begin{tabular}{c}
$V'\otimes V'\otimes J^{-1}
\stackrel{\alpha\otimes\alpha\otimes\hbox{\rm\tiny Id}
\thinspace(\cong)}{\longrightarrow}
V\otimes L\otimes V\otimes L\otimes J^{-1}
\stackrel{\hbox{\rm\tiny SWAP}(2,3)\thinspace(\equiv)}
{\longrightarrow}$\\
$V\otimes V\otimes L^2\otimes J^{-1}
\stackrel{\hbox{\rm\tiny Id}\otimes h\otimes\hbox{\rm\tiny Id}
\thinspace(\cong)}{\longrightarrow}$\\
$\stackrel{\hbox{\rm\tiny Id}\otimes h\otimes\hbox{\rm\tiny Id}
\thinspace(\cong)}{\longrightarrow}V\otimes V\otimes J\otimes J^{-1}
\stackrel{\hbox{\rm\tiny CANON}\thinspace(\equiv)}{\longrightarrow}
V\otimes V.$ \end{tabular}\end{center}
Then the association $b'\mapsto b$ induces linear isomorphisms 
 shown by vertical upward arrows in the following diagram of associated 
 locally-free sheaves (with 
 exact rows) making it commutative: 
$$\hspace*{-10pt}\begin{CD}
0 @>>> \hbox{\rm Alt}^2_{({\cal V},{\cal I}\otimes {\cal J}^{-1})}
 @>>> \hbox{\rm Bil}_{({\cal V}, {\cal I}\otimes {\cal J}^{-1})}
 @>>> \hbox{\rm Quad}_{({\cal V}, {\cal I}\otimes {\cal J}^{-1})} @>>> 0\\
 & & @A{\cong}AA  @A{\cong}AA  @A{\cong}AA & & \\
0 @>>> \hbox{\rm Alt}^2_{({\cal V}', {\cal I})} @>>> \hbox{\rm Bil}_{({\cal V}',{\cal I})} @>>>
 \hbox{\rm Quad}_{({\cal V}', {\cal I})} @>>> 0.
\end{CD}$$ 
Therefore one also has the following commutative diagram of vector bundle 
 morphisms with the vertical upward arrows being isomorphisms: 
$$\begin{CD}
 \hbox{\rm Alt}^2_{(V, I\otimes J^{-1})} 
@>{\hbox{\rm closed}}>{\hbox{\rm immersion}}> 
\hbox{\rm Bil}_{(V, I\otimes J^{-1})}
 @>{\hbox{\rm locally}}>{\hbox{\rm trivial}}>
 \hbox{\rm Quad}_{(V, I\otimes J^{-1})} \\
 @A{\cong}AA  @A{\cong}AA @A{\cong}AA \\
\hbox{\rm Alt}^2_{(V', I)} 
@>{\hbox{\rm closed}}>{\hbox{\rm immersion }}>
 \hbox{\rm Bil}_{(V', I)}
 @>{\hbox{\rm locally}}>{\hbox{\rm trivial}}>
 \hbox{\rm Quad}_{(V', I)} 
\end{CD}$$ 
\item[(2)]
 Let $b'\in
\Gamma\left(X, \hbox{\rm Bil}_{(V', I)}\right)$ be a global 
 bilinear form and let it induce 
$$b\in\Gamma\left(X, \hbox{\rm Bil}_{(V, I\otimes J^{-1})}\right)$$
 via $\alpha$ and $h$  as defined in (1) above.
 Let 
$$\Psi_{b'}\in\hbox{\rm Aut}_{{\cal O}_X}(TV'\otimes L[I])\textrm{  
(resp. }\Psi_{b}\in
\hbox{\rm Aut}_{{\cal  O}_X}(TV\otimes L[I\otimes J^{-1}]))$$ 
be the $\mathbb Z$-graded linear isomorphism   
 induced by 
 $b'$ (resp. by $b$) defined locally (and hence globally) 
 as in (2), Theorem \ref{Bourbaki} above.
 Let $$Z_{(\alpha, h)}:\oplus_{n\geq 0}(T^{2n}_{{\cal  O}_X}(V')\otimes I^{-n})
\cong\oplus_{n\geq 0}(T^{2n}_{{\cal  O}_X}(V)\otimes I^{-n}\otimes J^n)$$
 be the ${\cal  O}_X$-algebra isomorphism induced via the isomorphism
 $\zeta_{(\alpha,h)}$ defined in (1) above. Then, taking into account 
 (2e), Theorem \ref{Bourbaki},  
 the following diagram commutes:
$$\begin{CD}
\oplus_{n\geq 0}(T^{2n}_{{\cal  O}_X}(V')\otimes I^{-n}) 
@>{Z_{(\alpha, h)}}>{\cong}>
\oplus_{n\geq 0}(T^{2n}_{{\cal  O}_X}(V)\otimes I^{-n}\otimes J^n)\\
@V{\Psi_{b'}}V{\cong}V  @V{\cong}V{\Psi_b}V \\
\oplus_{n\geq 0}(T^{2n}_{{\cal  O}_X}(V')\otimes I^{-n}) 
@>{\cong}>{Z_{(\alpha, h)}}>
\oplus_{n\geq 0}(T^{2n}_{{\cal  O}_X}(V)\otimes I^{-n}\otimes J^n)
\end{CD}$$ 
thereby inducing by (2d), Theorem \ref{Bourbaki} 
the following commutative diagram of ${\cal  O}_X$-linear 
 isomorphisms
$$\begin{CD}
C_0(V', q_{b'}, I) @>{\hbox{\rm via }Z_{(\alpha, h)}}>{\cong}>C_0(V, q_b, I\otimes J^{-1})\\
@V{\psi_{b'}}V{\cong}V @V{\cong}V{\psi_b}V \\
\oplus_{n\geq 0}(\Lambda^{2n}_{{\cal  O}_X}(V')\otimes I^{-n})
 @>{\cong}>{\hbox{\rm via }Z_{(\alpha, h)}}>
\oplus_{n\geq 0}(\Lambda^{2n}_{{\cal  O}_X}(V)\otimes I^{-n}\otimes J^n)
\end{CD}$$ 
\item[(3)]
 Let $b$ and $b'$ be as in (2) above. Then $\alpha:V'\cong V\otimes L$ 
 induces an isometry 
 of bilinear form bundles 
$$\alpha:(V',b', I)\cong(V,b, I\otimes J^{-1})
\otimes(L,h, J)$$ and also an isometry of the induced 
 quadratic bundles 
$$\alpha: (V',q_{b'}, I)\cong 
 (V, q_b, I\otimes J^{-1})\otimes (L,h, J).$$
 Moreover, if we are just given a global $I\otimes J^{-1}$-valued quadratic 
 form $q$ on $V$ (resp. an $I$-valued $q'$ on $V'$),
 then we may define the global $I$-valued quadratic 
 form $q'$ on $V'$ (resp. $I\otimes J^{-1}$-valued $q$ on $V$)
 via 
$$q':=(q\otimes h)\circ \alpha\textrm{  (resp. via }
q:=(q'\circ\alpha^{-1})\otimes{(h^\vee)}^{-1})$$ 
 and again we have an isometry of quadratic bundles  
$$\alpha:(V',q', I)\cong (V,q, I\otimes J^{-1})\otimes(L,h, J).$$
\end{description} 
\end{proposition}  
\begin{proposition}\label{simil-induces-iso-of-even-cliff}
Let $g:\vqi\thinspace{\cong}_{l}\thinspace\vprqpri$
 be an $I$-similarity with multiplier
 $l\in\Gamma(X, {\cal  O}_X^*).$ 
\begin{description}
\item[(1)] There exists a unique isomorphism of ${\cal  O}_X$-algebra bundles
$$C_0(g,l, I):C_0(V,q, I)\cong C_0(V',q', I)$$ such that  
for sections $v, v'$ of $V$ and 
 $s$ of $I^{-1}$ we have  $$C_0(g,l, I)(v.v'.s)=g(v).g(v').l^{-1}s.$$ 
\item[(2)] There exists a unique vector bundle isomorphism 
$$C_1(g,l, I):C_1(V,q, I)\cong C_1(V', q', I)$$ such that 
the following hold 
for any section $v$ of $V$ and any section $c$ of $C_0(V,q)$: 
\begin{description}
\item[(a)] 
 $C_1(g,l, I)(v.c)=g(v).C_0(g,l, I)(c)$ and 
\item[(b)] $C_1(g,l, I)(c.v)=C_0(g,l,I)(c).g(v).$ 
\end{description} 
 Thus $C_1(g,l, I)$ is $C_0(g,l, I)$-semilinear.  
\item[(3)]
 If $g_1:\vprqpri\thinspace{\cong}_{l_1}\thinspace(V'',q'',I)$ is another 
  similarity with multiplier $l_1$, then the composition 
 $$g_1\circ g:\vqi\thinspace{\cong}_{ll_1}\thinspace(V'',q'', I)$$ is 
 also a similarity with multiplier given by the product of 
 the multipliers. Further 
$$C_i(g_1\circ g,ll_1,I)=
 C_i(g_1,l_1, I)\circ C_i(g,l, I)\textrm{ for }i=0,1.$$ 
\end{description}
\end{proposition}
\noindent 
 A local computation 
 shows that tensoring by a twisted 
 discriminant bundle is locally the same as applying 
 a similarity. In this case also one gets a global
  isomorphism of even Clifford 
 algebras: 
\begin{proposition}\label{isom-covq-covlqh}
Let $\vqi$ be a quadratic bundle on a scheme $X$ and $(L,h, J)$ be a twisted 
 discriminant bundle. There exists a unique isomorphism of algebra bundles 
 $$\gamma_{(L,h,J)}:C_0\left(\vqi\otimes(L,h,J)\right)\cong C_0(V,q,I)$$ 
 such that for any  sections $v,v'$ of $V$, $\lambda,\lambda'$ of $L$, 
 $s$ of $I^{-1}\equiv I^{\vee}$ and $t$ of $J^{-1}\equiv J^{\vee}$ 
we have  
$$\gamma_{(L,h,J)}\left((v\otimes\lambda).
(v'\otimes \lambda').(s\otimes t)\right)=
t(h(\lambda\otimes \lambda'))v.v'.s.$$ 
\end{proposition}
\subsection{The Theorem of Max-Albert Knus}
\label{newsubsec2.6}
For any scheme $X,$ denote by ${\cal Q}_3^{sr}(X;{\cal O}_X)$ the set of 
isomorphism classes of {\em semiregular} ternary quadratic bundles
with values in the trivial line bundle 
on $X.$ Here isomorphism stands for isometry as defined
 in \S\ref{newsubsec2.2}. Recall the group $\textrm{Disc}(X)$ of 
 discriminant bundles on $X$ (Lemma~\ref{t-disc-defn}). By 
assertions (2) and (3) of 
 Prop.~\ref{tensoring-with-sym-bil-module}, this group
  acts on ${\cal Q}_3^{sr}(X;{\cal O}_X).$

 For a ternary quadratic form $q:V\rightarrow {\cal O}_X,$ recall from
 \S\ref{subsec2.3} the Bichsel-Knus even-Clifford algebra
 $C_0(V,q,{\cal O}_X)$, namely 
 the degree zero subalgebra of the 
 generalised Clifford algebra $\widetilde{C}(V,q,{\cal O}_X).$ 
Since the values of the form are in the trivial line bundle, 
 this even-Clifford algebra is the same as the 
  the classically-defined even-Clifford algebra.

Now the even-Clifford algebra of a semiregular form is Azumaya (see
for instance, Prop.3.2.4, \S 3, Chap.IV, \cite{Knus}). 
So if we let ${\cal AZU}_4(X)$ to denote the set of
  algebra-isomorphism classes of rank 4 Azumaya bundles over $X.$
 then by Prop.~\ref{isom-covq-covlqh}, 
 the association $(V,q,{\cal O}_X)\leadsto C_0(V,q,{\cal O}_X)$ 
 induces a map 
 $${\textrm{\sf Witt-Invariant}_{(X;{\cal O}_X)}^{sr}}:{\cal Q}_3^{sr}(X; {\cal O}_X)/\textrm{Disc}(X)\rightarrow{\cal AZU}_4(X),$$ 
where the left side represents the set of orbits. 
\begin{theorem}[Max-Albert Knus, \S~3, Chap.V,
 \cite{Knus}]\label{Max-Knus-Theorem} 
For any scheme $X$, \newline the map 
$\textrm{\sf Witt-Invariant}_{(X;{\cal O}_X)}^{sr}$
defined above is a bijection.
\end{theorem}
 By Prop.3.2.2, \S 3, Chap.III, \cite{Knus}, 
\label{coho-interpretation}
 the group $\textrm{Disc}(X)$
 is naturally isomorphic to the cohomology (abelian group)  
 ${\hbox{\rm\v H}}^1_{\hbox{\rm\tiny fppf}}(X, \mu_2).$
Further, 
  by Lemma 3.2.1, \S 3, Chap.IV, \cite{Knus}, the cohomology   
  ${\hbox{\rm\v H}}^1_{\hbox{\rm\tiny fppf}}(X, \hbox{\sf O}_3)$
 classifies the set of isomorphism classes of semiregular rank 3 
 quadratic bundles with values in the trivial line bundle, so it is 
the same as ${\cal Q}_3(X; {\cal O}_X).$  
 On the other hand, 
 the set of isomorphism classes of rank 4 Azumaya algebras 
 on $X$ may be interpreted as the cohomology 
 ${\hbox{\rm\v H}}^1_{\hbox{\rm\tiny \'etale}}(X, \hbox{\sf PGL}_2)$
 (see page 145, \S 5, Chap.III, \cite{Knus}). 
 Thus the bijection of Theorem~\ref{Max-Knus-Theorem} can be thought
 of as a statement in cohomology:  
 $${\hbox{\rm\v H}}^1_{\hbox{\rm\tiny fppf}}(X, \hbox{\sf O}_3)/{\hbox{\rm\v H}}^1_{\hbox{\rm\tiny fppf}}(X, \mu_2)\cong
 {\hbox{\rm\v H}}^1_{\hbox{\rm\tiny \'etale}}(X, \hbox{\sf PGL}_2).$$ 
\subsection{Results on the notion of Schematic Image}
\label{subsec2.6}
\begin{definition}[Defs.6.10.1-2, Chap.I, EGA I \cite{ega1}]\rm\label{sch-img-def}
Let $f:X\longrightarrow Y$ be a morphism of schemes. If there exists 
a smallest closed subscheme $Y'\hookrightarrow Y$ such that the inverse 
 image scheme 
$$f^{-1}(Y'):=Y'\times_Y {(}_fX)$$ is equal to $X$, 
one calls $Y'$ the {\em schematic image} of $f$ (or of $X$ in $Y$ under $f$).
If $X$ were a subscheme of $Y$ and $f$ the canonical immersion, and 
 if $f$ has a schematic image $Y'$, then $Y'$ is called the {\em schematic 
 limit} or the {\em limiting scheme}
 of the subscheme $X\stackrel{f}{\hookrightarrow}Y.$ 
\end{definition}
\begin{proposition}[Prop.6.10.5, Chap.I, EGA I]\label{sch-img-existence}
The schematic image $Y'$ of $X$ by a morphism $f:X\longrightarrow Y$ exists 
 in the following two cases:
\begin{description}
\item[(1)] $f_*({\cal O}_X)$ is a quasi-coherent ${\cal O}_Y$-module, which 
 is for example the case when $f$ is quasi-compact and
 quasi-separated;
\item[(2)]  
$X$ is reduced. 
\end{description}
\end{proposition}
\begin{proposition}\label{sch-img-props}
In each of the following statements whenever a schematic image is 
 mentioned, we assume that one of the two hypotheses  
 of the above Prop.~\ref{sch-img-existence}
 is true so that the schematic image does exist. 
\begin{description}
\item[1.] Let $Y'$ be the schematic image of $X$ under a morphism 
 $f:X\longrightarrow Y$ and let $f$ factor as 
$$X\stackrel{g}{\longrightarrow}Y'\stackrel{j}{\hookrightarrow}Y.$$
Then $Y'$ is topologically the closure of $f(X)$ in $Y$, the morphism 
 $g$ is schematically dominant i.e., 
$$g^{\#}:{\cal O}_{Y'}\longrightarrow 
 g_*({\cal O}_X)$$ is injective and the 
 schematic image of $X$ in $Y'$ (under $g$) is $Y'$ itself. If $X$ is 
 reduced (respectively integral) then the same is true of $Y'.$ 
\item[2.] The schematic image of a closed subscheme under its canonical 
 closed immersion is itself.
\item[3.] (Transitivity of Schematic Image.) Let there be given morphisms 
$$X\stackrel{f}{\longrightarrow}Y\stackrel{g}{\longrightarrow}Z,$$ 
such that the schematic image $Y'$ of $X$ under $f$ exists, and further 
 such that if $g'$ is the restriction of $g$ to $Y'$, the schematic image $Z'$ 
 of $Y'$ by $g'$ exists. Then the schematic image of $X$ under $g\circ f$ 
 exists and equals $Z'.$   
\item[4.] Let $f:X\longrightarrow Y$ be a morphism which factors through 
 a closed subscheme $
Y_1$ of $Y$ by a morphism $f_1:X\longrightarrow Y_1.$
Then the scheme-theoretic image $Y'$ of $X$ in $Y$ is the same as the 
 scheme-theoretic image $Y'_1$ of $X$ in $Y_1$ considered canonically 
 as closed subscheme of $Y.$ 
\item[5.] If $f:X\longrightarrow Y$ has a schematic image $Y'$ then 
 $f$ is schematically dominant 
 iff $Y'=Y.$
\item[6.] The formation of schematic image commutes with flat base 
 change: if $f:X\longrightarrow Y$ is a morphism of $S$-schemes 
 which has a schematic image $Y'$ then
 for a flat morphism $S'\longrightarrow S$,
 one has that the induced $S'$-morphism
 $$f\times_S S':X\times_S S'\longrightarrow Y\times_S S'$$
 has a schematic image and it 
 may be canonically identified with $Y'\times_S{S'}.$ In particular this means 
 that the formation of schematic image is local over the base. 
\end{description}
\end{proposition}
Assertions (1) and (3) are respectively Prop.6.10.5 and Prop.6.10.3 
 in EGA I. The defining property of schematic image gives (2),
 while (4) can be deduced from the first three. As for (5), 
 from (1) it follows that $Y'=Y$ implies $f=g$ is schematically 
 dominant. For the other way around, one uses 
 the following characterisation of a schematically dominant morphism 
 in Prop.5.4.1 of EGA I: if $f:X\longrightarrow Y$ is a morphism of schemes, 
then $f$ is schematically dominant iff for every open subscheme $U$ of $Y$ 
 and every closed subscheme $Y_1$ of $U$ such that there exists a 
 factorisation 
$$f^{-1}(U)\stackrel{g_1}{\longrightarrow}Y_1
\stackrel{j_1}{\hookrightarrow}U,$$ 
of the restriction $f^{-1}(U)\longrightarrow U$ of $f$ (where $j_1$ is 
 the canonical closed immersion), one has $Y_1=U.$ Given $f$ is 
 schematically dominant, one just has to take $U=Y$, $Y_1=Y'$ and $g_1=g.$ 
Assertion (6) follows from statement (ii) (a) of Theorem 11.10.5 of EGA IV
 \cite{ega4}.
\subsection{Specialisations of Rank 4 Azumaya Algebras: Recap}
\label{subsec2.7}
Until further notice we assume that $W$ is a vector bundle of fixed positive 
 rank on the scheme $X$ with associated coherent locally-free sheaf 
 $\cal W.$ 
Given any $X$-scheme $T$, 
  by a  $T$-algebra structure  
 on ${{W}}_T:={{W}}\times_X T$
 (also referred to as $T$-algebra bundle), we mean  
 a morphism 
 $${{W}}_T\times_T{{W}}_T\longrightarrow{{W}}_T$$ 
 of vector bundles on $T$ arising from a morphism of the associated 
 locally-free sheaves. Given such a $T$-algebra structure and  $T'
\longrightarrow T$ an $X$-morphism, it is clear that one  
 gets by pullback (i.e., by base-change) a canonical $T'$-algebra 
 structure on ${{W}}_{T'}.$ 
Thus one has a contravariant 
 ``functor of algebra structures on ${W}$'' from  
 $\{X-Schemes\}$ to $\{Sets\}$ 
 denoted $\falgW$ with 
$$\falgw(T)=\{T-\textrm{algebra structures on }
 {{W}}_T\}=\hbox{Hom}_{{\cal O}_T}\left({\cal W}_T\otimes{\cal W}_T,
 {\cal W}_T\right).$$
It follows from 
 Prop.9.6.1, Chap.I of EGA I \cite{ega1} that the functor
 $\falgW$ is represented by the $X$-scheme
 $$\algW:=\spec{\sym{X}{{\cal W}}}.$$ 
 Hence
 $\algW$ is affine (therefore separated),
  of finite presentation over $X$  and in fact smooth 
 of relative dimension $\hbox{rank}_X({{W}})^3$. If $X'\longrightarrow 
 X$ is an extension of base, then the construction $\algW$ base-changes 
 well i.e., one may canonically identify 
 $\algW\times_X X'$ with 
 $\hbox{Alg}_{{W}'}$ where ${{W}'}={{W}}\times_X X'$ 
 (cf. Prop.9.4.11, Chap.I, EGA I \cite{ega1}).

The general linear groupscheme associated to ${W}$ viz $\gl{{{W}}}$ 
 naturally acts on $\algW$ on the left, so that for each $X$-scheme $T$, 
 $\algW(T)$ mod $\gl{{{W}}}(T)$ is the set of isomorphism classes of 
 $T$-algebra structures on ${{W}}\times_X T$.

 We remark that 
 an algebra structure may fail to be associative and may fail 
 to have a (two-sided) identity element for multiplication. However, 
 a multiplicative identity for an 
 associative algebra structure must be a nowhere vanishing 
 section (Lemma 2.3, and (2)$\Rightarrow$(4) of Lemma 2.4, 
 Part A, \cite{tevb-paper1}).

 Let $w\in\Gamma(X, {\cal W})$ be a nowhere vanishing 
 section. 
For any  $X$-scheme $T$, let 
$\faaWw(T)$ denote the subset of $\algW(T)$ consisting of 
 associative algebra structures with multiplicative identity 
 the nowhere vanishing section $w_T$ over $T$ induced from $w$.
 Thus we obtain 
 a contravariant subfunctor 
$\faaWw$ of $\algW.$

Let $\stabw(T)\subset \gl{{{W}}}(T)$
 denote the stabiliser subgroup 
 of $w_T$, so that one gets a subfunctor in
 subgroups $\stabw\subset\gl{{{W}}}.$
 It is in fact represented by a closed subgroupscheme (also denoted
 by) $\stabw$ and further behaves well under base change relative 
 to $X$ i.e., $\stabw\times_X T$ can be canonically identified with
 $\hbox{\rm Stab}_{w_T}$ for any $X$-scheme $T.$
 These follow from para 9.6.6 of Chap.I, EGA I \cite{ega1}. 

 It is clear that the natural action of $\gl{{{W}}}$ on $\algW$ 
 induces one of $\stabw$ on $\faaWw.$
 It is easy to check
 (p.489, Part A, \cite{tevb-paper1})
 that the functor $\faaWw$ is a sheaf in the 
 big Zariski site over $X$ and further that this functor is 
 represented by a natural closed subscheme  
$$\aaWw\hookrightarrow
 \algW$$ which is $\stabw$-invariant. 
 Further the construction $\aaWw$ behaves well with respect to 
 base-change (relative to $X$).
Consider the  subfunctor 
 $$\fazuWw\hookrightarrow\aaWw$$ 
 corresponding to Azumaya algebras.
\begin{theorem}[Theorem 3.4, Part A, \cite{tevb-paper1}]
\label{repbility-smoothness-azu}\ 
\begin{description}
\item[(1)] $\fazuWw$
 is represented by a $\stabw$-stable
 open subscheme $$\azuWw\hookrightarrow\aaWw$$ and the canonical open 
 immersion is an affine morphism.
\item[(2)] $\azuWw$ is affine (hence separated) 
 and of finite presentation over $X$, and 
 $\azuWw$ behaves well with respect to base-change (relative to $X$). 
\item[(3)] Further, $\azuWw$ 
 is smooth of relative dimension $(m^2-1)^2$
 and geometrically irreducible over $X$, 
 where $m^2:=\hbox{rank}_X({\cal W})$. 
\end{description}
\end{theorem} 
\begin{theorem}[Theorem 3.8, Part A, \cite{tevb-paper1}]
\label{smoothness-spazu}\ 
\begin{description}
\item[(1)] The open immersion $\azuWw\hookrightarrow\aaWw$
 has a schematic image denoted 
 $$\spazuWw\hookrightarrow\aaWw$$ which is affine (hence separated) 
  and of finite type over $X$ and is naturally a $\stabw$-stable
 closed subscheme of $\aaWw$,
 the action 
 extending the natural one on the open subscheme $\azuWw$.
 \item[(2)] When the rank of ${W}$ over 
 $X$ is 4, $\spazuWw$ is locally (over $X$) isomorphic to relative 
 9-dimensional affine space; in fact over every open affine subscheme $U$ of 
 $X$ where ${W}$ becomes trivial and $w$ becomes part of a global 
 basis we have  $$\spazuWw|_U\cong{\mathbb A}^9_U.$$
 For the explicit isomorphism, see Theorem~\ref{Thetaisom},
 page~\pageref{Thetaisom}. 
 Thus $\spazuWw$ is smooth of relative 
 dimension 9 and geometrically irreducible over $X.$ In particular, 
 it is of finite presentation over $X.$  
\item[(3)] When $\hbox{rank}_X({{W}})=4$,
 the construction 
 $\spazuWw\longrightarrow X$ behaves well with respect to base change
 (relative to $X$). 
\end{description}
\end{theorem}


\section{Statements of the Main Results}   
The theory of semiregular/regular
 quadratic forms of low rank is well-known, as in 
 Chap.V of Knus' book \cite{Knus}. 
 This theory is 
 satisfactory for such forms, since it classifies them in terms
 of various invariants and also gives information on the groups of 
 generalised similarities, isometries and special isometries of 
 such forms. It is natural to look for a corresponding theory for 
  limiting or degenerate forms. We formulate in the
 following such a theory for the case of ternary quadratic forms.

\label{sec1}\subsection{A Limiting Version of a
  Theorem in Cohomology}
\label{subsec1.1}
\noindent For any scheme $X$, denote by ${\cal Q}_3(X)$ (respectively,
 by ${\cal Q}_3^{sr}(X)$) the set of 
isomorphism classes of line-bundle-valued ternary quadratic bundles
(respectively line-bundle-valued semiregular ternary quadratic
 bundles) 
on $X.$ Here isomorphism stands for isometry as defined
 in \S\ref{newsubsec2.2}. Consider the group 
 $\textrm{T-Disc(X)}$ of twisted discriminant bundles on $X$
 (cf.~Lemma~\ref{t-disc-defn}). By 
 Prop.~\ref{tensoring-with-sym-bil-module}, 
 this group 
 acts on ${\cal Q}_3(X)$ and on
 the subset ${\cal Q}_3^{sr}(X).$

 For a line-bundle-valued quadratic form $\vqi$, recall from
 \S\ref{subsec2.3} the Bichsel-Knus even-Clifford algebra $C_0(V,q,I),$
 which is the degree zero subalgebra of the 
 generalised Clifford algebra $\widetilde{C}(V,q,I)$ and reduces to  
  the usual even Clifford algebra for a quadratic form with 
 values in the structure sheaf.

 Let ${\cal SPAZU}_4(X)$ (respectively, ${\cal AZU}_4(X)$) denote the set of
 isomorphism classes of 
associative unital algebra structures on vector bundles 
 of rank 4 over $X$ that are Zariski-locally isomorphic to 
 even-Clifford algebras of rank 3 quadratic 
 bundles (respectively, that are Azumaya).
 Recall that the Theorem~\ref{Max-Knus-Theorem}
(page~\pageref{Max-Knus-Theorem}) of Max-Albert Knus gives a bijection 
$$\textrm{\sf Witt-Invariant}_{(X;{\cal O}_X)}^{sr}:{\cal Q}_3^{sr}(X; {\cal
  O}_X)/\textrm{Disc(X)}\stackrel{\cong}{\rightarrow}{\cal AZU}_4(X).$$ 
It follows that 
 ${\cal AZU}_4(X)\subset
 {\cal SPAZU}_4(X).$
By Prop.~\ref{isom-covq-covlqh}, 
 the association $\vqi\leadsto\covqi$  induces a map 
 $$\textrm{\sf Witt-Invariant}_X:{\cal
   Q}_3(X)/\textrm{T-Disc}(X)\rightarrow 
{\cal SPAZU}_4(X),$$ 
where the left side represents the set of orbits. 
Since the even-Clifford algebra of a semiregular quadratic module is 
Azumaya, it follows that 
the above map restricts to a map: 
$$\textrm{\sf Witt-Invariant}_X^{sr}:{\cal
  Q}_3^{sr}(X)/\textrm{T-Disc}(X)\rightarrow {\cal AZU}_4(X).$$ 
\begin{theorem}\label{bijectivity}
For any scheme $X$, both the map $\textrm{\sf Witt-Invariant}_X$
 as well as  its restriction $\textrm{\sf Witt-Invariant}_X^{sr}$ are
 bijections. 
\end{theorem}
Thus the above Theorem~\ref{bijectivity} may be viewed as a
``limiting version'' of the Theorem~\ref{Max-Knus-Theorem} of Knus which
as noted in page~\pageref{coho-interpretation} may be interpreted as a
statement in cohomology.

 We shall see later ((b1), Theorem \ref{structure-of-specialisation}) that it is necessary to consider line-bundle-valued 
 quadratic forms to obtain the surjectivity of Theorem \ref{bijectivity} in 
 those cases for which a given $A$ representing an element in 
${\cal SPAZU}_4(X)$ is such that 
 $\hbox{\rm det}(A)\not\in 2.\hbox{\rm Pic}(X).$ 

It was shown in Part A, \cite{tevb-paper1} that algebra bundles
 belonging to ${\cal SPAZU}_4(X)$  
 are precisely the scheme-theoretic specialisations 
 (or limits) of rank 4 Azumaya algebra bundles on $X$
 (Theorem~\ref{smoothness-spazu}).
 
 Thus one may also restate the surjectivity as {\em schematic
   specialisations of rank 4 Azumaya bundles arise as 
   even-Clifford algebras of ternary quadratic bundles} and the
 injectivity as follows: {\em if the even-Clifford algebras of two ternary
   quadratic bundles are isomorphic, then the quadratic bundles are isometric
   upto tensoring by a twisted discriminant bundle.} 

 The main result of \cite{tevb-paper1}  
was the smoothness of the schematic closure of Azumaya algebra
 structures on a fixed 
 vector bundle of rank 4 over any scheme. 
 Part B  of \cite{tevb-paper1}  
 had applied this result 
 to obtain the generalised Seshadri desingularisation of the moduli
 space of semistable rank 2 degree
 zero vector bundles on a smooth proper curve relative to 
 a locally universally-japanese
 (Nagata) base scheme and also to obtain the generalised Nori desingularisation
 of the Artin moduli space of invariants of several matrices in rank 2
 over such a base scheme 
 together with good specialisation properties over $\mathbb Z.$ 
 The present work is concerned with applications to ternary quadratic
 forms  and the results obtained follow from an analysis of the 
computations that lead to the 
 smoothness.

 The good algebraic properties of Azumaya algebras
 are reflected as good geometric properties of the scheme of 
 Azumaya algebra structures on a fixed vector bundle: this scheme is 
 separated, of finite type and smooth relative to the base scheme (over 
 which the vector bundle is fixed) and also base-changes well 
 relative to the base scheme (Theorem~\ref{repbility-smoothness-azu},
 page~\pageref{repbility-smoothness-azu}).
  When the vector bundle is of rank 4, 
  the nice thing that happens is that all these good 
 properties also pass over to the limit i.e., to the 
  scheme of specialisations, defined 
 to be the schematic image of the scheme of Azumaya algebra structures 
 (Theorem~\ref{smoothness-spazu}, page~\pageref{smoothness-spazu}.)
 In the same vein, the present work shows that the theories of rank 3 
 semiregular quadratic forms and of rank 4 Azumaya algebras and
 their inter-relationships extend to
  the limit.

\subsection{Study of Groups of Similitudes}
\label{subsec1.2}
A bilinear form $b$, with values in a line bundle $I$, defined  
 on a vector bundle $V$ over the scheme $X$ 
induces an $I$-valued  quadratic form $q_b$ given on sections by 
 $x\mapsto b(x,x).$
 Let  $$L[I]:={\cal O}_X\oplus
\left(\bigoplus_{n> 0}(T^n(I) \oplus T^n(I^{-1}))\right)$$ be the 
 Laurent-Rees algebra of  $I$, where sections of $V$ (resp. of $I$) 
 are declared to be of degree one (resp. of degree two). 
 Then, as we saw in (2d), Theorem \ref{Bourbaki}, 
 $b$ naturally defines an $\mathbb Z$-graded 
 linear isomorphism $$\psi_b:\widetilde{C}(V,q_b, I)\cong
 \Lambda(V)\otimes L[I].$$ 
  In fact we have $$\psi_0:\widetilde{C}(V,0,I)=\Lambda(V)\otimes L[I].$$ 
Since in general   
 a quadratic bundle $\vqi$ on a non-affine scheme $X$
may not be induced from a global $I$-valued bilinear 
 form, one is unable to identify the $\mathbb Z$--graded vector bundle 
 underlying its generalised Clifford algebra bundle with
 $\Lambda(V)\otimes L[I].$
 The following 
 result overcomes this problem.  
\begin{proposition}\label{transfer-to-lambda2}
To every isomorphism of algebra-bundles
 $$\phi:\covqi\cong\covprqpripr,$$ one may 
 naturally associate  an 
 isomorphism of bundles $$\phi_{\Lambda^2}:\Lambda^2(V)\otimes I^{-1}
\cong\Lambda^2(V')\otimes {(I')}^{-1}$$ inducing a map  
$$\zeta_{\Lambda^2}:\isomcovqicovprqpripr\longrightarrow\isolambdatwovivpripr: \phi\mapsto\phi_{\Lambda^2}$$ 
where  $$\isomcovqicovprqpripr$$ is the set of algebra bundle
isomorphisms.
\end{proposition}
\begin{definition}\rm\label{defn-of-iso-sets}   
When $ V= V'$ and $I=I'$, we may thus denote  
  the subset of those $\phi$ for which
 $$\hbox{\rm det}(\phi_{\Lambda^2})\in\autdetlambdatwovi\equiv
\Gamma(X, {\cal O}_X^*)$$ is 
 a square by $$\isomcovqcovqprprimei$$ and those for which 
 $\hbox{\rm det}(\phi_{\Lambda^2})=1$ 
by the smaller subset $$\sisomcovqcovqpri.$$ 
 Taking $ q= q'$ in these sets and replacing `` {\rm Iso}'' by ``{\rm Aut}''
 in their notations respectively defines the  
 groups $$\autcovqi\supset\autcovqprimei\supset\sautcovqi.$$  
\end{definition}
\begin{theorem}\label{lifting-of-isomorphisms}
 For $I$-valued quadratic forms $q$ and $q'$ on a rank 3 vector bundle 
  $V$ over a scheme $X$, we have the following 
 commuting diagram of natural maps of sets with the downward arrows 
 being the canonical inclusions, the horizontal arrows being surjective 
 and the top horizontal arrow being bijective:  
$$\begin{CD}
\sisomvqvqpri @>{\cong}>> \sisomcovqcovqpri \\
 @V{\hbox{\rm inj}}VV @VV{\hbox{\rm inj}}V \\
\isomvqvqpri @>{\hbox{\rm onto}}>> \isomcovqcovqprprimei \\
  @V{\hbox{\rm inj}}VV @VV{\hbox{\rm inj}}V \\
\simvqvqpri @>{\hbox{\rm onto}}>> \isomcovqcovqpri 
\end{CD}$$
With respect to the surjections of the 
 horizontal arrows in the diagram above, we further have the following (where 
 $l$  is the function that associates a similarity to its 
 multiplier, $\hbox{\rm det}(g,l):=\hbox{\rm det}(g)$ for 
 an $I$-similarity $g$ with multiplier $l$ and $\zeta_{\Lambda^2}$
 is the map of
 Prop.\ref{transfer-to-lambda2} above):
\begin{description}
\item[(a)] there is a family of sections
 $$s_{2k+1}: \isomcovqcovqpri\longrightarrow\simvqvqpri$$
 indexed by the integers 
 such that 
$$l\circ s_{2k+1}
=\hbox{\rm det}^{2k+1}\circ\zeta_{\Lambda^2}\textrm{ and }(\hbox{\rm det}^2
\circ s_{2k+1})\times(l^{-3}\circ s_{2k+1})=\hbox{\rm det}
\circ\zeta_{\Lambda^2};$$
\item[(b)] there is also a 
 section 
$$s':\isomcovqcovqprprimei\longrightarrow \isomvqvqpri$$
 such that $\hbox{\rm det}^2\circ s'=\hbox{\rm det}\circ\zeta_{\Lambda^2};$
\item[(c)]  there is a family of sections 
 $$s_{2k+1}^{+}:\isomcovqcovqpri\longrightarrow\simvqvqpri$$ indexed by 
 the integers which is 
 multiplicative when followed by the natural inclusions
 into $\hbox{\rm GL}(V)\times\Gamma(X, {\cal O}_X^*)$, i.e., 
 if $$\phi_i\in\hbox{\rm Iso}[C_0(V,q_i, I), C_0(V,q_{i+1}, I)]$$ 
 then 
$$s_{2k+1}^{+}(\phi_2\circ\phi_1)=s_{2k+1}^{+}(\phi_2)\circ s_{2k+1}^{+}(\phi_1)\in\hbox{\rm GL}(V)\times\Gamma(X, {\cal O}_X^*).$$
 Further, 
$$l\circ s_{2k+1}^{+}
=\hbox{\rm det}^{2k+1}\circ\zeta_{\Lambda^2}\textrm{ and }(\hbox{\rm det}^2
\circ s_{2k+1}^{+})\times(l^{-3}\circ s_{2k+1}^{+})=\hbox{\rm det}
\circ\zeta_{\Lambda^2}.$$
\item[(d)] The maps $s_{2k+1}$ and 
 $s'$  above  may not be 
 multiplicative  but are
 mutliplicative upto $\mu_2(\Gamma(X, {\mathcal O}_X))$
 i.e., these maps followed  by the canonical 
 quotient map, on taking the quotient of 
 $\hbox{\rm GL}(V)\times\Gamma(X, {\cal O}_X^*)$ by
 $\mu_2(\Gamma(X, {\mathcal O}_X)).\hbox{\rm Id}_V\times\{1\}$,  
 become multiplicative. 
\end{description}
\end{theorem}

\begin{theorem}\label{lifting-of-automorphisms}
For a rank 3 quadratic bundle $\vqi$ on a scheme $X$, 
 one has the following natural 
 commutative diagram of groups with exact rows, where the downward arrows 
 are the  canonical inclusions and where 
 $l$ 
 is the function that associates to any $I$-(self)similarity its 
 multiplier:
\begin{displaymath}
\xymatrix{ 
  &  & \sovqi \ar[d]_{\hbox{\rm inj}} \ar[r]^{\cong\hspace*{1cm}}  &
 \sautcovqi \ar[d]^{\hbox{\rm inj}} & \\
 1\ar[r] & \mu_2(\Gamma(X, {\mathcal O}_X)) \ar[r]  \ar[d]_{\hbox{\rm inj}} & 
 \ovqi \ar[r] \ar[d]_{\hbox{\rm inj}} & \autcovqprimei \ar[d]^{\hbox{\rm inj}} 
 \ar[r]  & 1\\
 1\ar[r] & \Gamma(X, {\mathcal O}_X^{*}) \ar[r] & 
 \govqi \ar[r] \ar[d]_{\hbox{\rm det}^2\times l^{-3}} & \autcovqi \ar[d]^{\hbox{\rm det}} 
  \ar[r] 
  & 1\\
 & & \Gamma(X, {\mathcal O}_X^{*}) \ar@{=}[r]  & \Gamma(X, {\mathcal O}_X^{*}) & 
}
\end{displaymath}
Further, we have: 
\begin{description}
\item[(a1)] There are splitting homomorphisms 
 $$s_{2k+1}^{+}: \autcovqi\longrightarrow\govqi$$ such that 
$$l\circ s_{2k+1}^{+}=\hbox{\rm det}^{2k+1}\textrm{ and }(\hbox{\rm det}^2
\circ s_{2k+1}^{+})\times(l^{-3}\circ s_{2k+1}^{+})=\hbox{\rm det}.$$
In particular, $\govqi$ is a semidirect product.
\item[(a2)] The 
 restriction of $s_{2k+1}^{+}$ to $\autcovqprimei$ does not necessarily take 
 values in $\ovqi$, but the further restriction to $\sautcovqi$ does 
 take values in $\sovqi.$ 
\item[(a3)] The maps $s_{2k+1}$ and 
 $s'$ of Theorem \ref{lifting-of-isomorphisms} above (under the 
 current hypotheses) may not be 
 homomorphisms  but are homomorphisms upto $\mu_2(\Gamma(X, {\mathcal O}_X)).$ 
\item[(b)] Suppose $X$ is integral and $q\otimes\kappa(x)$ is semiregular
 at some point $x$ of 
 $X$ with residue field $\kappa(x).$
 Then any automorphism of $\covqi$ has determinant 1. Hence   
   $$\autcovqi=\autcovqprimei=\sautcovqi$$ and $\ovqi$ is the 
 semidirect product of $\mu_2(\Gamma(X, {\cal O}_X))$ and 
 $\sovqi.$
\end{description}
\end{theorem}
\noindent The proofs of the above results, and of the injectivity part of 
 Theorem \ref{bijectivity}, will be given in \S\ref{sec3} and \S\ref{sec4}. 
\subsection{Study of Bilinear Forms and Interpretation as Specialisations}
\label{subsec1.3}
As for the proof of the surjectivity part of  
Theorem \ref{bijectivity}, we have the following which will be proved 
 in \S\ref{sec5}:  
\begin{theorem}\label{structure-of-specialisation}
Let $X$ be a scheme and $ A$ a specialisation of rank 4 Azumaya algebra 
 bundles on $X.$  
 Let ${\cal O}_X. 1_A\hookrightarrow A$  be the line sub-bundle generated 
 by the  nowhere-vanishing 
 global section of $A$ corresponding to the unit
 for algebra multiplication.
\begin{description}
\item[(a)]
 There exist a rank 3 vector bundle 
 $V$ on $X$,  
 a quadratic form $q$ 
 on $V$ with values in the line bundle $I:=\hbox{\rm det}^{-1}(A)$, and 
 an isomorphism of algebra bundles 
  $A\cong \covqi.$ This gives the surjectivity in the statement of 
 Theorem \ref{bijectivity}.  Further,  the following linear 
 isomorphisms may be deduced: 
\begin{description}
\item[(1)] $\hbox{\rm det}(A)\otimes \Lambda^2({ V})\cong 
{ A}/{\cal O}_X. 1_{ A},$ from which follow:  
\item[(2)]  $\hbox{\rm det}(\Lambda^2({ V}))\cong{(\hbox{\rm det}({ A}))}^{\otimes-2};$
\item[(3)]  ${ V}\cong{({ A}/{\cal O}_X. 1_{ A})}^\vee\otimes
\hbox{\rm det}({ V})\otimes\hbox{\rm det}({ A});$
\item[(4)] $\hbox{\rm det}({ A}^\vee)\cong {(\hbox{\rm det}({ A}))}^{\otimes-3}\otimes {(\hbox{\rm det}({ V}))}^{\otimes-2}$ which 
implies that $\hbox{\rm det}({ A})\otimes\hbox{\rm det}({ A}^\vee)\in
2.\hbox{\rm Pic}(X).$ 
\end{description}
\item[(b)] There exists a quadratic bundle 
 $(V',q',I')$ such that $\covprqpripr\cong A$ 
 and with 
\begin{description}
\item[(1)] $I'={\cal O}_X$ iff $\hbox{\rm det}({ A})\in
2.\hbox{\rm Pic}(X)$;
\item[(2)] with $q'$ induced from a global $I'$-valued bilinear form 
 iff ${\cal O}_X.1_A$ is an ${\cal O}_X$-direct summand of $A$;
\item[(3)] with both $I'={\cal O}_X$ and with $q'$ induced from a 
 global bilinear form (with values in $I'$) iff 
 ${\cal O}_X.1_A$ is an ${\cal O}_X$-direct summand of $A$ and 
 $\hbox{\rm det}({ A})\in
2.\hbox{\rm Pic}(X).$
\end{description} 
\end{description}  
\end{theorem}

\noindent If $2\in\Gamma(X, {\cal O}_X^*)$, then any quadratic form is induced
from a symmetric bilinear form. Therefore from assertion (b2) of
the above, we have the following: 

\begin{corollary}\label{char-not-2-lifts-to-bilinear}
Suppose $2\in\Gamma(X, {\cal O}_X^*).$ Then for any specialisation $A$
of rank 4 Azumaya algebras on $X$, ${\cal O}_X.1_A$ is an ${\cal
  O}_X$-direct summand of $A.$ 
\end{corollary}

There are two ingredients in the proof of part (a) of
 Theorem~\ref{structure-of-specialisation}. The first is 
 Theorem \ref{lifting-of-isomorphisms}. The second is the following 
theorem which describes specialisations as bilinear forms under 
 certain conditions. 
 As a preparation towards its statement, we briefly 
 remind the reader of a few  results from  from Part A, \cite{tevb-paper1} 
(cf.~\S\ref{subsec2.7}).  

For a rank $n^2$ vector bundle $W$ 
 on a scheme $X$ and $w\in\Gamma(X,W)$ 
 a nowhere-vanishing global section, 
 recall that if  
  $\azuWw$ is 
  the open $X$--subscheme of Azumaya 
 algebra structures on $ W$ with identity $ w$ then 
 its schematic image (or the scheme of specialisations or the limiting 
 scheme) in the bigger $X$--scheme 
$\aaWw$ of associative 
 $ w$-unital algebra structures on $ W$  is the $X$--scheme  
 $\spazuWw.$ 
By definition, the set of distinct specialised $w$-unital 
 algebra structures on 
 $W$ corresponds precisely to 
the set of global sections of this last scheme over 
 $X$. 
 
If 
 $\hbox{\rm Stab}_{w}\subset \hbox{\rm GL}_{W}$
 is the stabiliser subgroupscheme of $w$,  recall from Theorems 
 3.4 and 3.8, Part A, \cite{tevb-paper1}, that there exists a canonical 
 action of $\hbox{\rm Stab}_{w}$ on
 $\spazuWw$ such that the 
 natural inclusions \label{clubsuit} 
 $$(\clubsuit)\hspace*{1cm}\azuWw\hookrightarrow \spazuWw \hookrightarrow \aaWw$$
 are all $\hbox{\rm Stab}_{w}$-equivariant.
  Now let
  $ V$ be a rank 3 vector bundle on the scheme $X$ and 
 $\hbox{\rm Bil}_{(V, I)}$
 be the associated  rank 9 vector bundle of 
 bilinear forms on $ V$ with values in the line bundle $I.$ 
 We say that a bilinear 
 form $b$ over an open subset $U\hookrightarrow X$ 
 is semiregular if there is a trivialisation $\{U_i\}$ of $I|U$, 
 such that  
 over each open subscheme $U_i$,  the  quadratic form $q_i$ with values 
 in the trivial line bundle induced from $q_b|U_i$ is semiregular (it
  may turn out that a semiregular bilinear form may be degenerate).
 This definition is independent of the choice of a 
 trivialisation, since $q_i$ is semiregular iff $\lambda q_i$ is semiregular 
 for every $\lambda\in\Gamma(U_i, {\cal O}_X^*)$ 
 (for further details see \S\ref{subsec2.2}).  
 In this way we obtain the open subscheme 
$$\hbox{\rm Bil}^{sr}_{(V,I)}\hookrightarrow
 \hbox{\rm Bil}_{(V,I)}$$
 of semiregular bilinear forms on $V$ with values in $I.$ We next take
 for $W$ the following special choice: 
 $$W:=\Lambda^{even}( V, I):=\bigoplus_{n\geq 0}\Lambda^{2n}(V)\otimes
 I^{-\otimes n}$$ and we let $ w\in\Gamma(X, 
 W)$ be the nowhere-vanishing global section corresponding to the 
 unit for the natural multiplication in the twisted
  even-exterior algebra bundle $W.$
There is an obvious
  natural action of $\hbox{\rm GL}_{V}$ on $\hbox{\rm Bil}_{(V,I)}.$ 
 There is 
also a natural morphism of groupschemes  $\hbox{\rm GL}_{V}\longrightarrow
\hbox{\rm Stab}_{w}$ given on valued points
 by $$g\mapsto \bigoplus_{n\geq 0}\Lambda^{2n}(g)\otimes\hbox{\rm Id}$$ 
 and therefore  
 the natural inclusions marked by  $(\clubsuit)$ above are  
 $\hbox{\rm GL}_{V}$-equivariant. Finally, note that there is an obvious 
 involution $\Sigma$ on $\aaWw$
 given by $A\mapsto \hbox{\rm opposite}(A)$ which leaves the open 
 subscheme $\azuWw$ invariant. 
\begin{theorem}\label{bilinear-forms-as-specialisations}\
\begin{description}
\item[(1)]
 Let $V$ be a rank 3 vector bundle on the scheme $X$,
 $ W:=\Lambda^{even}( V, I)$ and  $ w\in\Gamma(X, 
 W)$ correspond to 1 in the twisted even-exterior algebra bundle. 
 There is a natural $\hbox{\rm GL}_{V}$-equivariant 
 morphism of $X$--schemes 
$$\Upsilon'=\Upsilon'_X: \hbox{\rm Bil}_{(V,I)}\longrightarrow
\aaWw$$ whose 
 schematic image is precisely the scheme of specialisations 
$$\spazuWw.$$ Further 
if  $\Upsilon'$ factors canonically through 
$$\Upsilon=\Upsilon_X: \hbox{\rm Bil}_{(V, I)}\longrightarrow
\spazuWw,$$ 
 then $\Upsilon$ is a $\hbox{\rm GL}_{V}$-equivariant
  isomorphism and it maps the $\hbox{\rm GL}_{V}$-stable 
 open subscheme $ \hbox{\rm Bil}^{sr}_{(V, I)}$
 isomorphically onto the $\hbox{\rm GL}_{V}$-stable open 
 subscheme \newline $\azuWw.$
\item[(2)] The involution 
 $\Sigma$ of  $\aaWw$ 
 defines a unique involution (also denoted by $\Sigma$) on 
 the scheme of specialisations $\spazuWw$
 leaving the open subscheme $\azuWw$
 invariant, and therefore via the isomorphism $\Upsilon$, it defines 
 an involution on $\hbox{\rm Bil}_{(V,I)}.$ This involution is none other
 than the one on valued points given by $B\mapsto \hbox{\rm transpose}(-B).$  
\item[(3)]
For an $X$-scheme $T$, let $V_T$ (resp.$\thinspace W_T$, resp.$\thinspace I_T$) denote the 
 pullback of $V$ (resp.$\thinspace W$, resp.$\thinspace I$) to $T$,
  and let $w_T$ be the global section of $W_T$ induced by 
 $w.$ Then the base-changes of $\Upsilon'_X$ and $\Upsilon_X$ to $T$, 
 namely  
$$\Upsilon'_X\times_X T: \hbox{\rm Bil}_{(V, I)}\times_X T\longrightarrow
\aaWw\times_X T$$ and 
$$\Upsilon_X\times_X T: \hbox{\rm Bil}_{(V, I)}\times_X T\cong
\spazuWw\times_X T$$ may be  
 canonically identified with the corresponding ones over $T$ namely
 with 
$$\Upsilon'_T: \hbox{\rm Bil}_{(V_T, I_T)}\longrightarrow
\hbox{\rm Assoc}_{W_T,w_T}$$ and 
$$\Upsilon_T: \hbox{\rm Bil}_{(V_T, I_T)}\cong
\hbox{\rm SpAzu}_{W_T,w_T}\textrm{ respectively.}$$
\end{description} 
\end{theorem}
The explicit computation of the morphism $\Upsilon$ locally over $X$ is
 an important step in proving the above theorem. To describe this,  
 suppose that $I$ is trivial and  $V$ is free of rank 3 over $X$, so that 
 we may fix a basis $\{e_1, e_2, e_3 \}$ for $V$, which naturally gives rise 
 to a basis of $\hbox{\rm Bil}_{V}.$

  For any $X$-scheme $T$, 
 a $T$-valued point $B$ of $\hbox{\rm Bil}_{V}$ is just   
 a global bilinear form with values in ${\cal O}_T$
 on the pull-back $V\otimes_X T$ of 
 $V$ to $T.$  Such a $B$ 
 is given uniquely 
 by a $(3\times 3)$-matrix
 $(b_{ij})$ with the $b_{ij}$ being global sections of the trivial 
 line bundle ${\mathbb A}^1_T$ (or equivalently, elements of 
 $\Gamma(T, {\cal O}_T)$).
 The chosen basis for $V$ also gives rise to 
 the basis $$\{\epsilon_0:=w=1\hspace*{1mm};\hspace*{1mm}
\epsilon_1:=e_1\wedge e_2\hspace*{1mm},\hspace*{1mm}
\epsilon_2:=e_2\wedge e_3\hspace*{1mm},\hspace*{1mm}
\epsilon_3:=e_3\wedge e_1\}$$ 
of $W=\Lambda^{even}(V).$ A $T$-valued point $A$ of 
$\aaWw$ is just a
 $w_T:=(w\otimes_X T)$-unital 
 associative algebra structure on the bundle $W_T:=W\otimes_X T.$ 
 Let $\cdot_A$ denote the multiplication in the algebra bundle $A$, 
 and for ease of notation, let $s^{\circ}$ denote the section 
 $s\otimes_X T$ induced from a section $s$ (for example, $(w\otimes_X T)=
 w^{\circ}, \epsilon_i\otimes_X T=\epsilon_i^\circ$ etc).
\begin{theorem}\label{multiplication-table-wrt-bilinear}
In addition to the hypothesis of
 Theorem \ref{bilinear-forms-as-specialisations}, assume that $V$ is 
 free of rank 3 and that $I={\cal O}_X.$
 Then fixing a basis for $V$ and adopting 
 the notations above, 
 the map $\Upsilon(T)$ takes $B=(b_{ij})$ to 
$(A, 1_A, \cdot)=(W_T, w_T=w^\circ , \cdot_A)$
 with multiplication given as follows, where $M_{ij}(B)$ is the 
 determinant of the minor of the element $b_{ij}$ in $B:$
\begin{enumerate}
\item[$\bullet$] 
$\epsilon_1^{\circ}\cdot_A\epsilon_1^{\circ}=-M_{33}(B)w^\circ + (b_{21}-b_{12})\epsilon_1^{\circ}$
\item[$\bullet$] 
$\epsilon_2^{\circ}\cdot_A\epsilon_2^{\circ}=-M_{11}(B)w^\circ + (b_{32}-b_{23})\epsilon_2^{\circ}$
\item[$\bullet$] 
$\epsilon_3^{\circ}\cdot_A\epsilon_3^{\circ}=-M_{22}(B)w^\circ + (b_{13}-b_{31})\epsilon_3^{\circ}$
\item[$\bullet$] 
$\epsilon_1^{\circ}\cdot_A\epsilon_2^{\circ}=-M_{31}(B)w^\circ -b_{23}\epsilon_1^{\circ} -b_{12}\epsilon_2^{\circ} 
-b_{22}\epsilon_3^{\circ}$
\item[$\bullet$] 
$\epsilon_2^{\circ}\cdot_A\epsilon_3^{\circ}=+M_{12}(B)w^\circ -b_{33}\epsilon_1^{\circ} -b_{31}\epsilon_2^{\circ} 
-b_{23}\epsilon_3^{\circ}$
\item[$\bullet$] 
$\epsilon_3^{\circ}\cdot_A\epsilon_1^{\circ}=+M_{23}(B)w^\circ -b_{31}\epsilon_1^{\circ} -b_{11}\epsilon_2^{\circ} 
-b_{12}\epsilon_3^{\circ}$
\item[$\bullet$] 
$\epsilon_1^{\circ}\cdot_A\epsilon_3^{\circ}=+M_{32}(B)w^\circ +b_{13}\epsilon_1^{\circ} +b_{11}\epsilon_2^{\circ} 
+b_{21}\epsilon_3^{\circ}$
\item[$\bullet$] 
$\epsilon_2^{\circ}\cdot_A\epsilon_1^{\circ}=-M_{13}(B)w^\circ +b_{32}\epsilon_1^{\circ} +b_{21}\epsilon_2^{\circ} 
+b_{22}\epsilon_3^{\circ}$
\item[$\bullet$] 
$\epsilon_3^{\circ}\cdot_A\epsilon_2^{\circ}=-M_{21}(B)w^\circ +b_{33}\epsilon_1^{\circ} +b_{13}\epsilon_2^{\circ} 
+b_{32}\epsilon_3^{\circ}$
\end{enumerate}
\end{theorem}
The key to the proofs of Theorems \ref{bijectivity}, 
\ref{lifting-of-isomorphisms} and \ref{lifting-of-automorphisms} 
lies in an analysis of a different 
 identification 
 of the scheme of specialisations, 
 namely one related to the scheme of ${\cal O}_X$-valued quadratic 
 forms on a trivial rank 3 bundle in the special situation when $W$ is free 
 and $w$ part of a global basis. Without loss of 
 generality we may in this situation therefore 
 take $V$ to be a free rank 3 vector bundle on $X$ and 
 $$(W, w)=(\Lambda^{even}(V), 1),$$ so that we are in the situation of 
 Theorem \ref{multiplication-table-wrt-bilinear} above. This relationship with 
 quadratic forms  was shown in Theorem 5.3, Part A, \cite{tevb-paper1}, 
  which we briefly recall next.
 Let $\hbox{\rm Quad}_{V}$ denote the 
 bundle of quadratic forms on $V$ (with values in ${\cal O}_X$) and 
 $\hbox{\rm Quad}^{sr}_{V}$ the open subscheme of semiregular quadratic forms. 
 Let $A_0$ denote the algebra bundle structure (with unit $w=1$) on
 $W=\Lambda^{even}(V)$ given by $\Lambda^{even}(V)$ itself. Fix a basis for 
 $V$ and adopt the notations preceding Theorem \ref{multiplication-table-wrt-bilinear} above. Then $\hbox{\rm Stab}_{w}$ is the semidirect product of a 
 commutative 3-dimensional subgroupscheme 
$$\hbox{\sf L}_w\cong ({\mathbb A}^3_X, +)$$
 with the stabiliser subgroupscheme $\hbox{\rm Stab}_{A_0}$ of $A_0$ in 
 $\hbox{\rm Stab}_{w}$ (Lemma 5.1, Part A, \cite{tevb-paper1}).  
\begin{theorem}[Definition 5.2 \& Theorem 5.3, Part A, \cite{tevb-paper1}]
\label{Thetaisom}
There is a \newline natural isomorphism 
$$\Theta: \hbox{\rm Quad}_{V}\times_X \hbox{\sf L}_w \cong
 \spazuWw$$ which maps the open subscheme 
 $\hbox{\rm Quad}^{sr}_{V}\times_X \hbox{\sf L}_w$ isomorphically onto the open subscheme 
 $\azuWw.$  
\end{theorem}
The isomorphism $\Theta$ was first defined by C.~S.~Seshadri in
\cite{css-desing} for the case $X=\textrm{Spec}(k)$, $k$ an
algebraically closed field of characteristic $\neq 2.$ 
Section \ref{sec4} is essentially devoted to studying $\Theta.$
 There we compute $\Theta$ explicitly and
 in Theorem \ref{mult-table-wrt-Theta} we write out the 
  multiplication table of every specialised algebra structure on any fixed
  free rank 4 vector bundle with fixed unit that is part of a global basis.
It turns out that $\Theta$ is not equivariant 
 with respect to $\hbox{\rm GL}_{V}$, but nevertheless satisfies a `twisted' 
 form of equivariance (Theorem \ref{tequiv}). 
A $T$-valued point $q$ of 
$$\hbox{\rm Quad}_{V}\cong{\mathbb A}^6_X$$ 
may be identified with a 6-tuple $(\lambda_1,\lambda_2,\lambda_3,
\lambda_{12},\lambda_{13},\lambda_{23})$ corresponding to the 
 quadratic form 
$$(x_1,x_2,x_3)\mapsto \Sigma_i \lambda_i x_i^2 
+ \Sigma_{i<j}\lambda_{ij}x_ix_j.$$
 A $T$-valued point $\underline{t}$
  of $\hbox{\sf L}_w\cong({\mathbb A}^3_X, +)$
 may be identified with a 3-tuple $(t_1, t_2, t_3)$ which 
 corresponds to the valued point of $\hbox{\rm Stab}_{w}$ given by the 
 $(4\times 4)$-matrix
 \begin{displaymath}
 \left(\begin{array}{c|ccc}
1 & t_1 & t_2 & t_3 \\
\hline
0 & & I_3 & 
\end{array}
\right)
 \end{displaymath} 
where $I_3$ is the $(3\times 3)$-identity matrix. 
 With these notations, 
the identification of Theorems \ref{bilinear-forms-as-specialisations} and 
 \ref{multiplication-table-wrt-bilinear} may be compared with that of 
 the above Theorem \ref{Thetaisom} as follows.

\begin{theorem}\label{Thetaisom-Upsilonisom}
The isomorphism $\Upsilon^{-1}\circ\Theta
:\hbox{\rm Quad}_{V}\times_X \hbox{\sf L}_w
\cong \hbox{\rm Bil}_{V}$ takes the valued point 
$$(\hspace*{2mm}q,\hspace*{2mm} \underline{t}\hspace*{2mm})=(\hspace*{2mm}(\lambda_1,\lambda_2,\lambda_3,
\lambda_{12},\lambda_{13},\lambda_{23}),\hspace*{2mm} (t_1, t_2,
t_3)\hspace*{2mm})$$ 
to the valued point $B=(b_{ij})$ given by 
\begin{displaymath}
B=\left(\begin{array}{ccc}
\lambda_1 & t_1 & \lambda_{13}-t_3 \\
\lambda_{12}-t_1 & \lambda_2 & t_2 \\
t_3 & \lambda_{23}-t_2 & \lambda_3
\end{array}
\right).
 \end{displaymath} 
Moreover, under this identification, the involution $B\mapsto {(-B)}^t$ on 
  $\hbox{\rm Bil}_{V}$ (induced from the isomorphism 
 $\Upsilon$ of Theorem \ref{bilinear-forms-as-specialisations})
 translates into the involution on 
 $$\hbox{\rm Quad}_{V}\times_X \hbox{\sf L}_w$$ given by 
 $$(\hspace*{2mm}q,\hspace*{2mm} (t_1, t_2, t_3)\hspace*{2mm})\mapsto (\hspace*{2mm}-q,\hspace*{2mm} (t_1-\lambda_{12}, t_2-\lambda_{23}, t_3-\lambda_{13})\hspace*{2mm}).$$ 
\end{theorem}
\subsection{Degenerations of Azumaya Bundles as Quaternion Algebra Bundles}
\label{subsec1.4}
We next make some comments on specialised algebras.
 Let $R$ be a unital commutative ring and $A$ a unital associative $R$-algebra.
  Given any involution $\sigma$ on $A$, we may 
 define the trace and norm associated to this involution by  
$$tr_\sigma:x\mapsto 
 x+\sigma(x)\textrm{ and }n_{\sigma}: x\mapsto x.\sigma(x).$$  
 In para.1.3, Chap.I, \cite{Knus}, Knus calls $\sigma$
 {\em standard}  if 
 $\sigma$ fixes $R.1_A$ and both $tr_\sigma$ and $n_{\sigma}$ take 
 values in $R.1_A.$ In  Prop.1.3.4 of the same chapter, he proves that
  a standard  involution is unique if it exists, provided the $R$-module
 underlying  $A$  is finitely generated projective and faithful.

 In para.1.3.7, {\em op.cit.},   Knus defines  
 $A$ to be 
 a {\em quaternion} algebra if $A$ is a projective $R$-module of rank 4 and 
 $A$ has a standard involution. Thus we may define a rank 4 algebra bundle on 
 a scheme $X$ to be a quaternion algebra bundle if it is locally 
 (in the Zariski topology) a quaternion algebra in Knus' sense, and  
 it would follow that the local standard involutions glue to define a
 unique global standard involution on the bundle. 
  \begin{proposition}\label{spln-is-quaternion}  
Any specialised algebra bundle is a quaternion algebra bundle.
\end{proposition}
 This result 
 can be deduced from the following two facts:
\begin{description}
\item[(1)] Any specialised algebra is 
 locally (in the Zariski topology) the even Clifford algebra of a rank 
 3 quadratic bundle (assertion (3), Theorem~\ref{smoothness-spazu} and Theorem~\ref{Thetaisom}).
\item[(2)] The even Clifford algebra of a quadratic 
  module of rank 3 over a commutative ring
 has a standard involution which is none other than the restriction of the `standard' involution on the full 
 Clifford algebra (Prop.3.1.1, Chap.V, \cite{Knus}). {\bf Q.E.D.}
\end{description}
We hasten to remark that even 
 over an algebraically closed field there are quaternion 
 algebras that are not even-Clifford algebras of quadratic forms.

The proof of the bijection stated in Theorem \ref{Max-Knus-Theorem} 
  follows from Prop.3.2.3 and Prop.3.2.4, Chap.V, \cite{Knus}
 generalised to the scheme-theoretic setting.
 We recall how the surjectivity is established. 

 Let $A$ be a specialised algebra bundle on the scheme $X$. 
 By the results just quoted if $A$ is Azumaya, or more generally by
 Prop.\ref{spln-is-quaternion}, we have the existence of a unique 
 standard involution $\sigma_A$  on $A$, to which are associated the 
 norm 
$$n_{\sigma_A}:A\longrightarrow {\mathbb A}^1_X\textrm{ 
given on sections by }x\mapsto x.\sigma_A(x)$$ and 
 the trace 
$$tr_{\sigma_A}:A\longrightarrow {\mathbb A}^1_X\textrm{
 given on sections by }x\mapsto x+\sigma_A(x).$$ 
 Let $A':=\hbox{\rm kernel}(tr_{\sigma_A})\hookrightarrow A$ 
be the subsheaf of  trace zero elements.
As the calculations in para.3.2, Chap.V, \cite{Knus} show,  
the trace map is surjective if 
 $A$ is itself an Azumaya algebra; if this is the case, then 
 it is further shown there that
 the rank 3 quadratic bundle 
$$\vqi:=(A', n_{\sigma_A}|A', {\cal
 O}_X.1_A)$$
 is  semiregular and its even Clifford algebra  $\covqi\cong A.$
 However the above method of retrieving a canonical 
 rank 3 quadratic bundle fails badly for specialised non-Azumaya algebras. 
 Consider even the case of $X=\hbox{\rm Spec}(k)$ where $k$
 is a field of characteristic 
 two and the Clifford algebra $A=C_0(V,q)$ 
 of a quadratic form $q$ on $V=k^{\oplus 3}$
 which is a perfect 
 square (i.e., a square of a linear form or equivalently a sum of squares). 
 In this case an easy computation shows that 
 the subspace $A'$ of trace zero elements is the full space $A.$ 
\begin{proposition}\label{semilocal-unique-semireg}
Let $S$ be a commutative semilocal ring that is 2-perfect i.e., such that  
the square map 
$$S\longrightarrow S: s\mapsto s^2$$ is surjective, and 
 $V$ a free rank 3 $S$-module. Then the set of semiregular quadratic $S$-forms 
 on $V$ forms a single $\hbox{\rm GL}(V)$--orbit; in other words, upto 
 isometry, $\exists$ only one semiregular quadratic $S$-module structure 
 on $V.$ 
\end{proposition} 
\begin{corollary}\label{local-unique-azumaya}
Let $S$ be a commutative local ring that is 2-perfect. Then any two rank 4 
 Azumaya $S$-algebras are isomorphic. If $S$ were only semilocal, the 
 conclusion still holds provided the identity elements for multiplication 
 for each of the two Azumaya $S$-algebras can be completed to an $S$-basis.
\end{corollary}
 The proof of the above Proposition will be given in \S \ref{sec7}. 
 In view of Theorem \ref{Max-Knus-Theorem},
 taking $X=\hbox{\rm Spec}(S)$ with $S$ as in 
 Prop.\ref{semilocal-unique-semireg} proves the first
 assertion of the above corollary.
 The second may be deduced by an application of Theorem \ref{Thetaisom} 
 alongwith Prop.\ref{semilocal-unique-semireg}.

\subsection{Non-existence of Azumaya Structures and Semiregular Forms}
\label{newsubsec1.5}
\noindent In this and the next subsections, we study what happens when we
impose the condition on self-duality on the bundle underlying an algebra
(or on one on which a quadratic form is defined). In this subsection, our
aim is to indicate examples of rank 4 vector bundles on which there
do not exist any global Azumaya structures; 
for these examples it also turns out that  
 there do not exist  
global regular quadratic forms with values in the trivial line bundle. 
However, in these examples there do exist algebra structures which are 
Azumaya on a nontrivial dense open subscheme and there do exist 
quadratic forms with values in the trivial bundle which are
semiregular on a dense open subscheme.  
We use Theorem~\ref{bijectivity} to obtain examples of rank three bundles on
which there do not exist any global semiregular quadratic forms with
values in certain line bundles; however, there do exist forms which
are semiregular on a dense open subscheme. We shall only outline the
results and the proofs will be given in \cite{tevb-paper3}. 

In the following we let $X$ be a scheme and $W$ a rank $n^2$ vector bundle 
 over $X.$
 The following result describes the behaviour of 
 the locus where an algebra structure is Azumaya. 
\begin{proposition}[Prop.3.3, Part A, \cite{tevb-paper1}]
\label{sch-azu-props}\ 
\begin{description}
\item[1.] 
Let $T$ be an $X$-scheme and $A$ an associative unital algebra 
 structure on $W_T:=W\otimes_X T.$ 
 Then the subset
$$U(T,A):=\{ t\in T\ |\ {\cal A}_t \hbox{ is an Azumaya }{\cal O}_{T,t}-\hbox{algebra}
\}$$ is an open (possibly empty) subset. When $U(T, A)$  is nonempty,
 denote by the same symbol the  
  canonical open subscheme structure. Then if 
 $f:T'\longrightarrow T$ is an $X$-morphism such that 
 the topological image intersects $U(T, A)$, then 
 $$U(T', f^*(A)=A\otimes_T T')\cong U(T, A)\times_{T}T'$$ as open 
 subschemes of $T'.$  Further  
  $U(T, A)\hookrightarrow  T$ is an affine morphism. 
\item[2.] $U(T,A)$ is the maximal open subset restricted to which 
 $A$ is  Azumaya. 
\item[3.] Further let $f:T'\longrightarrow T$ be a morphism of $X$-schemes 
 such that $f^*(A)$ is Azumaya. Then $f$ factors 
 through the open subscheme $U(T,A)$ defined above. 
\end{description}
\end{proposition}
\noindent Next we let $V$ denote a vector bundle and $I$ a
 line-bundle 
 over $X.$ The following result describes the behaviour of the locus
 where a quadratic form is good. 
\begin{proposition}
\label{sch-semireg-props}\ 
\begin{description}
\item[1.] 
Let $T$ be an $X$-scheme and $q:V_T\rightarrow I_T$ a 
 quadratic bundle with $V_T:=V\otimes_X T$ and $I_T:=I\otimes_X T.$
 Then the subset
$$U(T,q):=\{ t\in T\ |\
 q\otimes\kappa(t):V_T\otimes\kappa(t)\rightarrow
 I_T\otimes\kappa(t)\textrm{ is ``good''}\},$$
 where ``good'' means semiregular if $V$ is of odd rank and regular 
 if $V$ is of even rank, 
 is an open (possibly empty) subset. When $U(T, q)$  is nonempty,
 denote by the same symbol the  
  canonical open subscheme structure. Then if 
 $f:T'\longrightarrow T$ is an $X$-morphism such that 
 the topological image intersects $U(T, q)$, and if
 $f^*(q):=q_{T'}:V_{T'}\rightarrow I_{T'}$ is the induced quadratic bundle on 
$T'$,  then 
 $$U(T', q_{T'})\cong U(T, q)\times_{T}T'$$ as open 
 subschemes of $T'.$  Further  
  $U(T, q)\hookrightarrow  T$ is an affine morphism. 
\item[2.] $U(T,q)$ is the maximal open subset restricted to which 
 $q$ is good. 
\item[3.] Further let $f:T'\longrightarrow T$ be a morphism of $X$-schemes 
 such that $q_{T'}$ is good. Then $f$ factors 
 through the open subscheme $U(T,q)$ defined above. 
\end{description}
\end{proposition}
\noindent 
We now let $S$ be a scheme, and let $X\rightarrow S$ be an $S$-scheme which is
proper, of finite-type and has connected fibers relative to $S.$
Again, $W$ denotes a vector bundle on $X$, $A$ an algebra structure on
$W$, and $q:V\rightarrow I$ a quadratic form on the vector bundle $V$
with values in the line bundle $I.$ 
\begin{theorem}\label{relations-with-self-duality}
(With the above notations)\begin{description}
\item[1.] Suppose $W$ is self-dual and  $U(X,A)\rightarrow S$ is
  surjective.
Then $U(X,A)=X$ i.e., $A$ is Azumaya. 
\item[2.] Suppose that the rank of $W$ is 4, that there exists an algebra
  structure $A'$ on $W$ such that $U(X,A')\rightarrow S$ is surjective
  and that $U(X,A')\neq X.$ Then there does not exist any global
  Azumaya algebra structure on $W.$ Neither does there exist any
  global regular quadratic form on $W$ with values in the trivial line
  bundle.  
\item[3.] Let $V$ be of rank 3 and suppose that 
the underlying bundle of $\covqi$ is self-dual. If $U(X,q)\rightarrow
S$ is surjective, then $q$ is semiregular. 
\item[4.] Let $V$ be of rank 3 such that 
$U(X,q)\rightarrow S$ is surjective but $U(X,q)\neq X.$ Then 
there does not exist any global semiregular quadratic form on $V$ with
values in $I.$ If $(L,h,J)$ is a twisted discriminant bundle, then 
there does not exist any global semiregular quadratic form on
$V\otimes L$ with values in $I\otimes J.$ 
\end{description}
\end{theorem}
\noindent Assertion (3) follows from (1) and
Theorem~\ref{bijectivity}. 
 The proofs may be reduced to the case
$B=\textrm{Spec}(k)$ where $k$ is a field. In this case, assertions
(1) and (2) may be restated as follows. 
Let $X$ be a connected proper scheme of finite type 
  over a field and let 
  $W$ be a  vector bundle on $X$ of rank $n^2$ for some $n\geq 2.$  
\begin{description}
\item[(a)]
 Let $W$ be self-dual 
  and $A$ an associative unital algebra structure on $W$.
If a section to $\aaWw$ over $X$  
 meets $\azuWw$ topologically, then it
factors as a morphism through the open subscheme 
 $\azuWw$, where 
 $w:=1_A$ and $A$ corresponds to the given section;
\item[(b)] Let the rank of $W$ be 4. 
if there is a section over $X$ of $\aaWw$ that topologically meets both the 
 open subscheme $\azuWw$ and its complement (with 
 $w=1_A$ where $A$ corresponds to the given section),
 then the $X$-schemes \newline $\hbox{\rm Assoc}_{W,w'}$ (with 
 $w'$ global nowhere-vanishing) cannot have sections that land 
 topologically inside $\hbox{\rm Azu}_{W,w'}$ and 
 hence in particular the $X$-schemes 
  $\hbox{\rm Azu}_{W,w'}$ have 
  no sections over $X.$ 
\end{description}
\noindent We leave it to the reader to formulate 
similar statements corresponding to the other assertions of
Theorem~\ref{relations-with-self-duality}.   
 We next indicate situations where the above results apply 
 to give examples of 
 rank 4 bundles that do not admit
any Azumaya algebra structures and of
 rank 3 vector bundles which do not admit any
semiregular form with values in certain line bundles. The examples 
 occur naturally on certain moduli
 spaces of vector bundles over a relative curve. 
 As remarked 
earlier, 
we quickly mention the relevant objects, but we shall not get
 into definitions, proofs etc, which will appear in
 \cite{tevb-paper3}. 

We first recall some facts about Nagata Rings.
The standard reference is Chap.12 of Matsumura's book \cite{matsumura}.
An integral domain $A$ is said to satisfy condition N-1 if its integral 
 closure $A_K$ in its quotient field $K$ is a finite $A$-module. It is 
 said to satisfy condition N-2 if for every finite extension field 
 $L/K$, the integral closure $A_L$ of $A$ in $L$ is a finite $A$-module. 
The properties N-1 and N-2 are preserved under localisation and N-2 $\implies$ 
 N-1 whereas noetherianness with N-1 $\implies$ N-2 only in char.0; there 
 exists an example of Y. Akizuki of a noetherian domain of positive char. 
 which is not N-1. A commutative ring $B$ is called a {\em  Nagata} ring 
 ({\em  pseudo-geometric} ring in Nagata's own terminology and 
 {\em  universally japanese} ring in Grothendieck's) if it is noetherian 
 and $B/{\mathfrak p}$ is N-2 for each prime $\mathfrak p$ of $B.$ 
 Every localisation of $B$ and every finitely generated (commutative) 
 $B$-algebra are then also Nagata, and complete noetherian local rings 
 are Nagata as well. Dedekind domains of characteristic zero
  such as $\mathbb Z$ are 
 Nagata.

Let $S$ be a normal integral scheme, which we shall assume is
 locally-universally Japanese
(Nagata). This means that there is a covering of $S$ by affine open
 subschemes that are isomorphic to spectra of Nagata rings. 
 Such an assumption on $S$ is necessary to be able to obtain
 finite-type quotients by reductive groupschemes following Seshadri's 
 Geometric Invariant Theory over a General Base as in \cite{css-gitoab}.

   Let $C\rightarrow S$ be 
smooth projective of relative dimension 1 (i.e., a curve relative to
$S$) with geometric fibres irreducible and of constant genus $g\geq
2.$
 Let $\ucss$
denote 
the Seshadri moduli space of 
semistable vector bundles on $C$ of rank 2 and degree zero,
generalising the construction in 
 \cite{css-moduli} to the base $S$ using the methods of
 \cite{css-gitoab}.
 It is a
coarse moduli scheme, obtained as a Geometric Invariant Theory
quotient, and is a normal geometrically 
irreducible projective scheme of relative dimension 
$4g-3$ over $S.$ 
 Let $\ucs\subset\ucss$ denote the open subscheme corresponding to 
 stable vector bundles. Then $\ucs$ is precisely
the locus where $\ucss\rightarrow S$ is smooth, unless $g=2$ in which
case $\ucss\rightarrow S$ is itself smooth.
The Seshadri-desingularisation 
$$\pi_2: \nc\rightarrow \ucss$$
may be constructed, extending the case $S=\textrm{Spec}(k)$, $k$ a
field, carried out in \cite{css-desing} for $\textrm{Char.}(k)\neq 2$ and 
in a characteristic-free way 
 in \S~6,
Part B of \cite{tevb-paper1}. 
$\pi_2$  is an isomorphism over the open 
subscheme $\ucs.$
The scheme $\nc$ appears as a
 fine moduli space for certain parabolic stable vector
bundles of rank 4 and degree 0 on $C$ and it turns out to be a
geometrically 
irreducible smooth projective scheme relative to $S.$ 
 By the very construction of
$\nc$, there is a naturally defined rank 4 vector bundle $W$ on $\nc$
which is equipped with an associative unital algebra structure $A.$
Further, the locus where $A$ is Azumaya is precisely the dense open
subscheme $\ncs:=\pi_2^{-1}(\ucs).$  Let $\vqi$ be any ternary quadratic
  bundle on $\nc$ such that  $\covqi\cong A$ (such a choice exists by 
the surjectivity part of Theorem~\ref{bijectivity}). Applying
Theorem~\ref{relations-with-self-duality} to the present situation
gives the following. 
\begin{theorem}\label{non-existence-examples} (With the above notations:) 
\begin{description}
\item[1.]
$W$  does not admit any
(global) Azumaya algebra structure. 
\item[2.] There does not exist any global
  regular quadratic form on $W$ with values in the trivial line bundle. 
\item[3.] There does not exist any global
  semiregular quadratic form on $V$ with values in $I.$ If  
$(L,h,J)$ is a twisted discriminant bundle, there does not exist any 
global semiregular quadratic form on $V\otimes L$ with values in 
 $I\otimes J.$ 
\end{description} 
\end{theorem}
\subsection{Algebra Structures on Self-Dual Bundles}
\label{subsec1.5}
\noindent In this subsection, we investigate again the effect of the
 condition of self-duality on the vector bundle underlying an
 algebra. Our aim is to study the Picard group of the scheme of
 specialisations of Azumaya structures on a fixed rank 4 vector bundle
 over suitable schemes. Recall that an integral separated noetherian
 scheme is said to be locally-factorial if each of its 
 local rings is a unique factorisation domain (=UFD=factorial ring). 
The proofs of the following results will be given in \S\ref{sec6}.
\begin{theorem}\label{surjectivity-affine-limited}
Let $X$ be a scheme and  $W$ a rank 4 vector bundle on $X$ with  
 a global nowhere-vanishing section $w.$ Let $D_X$ denote  
 the closed subset 
$$\spazuWw\backslash\azuWw.$$
\begin{description}
\item[(a1)] $X$ is irreducible iff $\spazuWw$ is irreducible  iff 
 $\azuWw$ is irreducible iff $D_X$ is irreducible.
\item[(a2)] The set of 
 irreducible components of $X$ is locally finite---for example this 
 happens when  $X$ is locally noetherian---iff the same is 
 true of the corresponding set for $\spazuWw$ or for $\azuWw.$
\item[(a3)] If 
 $X$ is noetherian and finite-dimensional then the same are true for 
 \newline $\spazuWw$ and  
 $\azuWw.$
\item[(b)] If $X'\longrightarrow X$ is a morphism of schemes,  
 and if $(W',w')$ denotes the pullback of $(W,w)$, then we have  a 
 canonical isomorphism 
$$\azuWprwpr\cong\azuWw\times_{\spazuWw} \spazuWprwpr$$
In particular, the topological image of $D_{X'}$ is $D_{X}.$
Moreover, when $X'\longrightarrow X$ is a homeomorphism onto its 
 topological image---which is for example the case when it is a closed 
 or an open immersion, then $$D_X\cap\spazuWprwpr$$ can be identified with 
 $D_{X'}.$ 
\item[(c)] 
If $X$ is a scheme which is finite dimensional and 
 whose set of irreducible components is locally finite,
 then the closed subset $D_X$ is a divisor i.e., it has codimension 
 1 in $\spazuWw.$
\item[(d)]  $X$ is affine iff 
 $\spazuWw$ is affine iff 
$\azuWw$
 is affine. If $X$ is regular in codimension 1 (respectively locally-factorial)
 then so are $\spazuWw$ and $\azuWw.$ 
\item[(e)] Assume that $X$ is locally-factorial. Then the canonical
  homomorphism $$\hbox{\rm Pic}(X)\rightarrow\hbox{\rm Pic}(\spazuWw)$$ is
  an isomorphism. 
\item[(f)] 
 Assume that $X$ is 
 locally-factorial and $W$ is self-dual (i.e., $W\cong W^\vee$).  
 Then the (Weil) divisor $n.(D_X)$ is principal for some positive integer 
 $n$,   
 so that the natural homomorphism given by restriction of line bundles 
 $$\hbox{\rm Pic}(\spazuWw)\longrightarrow 
 \hbox{\rm Pic}(\azuWw)$$ is an isomorphism iff 
 $n=1.$ In this case, for any $X$-scheme $T$, any specialised algebra
 structure on $W\otimes_X T$ arises from a quadratic form with values
 in the trivial line bundle. 
\item[(g)] Assume that $X$ is locally-factorial and that there exists
  an Azumaya algebra structure on $W$ with unit $w.$ 
 Then the canonical homomorphisms 
  $$\hbox{\rm Pic}(X)\rightarrow\hbox{\rm Pic}(\azuWw)\textrm{  and   }
  \hbox{\rm Pic}(\spazuWw)\rightarrow\hbox{\rm Pic}(\azuWw)$$ are
  isomorphisms and the divisor $D_X$ is principal. 
\end{description}
 \end{theorem}
\subsection{Stratification of the Variety of Specialisations} 
\label{subsec1.6}
Let $W$ be a rank 4 vector bundle on a scheme $X$, $w\in\Gamma(X,W)$ 
 a nowhere-vanishing global section and $\hbox{\rm Stab}_{w}\subset \hbox{\rm 
 GL}_{W}$ the stabiliser subgroupscheme of $w.$ Recall from page
 \pageref{clubsuit} that the 
 natural inclusion $$\azuWw\hookrightarrow\spazuWw$$ is 
 $\hbox{\rm Stab}_{w}$-equivariant.  
 When $X=\hbox{\rm Spec}(k)$ where $k$ is an algebraically closed field,
 there is a canonical $\hbox{\rm Stab}_{w}$-stratification of the 
 $k$-variety underlying $\spazuWw$ as
 follows (the proof will be given in \S\ref{sec7}). 
\begin{theorem}\label{stratification}\ 
\begin{description}
\item[(1)] Let $k$ be a quadratically closed field and $X=\hbox{\rm Spec}(k).$ 
Then the set of ternary quadratic modules upto similarity 
  has 4 elements which correspond 
 to
\begin{description}
\item[(a)] semiregular quadratic modules; 
\item[(b)] rank 2 quadratic modules i.e., 
 those that are not semiregular but which are regular on a two-dimensional 
 subspace;
\item[(c)] nonzero perfect squares (= squares of linear forms) and 
\item[(d)] the zero form.
\end{description} If $V$ is a 3-dimensional
 vector space over $k$ and $\{e_1, e_2, e_3\}$ a $k$-basis for $V$, then 
 representatives for these 4 $\hbox{\rm GL}_{V}$-orbits
 in the space $ \hbox{\rm Quad}_{V}$ 
 of quadratic forms on $V$ can respectively be taken to be: 
\begin{description}
\item[(a)] $q^{(1)}(\Sigma_{i=1}^3x_ie_i)=x_1x_2+x_3^2;$
\item[(b)] $q^{(2)}(\Sigma_{i=1}^3x_ie_i)=x_1x_2;$
\item[(c)] $q^{(3)}(\Sigma_{i=1}^3x_ie_i)=x_3^2\textrm{ and }$
\item[(d)] $q^{(4)}=0.$
\end{description}
\item[(2)] In addition to the hypotheses and notations of (1) above, 
 assume that $k$ is an algebraically closed field. Then the 
 four orbits 
$$\hbox{\rm Quad}_{V}^{(i)}:=\hbox{\rm GL}_{V}\cdot q^{(i)}\textrm{ for }
 1\leq i\leq 4$$
 form a stratification of the $k$-variety $\hbox{\rm Quad}_{V}$ 
 in the sense that we have  
 $$\overline{\hbox{\rm Quad}_{V}^{(1)}}=\hbox{\rm Quad}_{V}\hspace*{0.25cm}and\hspace*{0.25cm}
  \overline{\hbox{\rm Quad}_{V}^{(i+1)}}=\overline{\hbox{\rm Quad}_{V}^{(i)}}\hspace*{2mm}\backslash\hspace*{2mm}
\hbox{\rm Quad}_{V}^{(i)},\hspace*{0.25cm}1\leq i\leq 3$$  and further we also have  
 $$\hbox{\rm Sing}(\overline{\hbox{\rm Quad}_{V}^{(i+1)}})=\overline{\hbox{\rm Quad}_{V}^{(i+1)}}\hspace*{2mm}\backslash\hspace*{2mm}
\hbox{\rm Quad}_{V}^{(i+1)}\hspace*{0.25cm}for\hspace*{0.25cm}1\leq i\leq 2$$
 unless the characteristic of 
 $k$ is 2 in which case $\overline{\hbox{\rm Quad}_{V}^{(3)}}$ is itself
 smooth  
 (the notation $\overline{T}$
 denotes the orbit closure and $\hbox{\rm Sing}(T)$ denotes the subset of 
 singular (non-smooth) points of $T$, each given the canonical reduced 
 induced closed subscheme structure).  
\item[(3)] Continuing with the notations and hypotheses of (2) above, 
 set $$(W, w):=(\Lambda^{even}(V), 1).$$ For 
 ease of notation denote $\spazuWw$ by $\hbox{\rm SpAzu}$ and 
 $\hbox{\rm Stab}_{w}$ by $H.$
 Then the four orbits  
$$\hbox{\rm SpAzu}^{(i)}:=H\cdot\Theta(q^{(i)}, I_4)\textrm{ for }
 1\leq i\leq 4$$ form a stratification of the $k$-variety $\hbox{\rm SpAzu}$ 
 in the sense that we have 
 $$\overline{\hbox{\rm SpAzu}^{(1)}}=\hbox{\rm SpAzu}\hspace*{0.25cm}and\hspace*{0.25cm}
  \overline{\hbox{\rm SpAzu}^{(i+1)}}=\overline{\hbox{\rm SpAzu}^{(i)}}\hspace*{2mm}\backslash\hspace*{2mm}
\hbox{\rm SpAzu}^{(i)},\hspace*{0.25cm}1\leq i\leq 3$$  and further we also have  
 $$\hbox{\rm Sing}(\overline{\hbox{\rm SpAzu}^{(i+1)}})=\overline{\hbox{\rm SpAzu}^{(i+1)}}\hspace*{2mm}\backslash\hspace*{2mm}
\hbox{\rm SpAzu}^{(i+1)}\hspace*{0.25cm}for\hspace*{0.25cm}1\leq i\leq 2$$
 unless the characteristic of 
 $k$ is 2 in which case $\overline{\hbox{\rm SpAzu}^{(3)}}$ is itself smooth. 
\end{description}
\end{theorem}

\section{Injectivity: Reduction to Lifting to Similarities in the Free
Case} 
\label{sec3}
\paragraph*{\sc  Proof of Prop.\ref{transfer-to-lambda2}.}
Start with an isomorphism of algebra-bundles $$\phi:\covqi\cong\covprqpripr.$$
Let ${\{U_i\}}_{i\in \mathcal{I}}$
 be an affine open covering of $X$ (which may 
also be chosen so as to trivialise some or any of 
 the involved bundles if needed). Choose bilinear forms 
 $$b_i\in\Gamma(U_i, \hbox{\rm Bil}_{(V, I)})\textrm{ and }
b'_i\in\Gamma(U_i, \hbox{\rm Bil}_{(V', I')})$$ 
such that 
$$q|U_i=q_{b_i}\textrm{ and }q'|U_i=q_{b'_i}\textrm{ for each }
i\in \mathcal{I}.$$ 
 By (2d), Theorem \ref{Bourbaki}, we have  isomorphisms of vector bundles 
$\psi_{b_i}$ and 
 $\psi_{b'_i}$, which preserve 1 by (2a) of the same Theorem, 
  and we define the isomorphism of vector bundles ${(\phi_{\Lambda^{ev}})}_i$ 
 so as to make the following diagram commute: 
$$\begin{CD}
C_0(V,q, I)|U_i @>{\phi|U_i}>{\cong}> C_0(V', q', I')|U_i\\
@V{\psi_{b_i}}V{\cong}V @V{\cong}V{\psi_{b'_i}}V\\
({\cal O}_X\oplus\Lambda^2(V)\otimes I^{-1})|U_i 
@>{\cong}>{{(\phi_{\Lambda^{ev}})}_i}> 
({\cal O}_X\oplus\Lambda^2(V')\otimes {(I')}^{-1})|U_i 
\end{CD}$$
The linear isomorphism ${(\phi_{\Lambda^{ev}})}_i$ preserves 1  and 
 therefore it induces a linear isomorphism 
$${({\phi_{\Lambda^2}})}_i:(\Lambda^2(V)\otimes I^{-1})|U_i
\stackrel{\cong}{\longrightarrow} 
 (\Lambda^2(V')\otimes {(I')}^{-1})|U_i.$$
Observe that ${({\phi_{\Lambda^2}})}_i$ is independent of the choice of 
 the bilinear forms $b_i$ and $b'_i.$ For, replacing these respectively 
 by $\widehat{b_i}$ and $\widehat{b'_i}$, it follows from (2f),
 Theorem \ref{Bourbaki}, that 
$$\psi_{b_i}\circ{(\psi_{\widehat{b_i}})}^{-1}\textrm{ 
 (resp. }\psi_{b'_i}\circ{(\psi_{\widehat{b'_i}})}^{-1})$$ followed by the 
 canonical projection onto 
$$(\Lambda^2(V)\otimes I^{-1})|U_i\textrm{ 
(resp. onto }(\Lambda^2(V')\otimes {(I')}^{-1})|U_i)$$ 
 is the same as the projection itself. By this observation, it is also 
 clear that the isomorphisms 
${\{{({\phi_{\Lambda^2}})}_i\}}_{i\in \mathcal{I}}$
 agree on (any open affine subscheme of, and hence on all of)
 any intersection $U_i\cap U_j.$ Therefore they glue to give a global 
 isomorphism of vector bundles 
$$\phi_{\Lambda^2}:\Lambda^2(V)\otimes I^{-1}\cong
\Lambda^2(V')\otimes {(I')}^{-1}$$
as required.  
{\bf Q.E.D, Prop.\ref{transfer-to-lambda2}.}
\paragraph*{\sc  Reduction of Proof of Injectivity of Theorem
 \ref{bijectivity} to 
 Theorem \ref{lifting-of-isomorphisms}.}
 We start with an isomorphism of algebra-bundles
 $$\phi:\covqi\cong\covprqpripr,$$ construct the isomorphism of vector
 bundles 
$$\phi_{\Lambda^2}:\Lambda^2(V)\otimes I^{-1}
\cong\Lambda^2(V')\otimes {(I')}^{-1}$$ 
 and keep the notations
 introduced in the proof of Prop.\ref{transfer-to-lambda2}. 
 Firstly we deduce a linear isomorphism 
  $$\hbox{\rm det}{({(\phi_{\Lambda^2})}^\vee)}^{-1}:
\hbox{\rm det}({(\Lambda^2(V)\otimes I^{-1})}^\vee)\cong
\hbox{\rm det}({(\Lambda^2(V')\otimes {(I')}^{-1})}^\vee).$$ 
 Since $V$ and  $V'$ are of rank 3,  
  there are canonical isomorphisms
 $$\eta:\Lambda^2(V)\equiv V^\vee\otimes\hbox{\rm det}(V)\textrm{ and }
 \eta':\Lambda^2(V')\equiv (V')^\vee\otimes\hbox{\rm det}(V').$$
It follows therefore that if we set 
$$L:=\hbox{\rm det}(V')
\otimes{(\hbox{\rm det}(V))}^{-1}\textrm{ and }J:=I'\otimes I^{-1}$$ 
 then we get a twisted discriminant line bundle $(L\otimes J^{-1},h, J)$
 and a vector 
 bundle isomorphism 
$$\alpha:V'\cong V\otimes (L\otimes J^{-1}).$$
 Now for each $i\in \mathcal{I}$, 
 the bilinear form 
 $b_i\in\Gamma(U_i, \hbox{\rm Bil}_{(V,I)})$
 induces, via $\alpha|U_i$ and $(LJ^{-1},h, J)|U_i$ and
 (1), Prop.\ref{tensoring-with-disc-module}, a bilinear form 
 $$b''_i\in\Gamma(U_i, \hbox{\rm Bil}_{(V', IJ)}).$$
 By (3) of the same Proposition, 
 over each $U_i$  we get an isometry  of bilinear form bundles 
$$\alpha|U_i:(V'|U_i, b''_i, IJ|U_i)\cong (V|U_i, b_i, I|U_i)
\otimes (LJ^{-1}, h, J)|U_i$$
 and also an isometry of quadratic bundles
$$\alpha|U_i:(V'|U_i, q_{b''_i}, IJ|U_i)\cong 
(V|U_i, q_{b_i}=q|U_i, I|U_i)\otimes (LJ^{-1}, h, J)|U_i.$$
 On the other hand,
 by an assertion in (3), Prop.\ref{tensoring-with-disc-module}, we could also 
 define the global
 quadratic bundle $(V',q'', IJ)$ using $(V,q, I)$, $\alpha$ and
 $(LJ^{-1},h, J)$, 
 so that we have an isometry of quadratic bundles 
$$\alpha:(V',q'', IJ)\cong(V,q, I)\otimes(LJ^{-1},h, J).$$
 It follows therefore that 
 the $q_{b''_i}$ glue to give $q''.$ 
Notice that in general 
 the $b''_i$ (resp. the $b_i$) need not glue to give a global bilinear 
 form $b''$ (resp. $b$) such that $q_{b''}=q''$ (resp. $q_b=q$). 
 By (1), Prop.\ref{simil-induces-iso-of-even-cliff}, there exists 
 a  unique isomorphism of algebra bundles
$$C_0(\alpha,1,IJ):C_0(V',q'', IJ)\cong C_0\left((V,q, I)
\otimes(LJ^{-1},h, J)\right)$$ and by 
 Prop.\ref{isom-covq-covlqh} 
 we have a unique isomorphism of algebra bundles 
 $$\gamma_{(LJ^{-1},h,J)}:C_0\left(\vqi\otimes(LJ^{-1},h,J)\right)
\cong\covqi.$$
 Therefore the composition 
 of the following sequence of isomorphisms of algebra bundles on X
\begin{center}\begin{tabular}{c}
$C_0(V',q'', I')\stackrel{\hbox{\rm\tiny using }I'\cong IJ}{\cong}
C_0(V',q'', IJ)\stackrel{C_0(\alpha,1,IJ)\thinspace{(\cong)}}{\longrightarrow}
C_0\left((V,q,I)\otimes(LJ^{-1},h,J)\right)$\\
$\stackrel{{\gamma_{(LJ^{-1},h,J)}}{\thinspace{(\cong)}}}{\longrightarrow}
\covqi\stackrel{\phi\thinspace{(\cong)}}{\longrightarrow}\covprqpripr$
\end{tabular}\end{center}
is an element of $$\hbox{\rm Iso}[C_0(V',q'', I'),C_0(V',q',I')],$$
 which, granting 
 Theorem \ref{lifting-of-isomorphisms}, is induced by a similarity 
 $$(g,l)\in\hbox{\rm Sim}[(V',q'', I'),(V',q',I')].$$ Hence we would have 
 $$g:(V',q'',I')\cong(V',q',I')\otimes({\cal O}_X,(s\otimes s'\mapsto 
 s.s'.l^{-1}), {\cal O}_X)$$ where $l\in\Gamma(X, {\cal O}_X^*).$
 This combined with the fact that 
   $(V,q,I)$ and 
$$(V',q'',I')=(V',q'',IJ)$$ are isomorphic 
 upto the twisted discriminant bundle $(LJ^{-1},h, J)$
 by the construction above, would imply 
 that $(V,q,I)$ and $(V',q',I')$ also differ by a
twisted  discriminant 
 bundle. Therefore the proof of the injectivity asserted in  
 Theorem \ref{bijectivity} reduces 
 to the proof of 
 Theorem \ref{lifting-of-isomorphisms}. 
\paragraph*{\sc  Reduction of Theorem \ref{lifting-of-isomorphisms}
 to the Case when $I$ is free.}
 For a similarity $g$ with multiplier $l$, we have $C_0(g,l,I)$
 given by (1), Prop.\ref{simil-induces-iso-of-even-cliff}, so that we 
 may define the map $$\simvqvqpri \longrightarrow \isomcovqcovqpri:
g\mapsto C_0(g,l,I).$$   
 The equality $$\hbox{\rm det}(C_0(g,l, I))
=\hbox{\rm det}\left((C_0(g,l, I))_{\Lambda^2}\right)
=l^{-3}\hbox{\rm det}^2(g)$$ will be shown to hold (locally, hence globally) 
 in (1), Lemma \ref{locally-det-coglambda-equals-detgbylcube}, 
 page \pageref{locally-det-coglambda-equals-detgbylcube}. 
 Thus $$\isomvqvqpri\textrm{ and }\sisomvqvqpri$$
 are respectively mapped 
 into  
 $$\isomcovqcovqprprimei\textrm{ and }\sisomcovqcovqpri$$ as claimed.
  We start with an isomorphism of algebra-bundles 
$$\phi:\covqi\cong\covqpri,$$
 which by Prop.\ref{transfer-to-lambda2} leads to the automorphism of vector
 bundles 
$$\phi_{\Lambda^2}\in\hbox{\rm Aut}(\Lambda^2(V)\otimes I^{-1}).$$
Firstly, define the global bundle automorphism 
 $$g'\in\hbox{\rm GL}\left(V\otimes{(\hbox{\rm det}(V))}^{-1}\otimes I\right)$$
 so that the following diagram commutes
$$\begin{CD}
{(\Lambda^2(V))}^\vee\otimes I 
@>{{({({\phi_{\Lambda^2}})}^\vee)}^{-1}}>{\cong}> 
{(\Lambda^2(V))}^\vee\otimes I\\
@V{{(\eta^\vee)}^{-1}\otimes I}V{\equiv}V 
@V{\equiv}V{{(\eta^\vee)}^{-1}\otimes I}V\\
V\otimes{(\hbox{\rm det}(V))}^{-1}\otimes I
 @>{\cong}>g'> V\otimes{(\hbox{\rm det}(V))}^{-1}\otimes I
\end{CD}$$ where 
 $\eta:\Lambda^2(V)\equiv V^\vee\otimes\hbox{\rm det}(V)$
 is the canonical isomorphism (since $V$ is of rank 3). Now let 
 $$g\in\hbox{\rm GL}(V)\stackrel{\cong}{\longleftarrow}\hbox{\rm GL}
(V\otimes{(\hbox{\rm det}(V))}^{-1}\otimes I)$$ be the image of $g'$ i.e., 
 the image of $g'\otimes \hbox{\rm det}(V)\otimes I^{-1}$
 under the canonical identification 
 $$\hbox{\rm GL}(V\otimes{(\hbox{\rm det}(V))}^{-1}\otimes I\otimes\hbox{\rm det}(V)\otimes I^{-1})
\equiv \hbox{\rm GL}(V).$$ 
 Next, let $l\in\Gamma(X, {\cal O}_X^*)$ be a global section such that 
 $$\gamma(l):=(l^3).\hbox{\rm det}(\phi_{\Lambda^2})$$ has a square root in 
 $\Gamma(X, {\cal O}_X^*).$ For example, we have the following 
independent special 
 cases when this is true: 
\begin{description}
\item[Case 1.]  \label{cases} If $\phi\in\sisomcovqcovqpri$ i.e., if 
 $\hbox{\rm det}(\phi_{\Lambda^2})=1,$ then  
set $l=1$ and $\sqrt{\gamma(l)}=1.$
\item[Case 2.]  If $\phi\in\isomcovqcovqprprimei$ i.e., 
$\hbox{\rm det}(\phi_{\Lambda^2})$ is a square, then set 
 $l=1$ and  let $\sqrt{\gamma(l)}$ denote any 
 fixed square root of $\hbox{\rm det}(\phi_{\Lambda^2}).$
\item[Case 3.] Given an integer $k$,
 take $l={(\hbox{\rm det}(\phi_{\Lambda^2}))}^{2k+1}$ and
 let $\sqrt{\gamma(l)}$ denote any fixed square root of 
 ${(\hbox{\rm det}(\phi_{\Lambda^2}))}^{6k+4}.$
\end{description}
 For each integer 
 $k$, we now associate to $\phi$ the element $$\label{glphi-def}
g_l^\phi:=(l^{-1}\sqrt{\gamma(l)}\thinspace)g$$ with 
 $g$ as defined above. We shall   
 show the following locally for the Zariski topology on $X$  
 (more precisely, for each open subscheme of $X$ over which 
 $V$ and $I$ are free): 
\begin{description}
\item[(1)]  
 that $g_l^\phi$ is an $I$-similarity from $(V,q,I)$ to $(V,q',I)$
 with multiplier $l$ (Lemma
 \ref{glphi-simil-q-qpr-with-mult-l}, page 
\pageref{glphi-simil-q-qpr-with-mult-l});  
\item[(2)] that $g_l^\phi$  
 induces $\phi$ i.e., with the notations of (1), 
Prop.\ref{simil-induces-iso-of-even-cliff}, 
 that $C_0(g_l^\phi, l, I)=\phi$ (Lemma 
 \ref{glphi-induces-phi}, page \pageref{glphi-induces-phi});  
\item[(3)] that $\hbox{\rm det}(g_l^\phi)=\sqrt{\gamma(l)}$
 so that 
 $\hbox{\rm det}^2(g_l^\phi)
=\hbox{\rm det}(\phi_{\Lambda^2})$ in cases 1 and 2 
(Lemma \ref{B-versus-Bglphi}, page \pageref{B-versus-Bglphi}) and  
\item[(4)]
 that the  map 
$$\sisomvqvqpri\longrightarrow\sisomcovqcovqpri$$
 is injective (Lemma \ref{bijection-of-special-isomorphisms}, page 
 \pageref{bijection-of-special-isomorphisms}).
\end{description}
 It would follow then that these statements are also true globally. The 
 maps 
$$s_{2k+1}:\phi\mapsto g_l^\phi\textrm{ with $l$ as in Case 3}$$ and 
 $$s': \phi\mapsto g_l^\phi\textrm{ with $l$ as in Case 2}$$
 will then give the 
 sections to the maps (which would imply their surjectivities) 
as mentioned in Theorem \ref{lifting-of-isomorphisms}.
 But these maps 
  are  not necessarily multiplicative since a computation reveals that 
 if $$\phi_i\in\hbox{\rm Iso}[C_0(V,q_i,I),C_0(V,q_{i+1},I)]$$ 
is associated to $$g_{l_i}^{\phi_i}\in\hbox{\rm Sim}
[(V,q_i,I),(V,q_{i+1},I)],$$ 
 and $\phi_2\circ\phi_1$ to $g_{l_{21}}^{\phi_2\circ\phi_1}$, then 
 $$g_{l_{21}}^{\phi_2\circ\phi_1}=\delta g_{l_2}^{\phi_2}\circ 
g_{l_1}^{\phi_1}\textrm{ for }
\delta\in\mu_2(\Gamma(X, {\cal O}_X))$$ because of the ambiguity in 
 the initial global choices of square roots for 
 $\gamma(l_i)$ and $\gamma(l_{21}).$  However 
 this can be remedied as follows. For any given 
 $$\phi\in\hbox{\rm Iso}[C_0(V,q,I),C_0(V,q',I)],$$ 
 irrespective of whether or not
 $\hbox{\rm det}(\phi_{\Lambda^2})$ is a square,
 take 
$$l={(\hbox{\rm det}(\phi_{\Lambda^2}))}^{2k+1}, 
 \gamma(l)=l^3\hbox{\rm det}(\phi_{\Lambda^2}),
  \sqrt{\gamma(l)}:={(\hbox{\rm det}(\phi_{\Lambda^2}))}^{3k+2}$$ and 
 $$s_{2k+1}^{+}(\phi):=g_l^\phi=\left(l^{-1}\sqrt{\gamma(l)}\right)g.$$
 Then it is clear that each $s_{2k+1}^{+}$ is multiplicative with the 
 properties as claimed in the statement.
 We thus reduce the proof of Theorem \ref{lifting-of-isomorphisms} 
 to the case when the rank 3 vector bundle 
 $V$ and the line bundle $I$ are free. \label{reduction-to-free-case}
 This will be treated in the next section. 

\section{The Free Case: Investigation of the Isomorphism Theta}\label{sec4}
Throughout this section, we work with $I={\cal O}_X$ and shorten our 
 earlier notations $\vqi$, $\covqi$, $C_0(g,l, I)$ etc respectively to
 $\vq$, $\covq$, $C_0(g,l)$ etc. 
We conclude the proofs of the injectivity of Theorem \ref{bijectivity} and
Theorem  \ref{lifting-of-isomorphisms}
 which were begun in \S \ref{sec3}  and  also 
 prove Theorem \ref{lifting-of-automorphisms}.

 As means to these ends, 
 we carry out two explicit computations. Firstly we compute 
the isomorphism $\Theta$ of Theorem \ref{Thetaisom}. This provides us with
 the multiplication table of every specialised algebra structure on any fixed
  free rank 4 vector bundle with fixed unit which is part of a global basis (Theorem \ref{mult-table-wrt-Theta} below). 
 This result  will also be used in \S \ref{sec5}   
 in the proof of Theorem \ref{Thetaisom-Upsilonisom}.
It turns out that $\Theta$ is not equivariant 
 with respect to $\hbox{\rm GL}_{V}$, but nevertheless satisfies a `twisted' 
 form of equivariance (Theorem \ref{tequiv}). 
 Secondly, we explicitly  compute the algebra bundle 
 isomorphism 
$$C_0(g,l):C_0(V,q)\cong C_0(V,q')$$ of 
 (1), Prop.\ref{simil-induces-iso-of-even-cliff} induced by a similarity $g:\vq\thinspace{\cong}_{l}\thinspace
\vqpr$
with multiplier $l\in\Gamma(X, {\cal O}_X^*)$ in the case when 
 $V$ is free of rank 3 (Theorem \ref{b-co-phi-lambda-indep-of-q}).
\subsection{The Action of GL on Forms}
Let $V$ be a vector bundle over a scheme $X$ with associated locally-free 
 sheaf $\cal V.$ The $X$-smooth $X$-groupscheme $\hbox{\rm GL}_{V}$ acts
 naturally on the left on  
 the sheaves $\hbox{\rm Alt}^2_{\cal V}$, $\hbox{\rm Bil}_{\cal V}$ and 
 $\hbox{\rm Quad}_{\cal V}$ of alternating, bilinear and quadratic forms on 
 $\cal V$ (with values in ${\cal O}_X$).
 Namely, for $U\hookrightarrow X$ an open subscheme, and for 
 $$b\in\Gamma(U, \hbox{\rm Bil}_{\cal V})\textrm{ (resp. }
a\in\Gamma(U, \hbox{\rm Alt}^2_{\cal V}),\textrm{  resp. }
q\in\Gamma(U, \hbox{\rm Quad}_{\cal V})),$$ and for 
 $g\in\Gamma(U,\hbox{\rm GL}_{V})=\hbox{\rm GL}(V|U)$, the corresponding 
 form of the same type $g.b$ (resp. $g.a$, resp. $g.q$) is defined on 
 sections (over open subsets of $U$) by 
\begin{equation}\begin{split}
(g.b)(v,v'):=b(g^{-1}(v), g^{-1}(v'))&\textrm{ (resp. }
(g.a)(v,v'):=a(g^{-1}(v), g^{-1}(v')),\notag \\
 &\textrm{  resp. }
 (g.q)(v):=q(g^{-1}(v))).\notag \end{split}\end{equation} 
 It is immediate that the following short-exact-sequence of sheaves, 
 indicated in \S\ref{subsec2.1},
 is  equivariant with respect to this action:
$$\label{spadesuit}(\spadesuit)\hspace*{8mm} 0\longrightarrow \hbox{\rm Alt}^2_{\cal V} \longrightarrow \hbox{\rm Bil}_{\cal V} \longrightarrow \hbox{\rm Quad}_{\cal V}\longrightarrow 0.$$
Equivalently, the $X$-groupscheme $\hbox{\rm GL}_{V}$  acts on the 
 corresponding geometric vector bundles such that  both of the $X$-morphisms 
 of $X$-vector bundles in the  following sequence are
  $\hbox{\rm GL}_{V}$-equivariant: 
 $$\hbox{\rm Alt}^2_{V} \hookrightarrow \hbox{\rm Bil}_{V} \twoheadrightarrow \hbox{\rm Quad}_{V}.$$ Notice that it is one and the same thing to require that 
  $$\hbox{\rm GL}(V|U)\ni g:(V|U, q)\thinspace{\cong}_{l}
\thinspace(V|U, q')$$ be a 
 similitude with multiplier $l\in\Gamma(U, {\cal O}_X)$, and to 
 require that $g.q={l}^{-1}q'.$ 
\subsection{Computation of the Isomorphism Theta}\label{Thetadef}
We briefly recall the definition of $\Theta$
 from Part A of \cite{tevb-paper1}. We keep the notations introduced just 
 before Theorem \ref{Thetaisom}; for ease of notation, 
 the pullback of a section $s$ (of a vector bundle or its associated sheaf)
 is denoted by $s^\circ.$  Since $V$ is free of rank 3 on $X$, we 
  choose an identification 
$${\cal V}\equiv
{\cal O}_X.e_1\oplus{\cal O}_X.e_2\oplus{\cal O}_X.e_3.$$
This gives the identification of the dual bundle as 
$${\cal V}^\vee\equiv
{\cal O}_X.f_1\oplus{\cal O}_X.f_2\oplus{\cal O}_X.f_3$$ (defined uniquely 
 by $f_i(e_j)=\delta_{ij}$, the Kronecker delta). Therefore the dual of 
 the sheaf 
 of quadratic forms on $V$, which is 
$${(\hbox{\rm Quad}_{\cal V})}^\vee:={(\hbox{\rm Bil}_{\cal V}/\hbox{\rm Alt}^2_{\cal V})}^\vee={({(T^2{\cal V})}^\vee/{(\Lambda^2{\cal V})}^\vee)}^\vee,$$
 has global ${\cal O}_X$-basis
given by 
$$\{e_i\otimes e_i; (e_i\otimes e_j+e_j\otimes e_i)\}.$$
 This leads to 
 an identification of the associated sheaf of 
 symmetric algebras 
$$\hbox{\rm Sym}_{{\cal O}_X}(\hbox{\rm Quad}_{\cal V}^\vee)\equiv
 {\cal O}_X[Y_1, Y_2, Y_3, Y_{12}, Y_{13}, Y_{23}],$$ where 
 $e_i\otimes e_i\equiv Y_i$ and $e_i\otimes e_j+e_j\otimes e_i\equiv Y_{ij}$,
 and therefore 
$$\hbox{\rm Quad}_{V}:=
 \hbox{\bf Spec}\left(\hbox{\rm Sym}_{{\cal O}_X}(\hbox{\rm Quad}_{\cal V}^\vee)\right)\equiv
 \hbox{\bf Spec}({\cal O}_X[Y_1, Y_2, Y_3, Y_{12}, Y_{13}, Y_{23}])={\mathbb A}^6_X.$$  Consider the 
 universal quadratic bundle  $({\mathbf V},{\mathbf q})$ where 
 $\mathbf V$ is the pullback of $V$ by $\hbox{\rm Quad}_{V}\longrightarrow X.$
 The universal quadratic form $\mathbf q$ is given by 
 $$(x_1, x_2, x_3)\mapsto \Sigma_i Y_i.(x_i)^2 + \Sigma_{i<j} Y_{ij}.x_i.x_j$$
 and moreover the global bilinear form on $\mathbf V$ given by 
 $${\mathbf b}({\mathbf q}): \left((x_1, x_2, x_3), (x'_1, x'_2, x'_3)\right)\mapsto
  \Sigma_i Y_i.x_i.x'_i + Y_{12}.x_2.x'_1+Y_{23}.x_3.x'_2+Y_{13}.x_1.x'_3$$
 induces 
 $\mathbf q$ (the bilinear form
 `associated in the usual sense' to $\mathbf q$, viz. $b_{\mathbf q}$ is 
 not $\mathbf b({\mathbf q})$ but in fact its symmetrisation). 
  Therefore, by (2d), Theorem \ref{Bourbaki}, we get an 
 isomorphism of vector bundles 
$$\psi_{{\mathbf b}({\mathbf q})}:C_0({\mathbf V}, {\mathbf
 q}=q_{{\mathbf b}({\mathbf q})})\cong \Lambda^{even}({\mathbf
 V})=:{\mathbf W}$$
 which, according 
 to (2a) and (2f) of the same Theorem,  carries  the ordered
 Poincar\'e-Birkhoff-Witt basis 
$$\{1;e^\circ_1.e^\circ_2, e^\circ_2.e^\circ_3, e^\circ_3.e^\circ_1\}$$  onto 
 the corresponding ordered basis of the even exterior algebra (=even 
 Clifford algebra of the zero quadratic form on $\mathbf V$) given by 
 $$\{w^\circ=1^\circ=1;e^\circ_1\wedge e^\circ_2, 
e^\circ_2\wedge e^\circ_3, e^\circ_3\wedge e^\circ_1\}.$$ 
The choices $e^\circ_3.e^\circ_1$ and $e^\circ_3\wedge e^\circ_1$ 
instead of the usual $e^\circ_1.e^\circ_3$ and $e^\circ_1\wedge
e^\circ_3$ are deliberate---for example, 
$\psi_{{\mathbf b}({\mathbf q})}$ would 
 carry 
$$\{1;e^\circ_1.e^\circ_2, e^\circ_2.e^\circ_3, e^\circ_1.e^\circ_3\}$$
 onto  
$$\{w^\circ=1^\circ=1;e^\circ_1\wedge e^\circ_2, e^\circ_2\wedge e^\circ_3, e^\circ_1\wedge e^\circ_3 + Y_{13}.w^\circ\}$$ which depends 
 on $Y_{13}.$
 Thus 
 the even Clifford algebra bundle $C_0({\mathbf V},
 {\mathbf q}=q_{{\mathbf b}({\mathbf q})})$
 induces  via $\psi_{{\mathbf b}({\mathbf q})}$ a
 $w^\circ$-unital algebra structure on the pullback bundle $\mathbf W$
 of $W:=\Lambda^{even}(V)$ (where $w$ corresponds to 1
 in ${\Lambda^{even}(V)}$).
 But by definition, this algebra structure corresponds precisely
 to an $X$-morphism \label{thetadef}  
$$\theta:\hbox{\rm Quad}_{V}
\longrightarrow \hbox{\rm Id-}w\hbox{\rm -Sp-Azu}_W.$$ The isomorphism 
 $\Theta$ is now given by the composition of the following $X$-morphisms (cf.
Def.5.2, Part A, \cite{tevb-paper1}):
$$\hbox{\rm Quad}_{V}\times_X\hbox{\sf
  L}_w\stackrel{\theta\times\hbox{\rm\tiny
    ID}}{\longrightarrow}\hbox{\rm Id-}w\hbox{\rm
  -Sp-Azu}_W\times_X\hbox{\sf L}_w\stackrel{\hbox{\rm\tiny
    SWAP}\thinspace(\equiv)}{\longrightarrow}$$
$$\hbox{\sf L}_w\times_X\hbox{\rm Id-}w\hbox{\rm -Sp-Azu}_W\stackrel{\hbox{\rm\tiny ACTION}}{\longrightarrow}\hbox{\rm Id-}w\hbox{\rm -Sp-Azu}_W.$$ 
 The association of $\mathbf q$ with ${\mathbf b}({\mathbf q})$
 also defines a splitting
 of the exact sequence  $(\spadesuit)$  of page \pageref{spadesuit} above,
 so that more generally,
 given a valued point $q\in(\hbox{\rm Quad}_{V})(T)$, we may 
 associate uniquely a valued point $b(q)\in(\hbox{\rm Bil}_{V})(T)$
 which induces it. That this association is not $\hbox{\rm GL}_V$-equivariant 
 is reflected in the lack of equivariance of the isomorphism $\Theta$
 (Theorem \ref{tequiv}). 
\begin{theorem}\label{mult-table-wrt-Theta} Let $T$ be an $X$-scheme. Let $q$
 be
a $T$-valued point of $\hbox{\rm Quad}_{V}\equiv{\mathbb A}^6_X$ which is
identified uniquely with a 6-tuple $$(\lambda_1,\lambda_2,\lambda_3,
\lambda_{12},\lambda_{13},\lambda_{23})$$ corresponding to the 
 quadratic form $$(x_1,x_2,x_3)\mapsto \Sigma_i \lambda_i x_i^2 + \Sigma_{i<j}\lambda_{ij}x_ix_j.$$ Let $\underline{t}$ be a  $T$-valued point  
  of $\hbox{\sf L}_w\equiv({\mathbb A}^3_X, +)$ which 
 is identified uniquely with a 3-tuple $(t_1, t_2, t_3)$ that 
 corresponds to the $T$-valued point of $\hbox{\rm Stab}_{w}$ given by the 
 $(4\times 4)$-matrix
 \begin{displaymath}
 \left(\begin{array}{c|ccc}
1 & t_1 & t_2 & t_3 \\
\hline
0 & & I_3 & 
\end{array}
\right)
 \end{displaymath} 
where $I_3$ is the $(3\times 3)$-identity matrix.
 Then in terms of the global basis 
 $$\{w^\circ=1^\circ=1\hspace*{1mm};\hspace*{1mm}
\epsilon^\circ_1:=e^\circ_1\wedge e^\circ_2\hspace*{1mm},\hspace*{1mm}
\epsilon^\circ_2:=e^\circ_2\wedge e^\circ_3\hspace*{1mm},\hspace*{1mm}
\epsilon^\circ_3:=e^\circ_3\wedge e^\circ_1\}$$ induced from that 
of $W=\Lambda^{even}(V),$ the multiplication 
 table for the specialised algebra structure $$\Theta(q,\underline{t})=
\underline{t}.\theta(q)$$ on the pullback bundle $W_T$ with 
 unit $w^\circ=w_T$ is given as follows:   
\begin{enumerate}
\item[$\bullet$] $\epsilon^\circ_1.\epsilon^\circ_1= (t_1\lambda_{12}-\lambda_1\lambda_2-t_1^2).w^\circ+(\lambda_{12}-2t_1).\epsilon^\circ_1\thinspace;$
\item[$\bullet$] $\epsilon^\circ_2.\epsilon^\circ_2= (t_2\lambda_{23}-\lambda_2\lambda_3-t_2^2).w^\circ+(\lambda_{23}-2t_2).\epsilon^\circ_2\thinspace;$
\item[$\bullet$] $\epsilon^\circ_3.\epsilon^\circ_3= (t_3\lambda_{13}-\lambda_1\lambda_3-t_3^2).w^\circ+(\lambda_{13}-2t_3).\epsilon^\circ_3\thinspace;$
\item[$\bullet$] $\epsilon^\circ_1.\epsilon^\circ_2=(\lambda_2\lambda_{13}-\lambda_2t_3-t_1t_2).w^\circ -t_2\epsilon^\circ_1-t_1\epsilon^\circ_2-\lambda_2\epsilon^\circ_3
\thinspace;$
\item[$\bullet$] $\epsilon^\circ_2.\epsilon^\circ_3=(\lambda_3\lambda_{12}-\lambda_3t_1-t_2t_3).w^\circ -\lambda_3\epsilon^\circ_1-t_3\epsilon^\circ_2-t_2\epsilon^\circ_3
\thinspace;$
\item[$\bullet$] $\epsilon^\circ_3.\epsilon^\circ_1=(\lambda_1\lambda_{23}-\lambda_1t_2-t_1t_3).w^\circ -t_3\epsilon^\circ_1-\lambda_1\epsilon^\circ_2-t_1\epsilon^\circ_3
\thinspace;$
\item[$\bullet$] $\epsilon^\circ_2.\epsilon^\circ_1=(\lambda_2t_3-(\lambda_{12}-t_1)(\lambda_{23}-t_2)).w^\circ +(\lambda_{23}-t_2)\epsilon^\circ_1+(\lambda_{12}-t_1)\epsilon^\circ_2+\lambda_2\epsilon^\circ_3\thinspace;$
\item[$\bullet$] $\epsilon^\circ_3.\epsilon^\circ_2=(\lambda_3t_1-(\lambda_{13}-t_3)(\lambda_{23}-t_2)).w^\circ +\lambda_3\epsilon^\circ_1+(\lambda_{13}-t_3)\epsilon^\circ_2+(\lambda_{23}-t_2)\epsilon^\circ_3\thinspace;$
\item[$\bullet$] $\epsilon^\circ_1.\epsilon^\circ_3=(\lambda_1t_2-(\lambda_{12}-t_1)(\lambda_{13}-t_3)).w^\circ +(\lambda_{13}-t_3)\epsilon^\circ_1+\lambda_1\epsilon^\circ_2+(\lambda_{12}-t_1)\epsilon^\circ_3$.
\end{enumerate}
\end{theorem}
\paragraph*{\sc Proof of Theorem \ref{mult-table-wrt-Theta}.} 
 For clarity, let $*_q$ denote the multiplication in the algebra $C_0(V_T,q)$, 
 and for uniformity, let $\epsilon_0:=w.$  Since $q=q_{b(q)}$, we have by 
 (2d), Theorem \ref{Bourbaki}, 
 the isomorphism $$\psi_{b(q)}:C_0(V_T,q)\cong \Lambda^{even}(V_T)=W_T.$$ 
 Let $*_{b(q)}$ denote the 
 product in the algebra structure $\theta(q)$ thus induced on $W_T.$  
 Since 
 the $\epsilon^\circ_i$ are a basis
for $W_T$,    
 it is  enough to compute the products  
$\epsilon^\circ_i*_{b(q)} \epsilon^\circ_j$ for $1\leq i,j \leq 3.$
 For example, consider 
the product $\epsilon^\circ_2*_{b(q)} \epsilon^\circ_1.$ 
Using the properties of the multiplication in $C(V_T, q)$, and the 
 properties of the isomorphism $\psi_{b(q)}$ from (2), Theorem \ref{Bourbaki},
 we get the following: 
\begin{equation}
\begin{split}
\epsilon^\circ_2*_{b(q)} \epsilon^\circ_1 
& =\psi_{b(q)}\left(
\{ \psi_{b(q)}^{-1}( e_2^{\circ}\wedge e_3^{\circ} ) \} *_q
\{ \psi_{b(q)}^{-1}( e_1^{\circ}\wedge e_2^{\circ} ) \} \right) \notag \\
& =\psi_{b(q)}\left(
  (e_2^{\circ}*_q e_3^{\circ})   *_q
  (e_1^{\circ}*_q e_2^{\circ})   \right) \notag \\
& =\psi_{b(q)}\left(
 (\lambda_{23}(1^{\circ})-e_3^{\circ}*_q e_2^{\circ})*_q
 (\lambda_{12}(1^{\circ})-e_2^{\circ}*_q e_1^{\circ}) \right) \\ \notag 
& =\psi_{b(q)}\left(
\lambda_{23}\lambda_{12}(1^{\circ})-\lambda_{23}e_2^{\circ}*_q
 e_1^{\circ}-\lambda_{12}e_3^{\circ}*_q e_2^{\circ}
  + (e_3^{\circ}*_q
 e_2^{\circ})*_q(e_2^{\circ}*_q
 e_1^{\circ})\right) \\ \notag 
 & =\psi_{b(q)}\left(
\lambda_{23}\lambda_{12}(1^{\circ})-\lambda_{23}(\lambda_{12}(1^{\circ})-e_1^{\circ}*_q
 e_2^{\circ})\right. \\ \notag
& \quad 
\left.-\lambda_{12}(\lambda_{23}(1^{\circ})-e_2^{\circ}*_q e_3^{\circ}) 
+ e_3^{\circ}*_q
 (e_2^{\circ}*_q e_2^{\circ})*_q
 e_1^{\circ}\right) \\ \notag
 & = (-\lambda_{12}\lambda_{23})
w^{\circ}+\lambda_{23}\epsilon^\circ_1
+\lambda_{12}\epsilon^\circ_2 +\lambda_{2}\epsilon^\circ_3. \notag 
\end{split}\end{equation}
In a similar fashion, the other products may be computed; this amounts to
 computing $\theta$ on $T$-valued points.  
The following result is needed to compute $\Theta$ from $\theta.$ 
\begin{lemma}\label{compute-Theta-from-theta}
Let  $*_{(b(q),\underline{t})}$ denote the multiplication in the algebra 
 $\Theta(q,\underline{t})=\underline{t}.\theta(q)$
 and as before, $*_{b(q)}$ denote 
 the multiplication in $\theta(q).$  
 Then we have 
\begin{description}
\item[1.]
 $\underline{t}(\epsilon^\circ_i)=t_{i}w^{\circ}+\epsilon^\circ_i$ for
 $1\leq i \leq 3;$
\item[2.]
  $(\underline{t})^{-1}(\epsilon^\circ_i)=-t_{i}w^{\circ}+\epsilon^\circ_i$ for
 $1\leq i \leq 3;$
\item[3.]
 $\epsilon^\circ_i*_{(b(q),\underline{t})} \epsilon^\circ_j=\underline{t}(\epsilon^\circ_i*_{b(q)} \epsilon^\circ_j) -t_{j}\epsilon^\circ_i-t_{i}\epsilon^\circ_j-t_it_j w^{\circ}.$
\end{description}
\end{lemma}
\noindent While the first two of the above
formulae follow easily by direct computation, 
the third follows by using the first two  
 alongwith the following: 
$$\epsilon^\circ_i*_{(b(q),\underline{t})}\epsilon^\circ_j=\underline{t}
 \left(\ 
(\underline{t}^{-1}(\epsilon^\circ_i))*_{b(q)}
(\underline{t}^{-1}(\epsilon^\circ_j))\ 
 \right).$$
\noindent We may now 
compute 
the multiplication in
the algebra $\Theta(q,\underline{t})=\underline{t}.\theta(q)$ by making use of
the formulas listed in the above lemma and the expressions for 
the products of the form
 $\epsilon^\circ_i*_{b(q)} \epsilon^\circ_j$ whose computation had
already been illustrated before the lemma. {\bf Q. E. D., Theorem \ref{mult-table-wrt-Theta}.}
\subsection{Computation of the Isomorphism arising from a Similarity}
We continue with the notations introduced above. In the following 
 we study the lack of equivariance of the isomorphism $\Theta$ relative to 
 $\hbox{\rm GL}_{V}$ and show  
 that it satisfies a curious
 `twisted' version of equivariance. Firstly we consider the 
 morphism of $X$-groupschemes
  $$\Lambda^{even}:\hbox{\rm GL}_{V}\longrightarrow\hbox{\rm Stab}_{w}\textrm{
(given on valued points by) }g\mapsto \Lambda^{even}(g).$$
Recall that $$\Lambda^{even}(V)
=:A_0\in\hbox{\rm Id-}w\hbox{\rm -Sp-Azu}_W(X)$$
 is the even graded part of the Clifford algebra of the zero quadratic form
 on $V.$ A simple computation reveals the following result.
 \begin{lemma}\label{bgdef} For each $X$-scheme $T$, 
define the map 
$$\hbox{\rm GL}(V_T)\longrightarrow 
\hbox{\rm Stab}({(A_0)}_T): 
g\mapsto \begin{pmatrix}1 & {\bf 0}\\{\bf 0}& B(g)\end{pmatrix}\in\hbox{\rm Stab}(w_T)$$
where $B(g):=\hbox{\rm det}(g)\left(E_{12}E_{23}(g^{-1})^t E_{23}E_{12}\right)$
with 
$$E_{12}=\begin{pmatrix}0 & 1 & 0\\ 1 & 0 & 0\\ 0 & 0 & 1\end{pmatrix}
\hbox{ and }  
E_{23}=\begin{pmatrix}1 & 0 & 0\\ 0 & 0 & 1\\ 0 & 1 & 0\end{pmatrix}.$$
Then the  above maps define a morphism of $X$-groupschemes 
 which in fact is  none other than 
$$\Lambda^{even}:\hbox{\rm GL}_{V}\longrightarrow \hbox{\rm Stab}_{w};$$  
  in other 
 words: $B(g)=\Lambda^2(g).$ 
\end{lemma}
Recall from Lemma 5.1, Part A, \cite{tevb-paper1}, that 
$\hbox{\rm Stab}_{w}$
 is the semidirect product of $\hbox{\rm Stab}_{A_0}$ and $\hbox{\sf L}_w$, 
 so that $\hbox{\rm Stab}_{A_0}$ naturally acts on $\hbox{\sf L}_w$ by ``outer conjugation''. Let $\hbox{\rm GL}_{V}$ act on $\hbox{\sf L}_w$
 via the homomorphism  $\Lambda^{even}$ i.e., 
 for $g\in\hbox{\rm GL}(V_T)$ and $\underline{t}\in\hbox{\sf L}_w(T)$, 
$$g.\underline{t}:=\Lambda^{even}(g).\underline{t}:= \Lambda^{even}(g)\thinspace\underline{t}\thinspace\Lambda^{even}(g^{-1}).$$
Any element $h\in\hbox{\rm Stab}(w_T)$ \label{hshl-decomp}
 can be uniquely written as 
$$h=h_sh_l=
h_l'h_s\textrm{ where }h_s\in \hbox{\rm Stab}({(A_0)}_T)\textrm{ and }h_l,\thinspace
h_l'\in\hbox{\sf L}_w(T).$$
 Then the relation between $h_l$ and $h_l'$ can be written
 as 
$$h_l'=h_s.h_l\textrm{ or }h_l=h_s^{-1}.h_l'$$ where ``$.$''
 stands for the action of $\hbox{\rm Stab}_{A_0}$ on $\hbox{\sf L}_w.$
Thus one has a $\hbox{\rm GL}_{V}$-action
 on $\hbox{\rm Quad}_{V}\times_X \hbox{\sf L}_w$ induced by the 
diagonal embedding 
$$\hbox{\rm GL}_{V}\stackrel{\Delta}{\hookrightarrow}\hbox{\rm GL}_{V}\times_X\hbox{\rm GL}_{V}.$$
Since $\hbox{\rm Id-}w\hbox{-Sp-Azu}_W$ comes with 
 a canonical action of $\hbox{\rm Stab}_{w}$ on it, we let 
 $\hbox{\rm GL}_{V}$ act on $\hbox{\rm Id-}w\hbox{-Sp-Azu}_W$ via 
 $\Lambda^{even}.$ 
The following result describes the lack of $\hbox{\rm GL}_{V}$-equivariance 
of the isomorphism $\Theta.$

\begin{theorem}\label{tequiv} Let $T$ be an $X$-scheme. 
For $T$-valued points $g, q, \underline{t}$ respectively of $\hbox{\rm GL}_{V}$, $\hbox{\rm Quad}_{V}$, and $\hbox{\sf L}_w$, there exists
 a unique $T$-valued point of $\hbox{\sf L}_w$ given by an 
  isomorphism $h_l'(g,q)$ of ${\cal O}_T$-algebra bundles 
$$h_l'(g,q): g.\Theta(q,\underline{t}) \stackrel{\cong}{\longrightarrow}
 \Theta(g.q,g.\underline{t}).$$
 Further, $h_l'(g,q)$ satisfies the formula
$$h_l'(gg',q)=h_l'(g,g'.q)(g.h_l'(g',q)).$$ 
\end{theorem}
Therefore $\Theta$ satisfies a `twisted' version of
 $\hbox{\rm GL}_{V}$-equivariance. The next theorem,
 which was originally motivated by the proof of this `twisted equivariance',
 will be of central importance to us for the rest of
 this section.
\begin{theorem}\label{b-co-phi-lambda-indep-of-q}
Given  a 
 similarity  
$$g:(V_T, q)\thinspace\cong_l\thinspace(V_T,q')$$ with multiplier 
 $l\in\Gamma(T, {\cal O}_T^*)$,  let $h(g,l,q,q')$ be the 
 automorphism of $(W_T, w_T)$ given by the composition of the 
 following isomorphisms:
$$W_T\stackrel{{(\psi_{b(q)})}^{-1}\thinspace(\cong)}{\longrightarrow}
C_0(V_T,q)\stackrel{C_0(g,l)\thinspace{(\cong)}}{\longrightarrow}
C_0(V_T,q')\stackrel{\psi_{b(q')}\thinspace(\cong)}{\longrightarrow} W_T$$
where the algebra bundle isomorphism $C_0(g,l)$ comes from (1), 
 Prop.\ref{simil-induces-iso-of-even-cliff} and the linear isomorphisms 
 $\psi_{b(q)}$ and  $\psi_{b(q')}$ come from (2d), Theorem \ref{Bourbaki}.
In terms of actions, this means that 
$$h(g,l,q,q').\theta(q)=\theta(q').$$ 
Write $h(g,l,q,q')\in\hbox{\rm Stab}(w_T)$ 
 uniquely as a product 
$$h(g,l,q,q')=h_s(g,l,q,q')h_l(g,l,q,q')$$ with 
 the first factor in $\hbox{\rm Stab}_{A_0}(T)$ and 
 the second in $\hbox{\sf L}_w(T)$ as explained earlier. 
 Then $h_s(g,l,q,q')$ depends only on $g$ and $l$ and 
 not on $q$ or $q'.$ In fact, one has 
$$h_s(g,l,q_1,q_2)=h_s(g,l):=\begin{pmatrix}1 & {\bf 0}\\{\bf 0}& l^{-1}\Lambda^2(g)\end{pmatrix} \forall q_1,q_2\in\hbox{\rm Quad}(V_T).$$
\end{theorem}
\paragraph*{\sc Proof:} We directly compute the ${\cal O}_T$-linear
 automorphism
 $h(g,l,q,q')$ of $W_T$ as follows.
 Of course, this automorphism fixes $w^\circ=w_T.$
 So we need to only compute the images of the three remaining basis elements  
 $$\epsilon^\circ_1=e^\circ_1\wedge e^\circ_2, \epsilon^\circ_2=e^\circ_2\wedge e^\circ_3\textrm{ and }\epsilon^\circ_3=e^\circ_3\wedge e^\circ_1$$ 
 in terms of the basis elements $w^\circ$ and $\epsilon^\circ_i.$
 Let $q$ and  $l(g.q)=q'$ respectively
 correspond to the 6-tuples 
$$(\mu_1,\mu_2,\mu_3,
\mu_{12},\mu_{13},\mu_{23})\textrm{  and  }  
({{\nu}}_1,{{\nu}}_2,{{\nu}}_3,
{{\nu}}_{12},{{\nu}}_{13},{{\nu}}_{23})
\in \Gamma(T, {{\cal O}_T}^{\oplus 6}).$$
(We caution the reader that 
 $l(g.q)\neq(lg).q=l^{-2}(g.q)$!)
 Let 
 $$g\in\hbox{\rm GL}(V_T)\equiv\hbox{\rm GL}_3(\Gamma(T,{\cal O}_T))$$
 be given by the matrix $(g_{ij}).$ 
 Observe that the ${{\nu}}$ are polynomials in the $\mu$ and $g_{ij}.$
In the following computation, for the sake of clarity, we denote the 
 product in $C(V_T,q)$ by $*_q.$
For example, we have
\begin{equation}
\begin{split} 
h(g,l,q,q')\epsilon_1
&= \psi_{b(q')}\circ C_0(g,l)\left(
{(\psi_{b(q)})}^{-1}(e^\circ_1 \wedge e^\circ_2)\right) \notag\\
&=\psi_{b(q')}\left( C_0(g,l)\left(e^\circ_1*_qe^\circ_2\right)\right) 
\hbox{\hspace*{5mm} (by (2f), Theorem \ref{Bourbaki})} \notag\\
&=\psi_{b(q')}\left(l^{-1}(g(e^\circ_1)*_{q'}g(e^\circ_2))\right)
 \hbox{\hspace*{5mm}(by (1), Prop.\ref{simil-induces-iso-of-even-cliff})} \notag \\
&=l^{-1}\psi_{b(q')}\left((g_{11}e^\circ_1+g_{21}e^\circ_2+g_{31}e^\circ_3)*_{q'} (g_{12}e^\circ_1+g_{22}e^\circ_2+g_{32}e^\circ_3)\right) \notag \\
&=l^{-1}\psi_{b(q')}\left(\right.(g_{11}g_{12}{{\nu}}_1+g_{21}g_{22}{{\nu}}_2+g_{31}g_{32}{{\nu}}_3+\notag
 \\
& \quad + g_{21}g_{12}{{\nu}}_{12}+g_{31}g_{22}{{\nu}}_{23}+g_{11}g_{32}{{\nu}}_{13})w^\circ +\notag \\
&\quad +(g_{11}g_{22}-g_{21}g_{12}) e^\circ_1*_{q'}e^\circ_2 +
  (g_{21}g_{32}-g_{31}g_{22}) e^\circ_2*_{q'}e^\circ_3+ \notag \\
&\quad +(g_{31}g_{12}-g_{11}g_{32}) e^\circ_3*_{q'}e^\circ_1\left.\right) \notag \\
&=l^{-1}\left(P_1(g,l,q,q')w^\circ+C_{33}(g)\epsilon^\circ_1+C_{13}(g)\epsilon^\circ_2+C_{23}(g)\epsilon^\circ_3\right)
 \notag \\
& \quad \quad 
\hbox{\hspace*{5mm}(by (2f), Theorem \ref{Bourbaki})} \notag
\end{split}\end{equation}
where $P_1(g,l,q,q')$ is the polynomial in the ${{\nu}}$ and $g_{ij}$ 
 (as computed in the previous step) and where $C_{ij}(g)$ represents 
 the cofactor determinant of the element $g_{ij}$ of the matrix $g=(g_{ij}).$
 Similarly one computes the values of $h(g,l,q,q')\epsilon_2$ and $h(g,l,q,q')\epsilon_3.$
 Then the matrix of $h(g,l,q,q')$ is given by 
$$h(g,l,q,q')=
\begin{bmatrix}
1 & l^{-1}P_1(g,l,q,q') & l^{-1}P_2(g,l,q,q') & l^{-1}P_3(g,l,q,q') \\
0 & l^{-1}C_{33}(g) &  l^{-1}C_{31}(g) &  l^{-1}C_{32}(g) \\
0 & l^{-1}C_{13}(g) &  l^{-1}C_{11}(g) &  l^{-1}C_{12}(g) \\
0 & l^{-1}C_{23}(g) &  l^{-1}C_{21}(g) &  l^{-1}C_{22}(g) 
\end{bmatrix}$$ 
which implies that 
$$h_s(g,l,q,q')=
\begin{bmatrix}
1 & 0 & 0 & 0 \\
0 & l^{-1}C_{33}(g) &  l^{-1}C_{31}(g) &  l^{-1}C_{32}(g) \\
0 & l^{-1}C_{13}(g) &  l^{-1}C_{11}(g) &  l^{-1}C_{12}(g) \\
0 & l^{-1}C_{23}(g) &  l^{-1}C_{21}(g) &  l^{-1}C_{22}(g) 
\end{bmatrix}\hbox{ depends only on }g\hbox{ and }l.$$
Next define the matrix
$$\widehat{g}=
\begin{bmatrix}
 g_{33} &  g_{13} &  g_{23} \\
 g_{31} &  g_{11} &  g_{21} \\
 g_{32} &  g_{12} &  g_{22} 
\end{bmatrix}\hbox{ so that }h_s(g,l,q,q')=
\begin{bmatrix}
1 & {\bf 0}\\
{\bf 0} & l^{-1}{C(\widehat{g})}^t
\end{bmatrix}$$
where $C(\widehat{g})$ is the cofactor matrix of $\widehat{g}.$
Now if $E_{12}$ and $E_{23}$ are the matrices defined in Lemma \ref{bgdef}
 above, premultiplying by $E_{ij}$ has the effect of 
 interchanging the $i$th and $j$th rows, while postmultiplying has a 
similar effect on the columns.
 Thus we get $$\widehat{g}=E_{12}E_{23}(g^t)E_{23}E_{12}$$
 from which it follows that 
\begin{equation}
{C(\widehat{g})}^t = \hbox{ Adjoint }(\widehat{g}) 
 = \hbox{ det }(\widehat{g}).(\widehat{g})^{-1}
 = \hbox{ det }(g).\left(E_{12}E_{23}{(g^{-1})}^t E_{23}E_{12}\right),\notag  
\end{equation} 
showing that ${C(\widehat{g})}^t=\Lambda^2(g)$ by Lemma \ref{bgdef}. {\bf Q.E.D., Theorem 
 \ref{b-co-phi-lambda-indep-of-q}.}
\paragraph*{\sc Proof of Theorem \ref{tequiv}.}
Note that $g$ is an isometry from $(V_T,q)$ to $(V_T, g.q)$ and hence 
according to (1),
 Prop.\ref{simil-induces-iso-of-even-cliff}, 
 induces the algebra isomorphism
 $$C_0(g,l=1):C_0(V_T,q)\cong C_0(V_T,g.q).$$  Let 
$h(g,q):=h(g,1,q,g.q)$ where $h(g,l,q,q')$ was 
 defined in Theorem \ref{b-co-phi-lambda-indep-of-q} above.
 As explained in page \pageref{hshl-decomp}, 
 there are two canonical decompositions of $\hbox{\rm Stab}(w_T)$,
 which 
  lead to 
 unique (ordered) decompositions of $h(g,q)$ 
 as $$h'_l(g,q)h_s(g,q)\textrm{ and }h_s(g,q)h_l(g,q).$$ 
  By Theorem
 \ref{b-co-phi-lambda-indep-of-q} above,  
 $$h_s(g,q_1)=h_s(g,q_2)=h_s(g,1)=\Lambda^{even}(g)=:h_s(g)\ \forall q_1,q_2\in\hbox{\rm Quad}(V_T)$$ and hence we get:
\begin{equation}
\begin{split}
\Theta(g.q, g.\underline{t}) &:= (g.\underline{t}).\theta(g.q) \notag \\
 &= (\Lambda^{even}(g)\underline{t}\Lambda^{even}(g^{-1})).(h(g,1,q,g.q).\theta(q)) \notag \\
 &= (h_s(g)\underline{t}h_s^{-1}(g)).(h(g,q).\theta(q)) \notag \\
 &= ((h_s(g)\underline{t}h_s^{-1}(g))(h_s(g,q)h_l(g,q))).\theta(q) \notag \\
 &= (h_s(g)\underline{t}h_l(g,q)).\theta(q) \notag \\
 &= (h_s(g)h_l(g,q)).(\underline{t}.\theta(q)) \notag \\
 &= (h_s(g,q)h_l(g,q)).(\Theta(q,\underline{t})) \notag \\
 &= (h'_l(g,q)h_s(g,q)).(\Theta(q,\underline{t})) \notag \\
 &= h'_l(g,q).(h_s(g,q).\Theta(q,\underline{t})) \notag \\
 &= h'_l(g,q).(g.\Theta(q,\underline{t})). \notag 
\end{split}
\end{equation}
Note that $h_l(g,q)$ was explicitly computed in the proof of 
 Theorem \ref{b-co-phi-lambda-indep-of-q} above to be 
$$h_l(g,q)=
\begin{bmatrix}
1 & P_1(g,1,q,g.q) & P_2(g,1,q,g.q) & P_3(g,1,q,g.q) \\
{\bf 0} &  & I_3 & 
\end{bmatrix}\in\hbox{\sf L}_w(T).$$ 
The formula for $h'_l(g_1g_2,q)$ stated in the theorem is gotten thus:
\begin{equation}
\begin{split}
h'_l(g_1g_2,q) &=  (h'_l(g_1g_2,q)h_s(g_1g_2)) h_s^{-1}(g_1g_2) \notag \\
 &= h(g_1g_2,q)h_s^{-1}(g_1g_2). \notag    
\end{split}
\end{equation}
Now by (3) of Prop.\ref{simil-induces-iso-of-even-cliff} it follows that 
\begin{equation}
\begin{split}
h'_l(g_1g_2,q) &=  (h(g_1,g_2. q)h(g_2,q))h_s^{-1}(g_1g_2) \notag \\
 &= h'_l(g_1,g_2. q)h_s(g_1,g_2. q)h_s(g_2,q)
h_l(g_2,q)h_s^{-1}(g_1g_2)\notag  \\
 &= h'_l(g_1,g_2. q)h_s(g_1)h_s(g_2)
h_l(g_2,q)h_s^{-1}(g_1g_2)\notag  \\
 &= h'_l(g_1,g_2. q)h_s(g_1g_2)
h_l(g_2,q)h_s^{-1}(g_1g_2)\notag  \\
 &= h'_l(g_1,g_2. q)((g_1g_2). h_l(g_2,q)).\notag    
\end{split}
\end{equation}
On the other hand 
$$g_2^{-1}.h'_l(g_2,q)=h^{-1}_s(g_2)(h'_l(g_2,q)h_s(g_2))=h^{-1}_s(g_2)(h_s(g_2)h_l(g_2,q))=h_l(g_2,q)$$ and  therefore  
\begin{equation}
\begin{split}
h'_l(g_1g_2,q) &= h'_l(g_1,g_2. q).((g_1g_2)
.(g_2^{-1}. h'_l(g_2,q))) \notag \\
 &= h'_l(g_1,g_2. q).(g_1
. h'_l(g_2,q)). \notag    
\end{split}
\end{equation}
Finally, one has to show the uniqueness of $h'_l(g,q)\in\hbox{\sf L}_w(T).$
 Suppose $h_l\in\hbox{\sf L}_w(T)$ is also an algebra isomorphism 
$$h_l: g.\Theta(q,\underline{t}) \stackrel{\cong}{\longrightarrow}
 \Theta(g.q, g.\underline{t}),$$ 
 which means 
 $$h_l.( g. \Theta(q,\underline{t})) = \Theta(g. q,g.\underline{t}).$$
Notice that while showing the `twisted' equivariance of $\Theta$ above, 
 we have also proved that 
$$\Theta(g.q,g.\underline{t})=h(g,q)\Theta(q,\underline{t}).$$
 Therefore we get 
\begin{equation}
\begin{split}
(h_lh_s(g)).(\underline{t}.\theta(q)) &= 
 h(g,q).\Theta(q,\underline{t}) \notag \\
\Rightarrow (h_lh_s(g)).(\underline{t}.\theta(q)) &= 
 (h_s(g)h_l(g,q)).\Theta(q, \underline{t}) \notag \\
\Rightarrow (h^{-1}_s(g)h_lh_s(g)).(\underline{t}.\theta(q)) &= 
 (h_l(g,q)\underline{t}).\theta(q) \notag \\
\Rightarrow \Theta(q, (g^{-1}. h_l)\underline{t}) &= \Theta(q, h_l(g,q)\underline{t}). 
\notag 
\end{split}
\end{equation}
But since $\Theta$ is an  isomorphism by Theorem~\ref{Thetaisom}, this 
 implies that $$(g^{-1}. h_l)\underline{t}=h_l(g,q)\underline{t}$$ which gives
$$h_l=h_s(g)h_l(g,q)h^{-1}_s(g) = h'_l(g,q).$${\bf Q.E.D., Theorem
   \ref{tequiv}.}
\subsection{Conclusion of Proof of Injectivity}
We remind the reader that towards the end of \S \ref{sec3}, we had 
 reduced the proof of the injectivity of Theorem \ref{bijectivity}
 to that of Theorem \ref{lifting-of-isomorphisms},
 and had indicated in page \pageref{reduction-to-free-case}
 that it would be enough to prove the latter in 
 the case when $V$ and $I$ are both 
 free---which has been the case in this section so far. 
 Starting with an isomorphism of algebra bundles 
$$\phi:C_0(V,q)\cong C_0(V,q')$$
 we arrive at the element $g_l^\phi\in\hbox{\rm GL}(V)$ as defined in page 
 \pageref{glphi-def}; to briefly recall this, firstly $g\in\hbox{\rm GL}(V)$
 was defined by the following commuting diagram: 
$$\begin{CD}
{(C_0(V,q))}^\vee @>{{({(\psi_{b(q)})}^\vee)}^{-1}}>{\cong}> 
{(\Lambda^{even}(V))}^\vee @>{\hbox{\rm\tiny surjection}}>> 
{(\Lambda^2(V))}^\vee @>{=}>>\\ 
@A{{\phi^\vee}}A{\cong}A @A{{{(\phi_{\Lambda^{ev}})}^\vee}}A{\cong}A
  @A{{{(\phi_{\Lambda^2})}^\vee}}A{\cong}A\\
{(C_0(V,q'))}^\vee @>{{({(\psi_{b(q')})}^\vee)}^{-1}}>{\cong}> 
{(\Lambda^{even}(V))}^\vee @<{\hbox{\rm\tiny inclusion}}<< 
{(\Lambda^2(V))}^\vee @>{=}>> 
\end{CD}$$
$$\begin{CD}
{(\Lambda^2(V))}^\vee @>{{(\eta^\vee)}^{-1}}>{\equiv}>
V\otimes{{(\hbox{\rm det}(V))}^{-1}} 
@>{\otimes{\hbox{\rm det}(V)}}>{\equiv}> V\\
@A{{{(\phi_{\Lambda^{2}})}^\vee}}A{\cong}A
  @A{{(g')}^{-1}}A{\cong}A
@A{g^{-1}}A{\cong}A \\
{(\Lambda^2(V))}^\vee @>{{(\eta^\vee)}^{-1}}>{\equiv}> 
V\otimes{{(\hbox{\rm det}(V))}^{-1}} 
@>{\otimes{\hbox{\rm det}(V)}}>{\equiv}> V
\end{CD}$$
Secondly, we had defined $$g_l^\phi:=(l^{-1}\sqrt{\gamma(l)}\thinspace)g.$$
Our current 
 special choices of bilinear forms $b(q)$ and $b(q')$ that induce $q$ and 
 $q'$ respectively do 
 not affect the generality,
 as was observed in the proof of Prop.\ref{transfer-to-lambda2}. We shall 
 now show that $g_l^\phi$ is a 
 similitude from $(V,q)$ to $(V,q')$ with multiplier $l$ and that 
 this similitude induces $\phi$ i.e., with the notations of (1), Prop. 
\ref{simil-induces-iso-of-even-cliff},  that $C_0(g_l^\phi, l)=\phi.$
We proceed with the proof which will follow from  several lemmas. 
\begin{lemma}\label{glphi-in-matrix-form}
Consider the element 
 $$h_sh_l=h'_lh_s=h:=\phi_{\Lambda^{ev}}\in(\hbox{\rm Stab}_{w})(X)$$ 
 written uniquely as an ordered product in two different ways with  
 $h_l,h'_l\in(\hbox{\sf L}_w)(X)$ and $h_s\in(\hbox{\rm Stab}_{A_0})(X)$ 
 as explained in page \pageref{hshl-decomp};
 let $B$ be the matrix corresponding to
 $\phi_{\Lambda^2}$, and let the matrices $E_{ij}$ be as defined 
 in Lemma \ref{bgdef}. Then we have matrix representations: 
$$h_s=\begin{pmatrix}1 & {\bf 0}\\{\bf 0}& B
\end{pmatrix} \hbox{ and }g_l^\phi=\left(l^{-1}\sqrt{\gamma(l)}\right)E_{23}E_{12}{((B)^t)}^{-1}E_{12}E_{23}$$
\end{lemma}
The proof of the above lemma follows from the fact that the matrix of 
 the canonical isomorphism 
$$\eta:\Lambda^2(V)\equiv V^\vee\otimes\hbox{\rm det}(V)$$ is given by 
$E_{23}E_{12},$ which can be verified by a simple computation.
\begin{lemma}\label{mult-in-covqpr-and-l-inv-covqpr}We have the 
 formulae 
$$h_s(\hbox{\rm Identity}, l^{-1}, q', l^{-1}q')
=\begin{pmatrix}1 & {\bf 0}\\{\bf 0}& l\times I_3
\end{pmatrix}$$ and 
$$h_s(\hbox{\rm Identity}, l^{-1}, q', l^{-1}q').\theta(q')=
\theta(l^{-1}q').$$
\end{lemma}
The identity map on $V$ is obviously a similarity with multiplier 
 $l^{-1}$ from $(V,q')$ to $(V,l^{-1}q').$ 
Hence the above lemma follows
 by taking $T=X$, $g=\hbox{Identity}$, and the $l^{-1}$
 and the $q'$ at hand for the $l$ and the  
$q$ of Theorem 
 \ref{b-co-phi-lambda-indep-of-q} (caution:  
 the $q'$ there would have to be replaced by $l^{-1}q'\thinspace$!).
 This can also be seen directly from 
 the multiplication tables for $$\theta(l^{-1}q')=\Theta(l^{-1}q',I_4)
\textrm{ and }\theta(q')=\Theta(q',I_4)$$
 written out in Theorem \ref{mult-table-wrt-Theta}, 
 where we must take $T=X$ and $\underline{t}=I_4$
 i.e., $t_i=0\thinspace\forall\thinspace i.$ 
 We observe from the multiplication table that each of the coefficients of 
 $\epsilon_i$ for $1\leq i\leq 3$ is a single $\lambda$, whereas each 
 coefficient of $w=1=\epsilon_0$ is a product of two $\lambda\thinspace$s, 
 and this observation implies the lemma above. 

As the reader might have noticed, there are two crucial facts
 about the identifications in this section; namely, firstly,  for 
 any $X$-scheme $T$, each of 
 the maps $\psi_{b(q)}$ (for different $q$) 
 identify $(C_0(V_T,q),1)$ with the same 
 $(W_T,w_T)$ and secondly, relative to the bases chosen, all these identifying 
 maps have trivial determinant. The latter is also true of the identification 
 $\eta$, since it is given by the matrix $E_{23}E_{12}$ 
 (as was noted after Lemma \ref{glphi-in-matrix-form}).
  It therefore follows that 
 $$\hbox{\rm det}(\phi)=\hbox{\rm det}(\phi_{\Lambda^{ev}})=\hbox{\rm det}(\phi_{\Lambda^2})=\hbox{\rm det}(g')=\hbox{\rm det}(g)=\hbox{\rm det}(B^{-1}).$$
 But we had chosen $l\in\Gamma(X, {\cal O}_X^*)$ such that 
$$\gamma(l):=(l^3).\hbox{\rm det}(\phi_{\Lambda^2})=l^3\hbox{\rm det}(B).$$ 
 Using these facts alongwith Lemma \ref{glphi-in-matrix-form} above, 
 a straightforward computation gives the following.
\begin{lemma}\label{B-versus-Bglphi} We have the equality 
$$\hbox{\rm det}(g_l^\phi)=\sqrt{\gamma(l)}$$ from which
 it follows that $$B(g_l^\phi)=l\times B$$ where $B(g_l^\phi)$ and $B$ are
 as defined in Lemmas \ref{bgdef} and \ref{glphi-in-matrix-form} respectively.
In particular, 
$$\hbox{\rm det}^2(g_l^\phi)=\hbox{\rm det}
(\phi_{\Lambda^2})$$
 when $\hbox{\rm det}(\phi_{\Lambda^2})$ is itself a square and 
 for the cases 1 and 2 of page \pageref{cases} where 
we had chosen $l:=1.$ 
\end{lemma}
\begin{lemma}\label{glphi-simil-q-qpr-with-mult-l}
 $g_l^\phi$ is a 
 similitude from $(V,q)$ to $(V,q')$ with multiplier $l.$
\end{lemma}
 The hypothesis $\phi:C_0(V,q)\cong C_0(V,q')$ is an algebra isomorphism 
 translates in terms of actions into $h.\theta(q)=\theta(q')$ where 
$h=\phi_{\Lambda^{ev}}\in(\hbox{\rm Stab}_w)(X).$   
Let 
$$h(g_l^\phi,q):=h(g_l^\phi,1,q,g_l^\phi.q)$$ where $h(g,l,q,q')$ was 
 defined in Theorem \ref{b-co-phi-lambda-indep-of-q} above. Then we have
\begin{equation}
\begin{split}
\Theta(g_l^\phi.q,I_4)&:=\theta(g_l^\phi.q)
=h(g_l^\phi,q).\theta(q)=h(g_l^\phi,q).(h^{-1}.\theta(q')) \notag \\
&=\left(h'_l(g_l^\phi,q)h_s(g_l^\phi,q){h_s}^{-1}{(h'_l)}^{-1}\right).\theta(q') \notag \\
&=\left(h'_l(g_l^\phi,q)\begin{pmatrix} 1 & {\bf 0}\\{\bf 0} & B(g_l^\phi)\end{pmatrix}
\begin{pmatrix} 1 &{\bf 0}\\ {\bf 0} & B^{-1}\end{pmatrix}
{(h'_l)}^{-1}\right).\theta(q') \notag \\
& \notag \\
&\quad  \hspace*{4mm}\hbox{\rm (by Theorem
 \ref{b-co-phi-lambda-indep-of-q}; 
 Lemmas \ref{bgdef} \& \ref{glphi-in-matrix-form})} \notag \\
& \notag \\
&=\left(h'_l(g_l^\phi,q)\begin{pmatrix} 1 &{\bf 0}\\{\bf 0} & l\times I_3\end{pmatrix}
{(h'_l)}^{-1}\right).\theta(q')
\hspace*{6mm}\hbox{\rm (by Lemma \ref{B-versus-Bglphi})} \notag \\
&=\left(h'_l(g_l^\phi,q)\thinspace h_s(\hbox{\rm Identity}, l^{-1}, q', l^{-1}q')
\thinspace {(h'_l)}^{-1}\right).\theta(q') 
\hspace*{2mm}\hbox{\rm (by Lemma \ref{mult-in-covqpr-and-l-inv-covqpr})} \notag \\
&=(h'_l(g_l^\phi,q)h''_l).\left(h_s(\hbox{\rm Identity}, l^{-1}, q',
l^{-1}q').\theta(q')\right) \notag \\
& \notag \\ 
& \quad \hspace*{6mm}\hbox{\rm (since $\hbox{\rm Stab}_{w}$ is a semidirect product)} 
 \notag \\
& \notag \\ 
&=(h'_l(g_l^\phi,q)h''_l).\theta(l^{-1}q')) 
\hspace*{6mm}\hbox{\rm (by Lemma \ref{mult-in-covqpr-and-l-inv-covqpr})}
  \notag \\
&=:\Theta(l^{-1}q',(h'_l(g_l^\phi,q)h''_l)). \notag
\end{split}
\end{equation}
But since $\Theta$ is an isomorphism (Theorem \ref{Thetaisom}), this implies the claim of the above 
 lemma namely, that 
$$g_l^\phi.q=l^{-1}q'\textrm{ and further that
 }h'_l(g_l^\phi,q)={(h''_l)}^{-1}.$$
\begin{lemma}\label{glphi-induces-phi} The similarity 
 $$g_l^\phi:(V,q){\cong}_l(V,q')$$ induces $\phi$ i.e.,
 with the notations of (1), Prop.\ref{simil-induces-iso-of-even-cliff},
 $C_0(g_l^\phi, l)=\phi.$
\end{lemma}
We have $C_0(g_l^\phi, l)=\phi$ iff 
\begin{equation}\begin{split}
h(g_l^\phi, l,q,q')&:=
\psi_{b(q')}\thinspace\circ\thinspace C_0(g_l^\phi, l)
 \thinspace\circ\thinspace \psi_{b(q)}^{-1} \notag \\
&=\psi_{b(q')}\thinspace\circ\thinspace\phi \thinspace\circ\thinspace
 \psi_{b(q)}^{-1}=: \phi_{\Lambda^{ev}}=:h.\notag \end{split}\end{equation}
Now using successively Theorem \ref{b-co-phi-lambda-indep-of-q},
Lemma \ref{bgdef}, Lemma \ref{B-versus-Bglphi} and Lemma \ref{glphi-in-matrix-form}, 
we get the following sequence of equalities: 
$$h_s(g_l^\phi,l,q,q')=\begin{pmatrix}1 & {\bf 0}\\{\bf 0}& l^{-1}\times\Lambda^2(g_l^\phi)\end{pmatrix}=\begin{pmatrix}1 & {\bf 0}\\{\bf 0}& l^{-1}\times B(g_l^\phi)\end{pmatrix}=\begin{pmatrix}1 & {\bf 0}\\{\bf 0}& B\end{pmatrix}=h_s.$$
Therefore the present hypotheses translated in terms of actions give
\begin{equation}\begin{split}
h(g_l^\phi,l,q,q').\theta(q)=\theta(q')&=h.\theta(q)\notag \\
\Rightarrow h_s(g_l^\phi,l,q,q').(h_l(g_l^\phi,l,q,q').\theta(q))&=
h_s.(h_l.\theta(q))\notag \\
\Rightarrow\Theta(q,h_l(g_l^\phi,l,q,q'))&=\Theta(q,h_l).\notag
\end{split}\end{equation}
But $\Theta$ being an isomorphism (Theorem \ref{Thetaisom}),
 the last equality implies that 
 $$h_l(g_l^\phi,l,q,q')=h_l\textrm{ which gives
 }h(g_l^\phi,l,q,q')=h.$$
\begin{lemma}\label{locally-det-coglambda-equals-detgbylcube}\
\begin{description}
\item[(1)]
For a similarity $g\in\simvqvqpr$
 with multiplier $l$ and the
 induced isomorphism 
$$C_0(g,l)\in\isomcovqcovqpr$$ 
 given by (1), Prop.\ref{simil-induces-iso-of-even-cliff},
 we have
 the equality 
$$\hbox{\rm det}({(C_0(g,l))}_{\Lambda^2})
=l^{-3}\hbox{\rm det}^2(g).$$
 Therefore the map $$\simvqvqpr \longrightarrow
 \isomcovqcovqpr:g\mapsto C_0(g,l)$$ 
maps the subsets 
$$\isomvqvqpr\textrm{ and }\sisomvqvqpr$$ respectively   
 into the subsets 
$$\isomcovqcovqprprime\textrm{ and }\sisomcovqcovqpr.$$ 
\item[(2)] In the case $q'=q$, if $C_0(g,l)$ is the identity on 
$\covq$, then 
$$g=l^{-1}\hbox{\rm det}(g)\times\hbox{\rm Id}_V,$$
 and 
 further if $g\in\ovq$ then 
$$g=\hbox{\rm det}(g)\times\hbox{\rm Id}_V\textrm{ with }
\hbox{\rm det}^2(g)=1.$$
\end{description}
\end{lemma} 
By definition, 
$${(C_0(g,l))}_{\Lambda^{ev}}=
\psi_{b(q')}\thinspace\circ\thinspace C_0(g,l)
\thinspace\circ\thinspace\psi_{b(q)}^{-1},$$ and the latter isomorphism is 
$h(g,l,q,q')$ from Theorem \ref{b-co-phi-lambda-indep-of-q} which 
further gives a formula for    
$h_s(g,l,q,q').$ 
Now using the facts that 
 $\psi_{b(q)}$ and $\psi_{b(q')}$ have trivial determinant (as noted 
 before Lemma \ref{B-versus-Bglphi}) 
we get assertion (1): 
\begin{equation}
\begin{split}
\hbox{\rm det}\left({(C_0(g,l))}_{\Lambda^2}\right)
&=\hbox{\rm det}\left({(C_0(g,l))}_{\Lambda^{ev}}\right)
=\hbox{\rm det}(h(g,l,q,q'))\notag \\
&=\hbox{\rm det}(h_s(g,l,q,q')) 
=l^{-3}\hbox{\rm det}^2(g).\end{split}\end{equation}
 If $q=q'$ and $C_0(g,l)$ is the identity, then
  the same argument in fact shows that 
$$l^{-1}\Lambda^2(g)=I_3$$
 and by using the formula in Lemma \ref{bgdef} for 
 $B(g)=\Lambda^{2}(g)$, we get 
$$g=l^{-1}\hbox{\rm det}(g)I_3;$$  
 taking determinants in the last equality 
gives $\hbox{\rm det}^2(g)=l^3$, so that when $g\in\ovq$ i.e., $l=1$, 
 $$\hbox{\rm det}^2(g)\in\mu_2(\Gamma(X, {\cal O}_X))$$ and 
assertion (2) follows. 
\begin{lemma}\label{bijection-of-special-isomorphisms}
 The following  map  is a bijection: 
$$\sisomvqvqpr\longrightarrow\sisomcovqcovqpr:g\mapsto C_0(g,1,q,q').$$ 
\end{lemma}
Given $$\phi\in\sisomcovqcovqpr,$$ by definition  
 \ref{defn-of-iso-sets} we have 
$\hbox{\rm det}(\phi_{\Lambda^2})=1$, so by Lemma \ref{B-versus-Bglphi}
$$\hbox{\rm det}(g_l^\phi)=\sqrt{\gamma(l)}:=1$$ for our 
 choice under Case 1 on page \pageref{cases}. Therefore the
 corresponding element 
 $$g_l^\phi\in\sisomvqvqpr$$ and is, according to Lemma 
\ref{glphi-induces-phi}, such that    
 $$C_0(g_l^\phi,l=1,q,q')=\phi$$ which gives the surjectivity.
As for the injectivity, if
 $$g_1,g_2:(V,q){\cong}_1(V,q')$$ are isometries with determinant 1 such 
 that $$C_0(g_1,1,q,q')=C_0(g_2,1,q,q'),$$ then we have 
 $$h(g_1,1,q,q')=h(g_2,1,q,q')\textrm{ so that }
h_s(g_1,1,q,q')=h_s(g_2,1,q,q')$$ 
 whence by Theorem
 \ref{b-co-phi-lambda-indep-of-q} and 
Lemma \ref{bgdef} $$B(g_1)=B(g_2)\Rightarrow g_1=g_2.$$ 
{\bf Q.E.D., Theorem \ref{lifting-of-isomorphisms} and injectivity of 
 Theorem \ref{bijectivity}.}
 \paragraph*{\sc Proof of Theorem \ref{lifting-of-automorphisms}.}
 Taking 
 $q'=q$ in Theorem \ref{lifting-of-isomorphisms} gives the commutative 
 diagram of groups and homomorphisms as asserted in the statement of 
 the theorem. We continue with the notations above. 
 For $g\in\govqi$ with multiplier $l$, that 
 the equality $$\hbox{\rm det}(C_0(g,l, I))
=\hbox{\rm det}\left((C_0(g,l, I))_{\Lambda^2}\right)
=l^{-3}\hbox{\rm det}^2(g)$$ holds (locally, hence globally) was 
 shown in (1),
 Lemma \ref{locally-det-coglambda-equals-detgbylcube}.
 Assertion (2) of the same Lemma shows the following (locally, and 
 hence globally):  
 if $C_0(g,l, I)$ is the identity on 
$\covqi$, then 
$$g=l^{-1}\hbox{\rm det}(g).\hbox{\rm Id}_V,$$
 and 
 further if $g\in\ovqi$ then 
$$g=\hbox{\rm det}(g).\hbox{\rm Id}_V\textrm{ with }\hbox{\rm det}^2(g)=1.$$ 
The map 
$$\Gamma(X, {\cal O}_X^*)\longrightarrow \govqi$$
is the natural one given by sending $\lambda$ to the similarity 
 $\lambda.\textrm{Id}_V$ with multiplier $\lambda^2.$ It follows
from  the formula in (1), Prop.\ref{simil-induces-iso-of-even-cliff} 
 that 
$$C_0(\lambda.\textrm{Id}_V, \lambda^2, I)=\textrm{Identity}.$$
 This gives 
 exactness at $\govqi$ and at $\ovqi.$ 
We proceed to prove assertion (b). 
 Let $$\phi\in\autcovqi,$$ and consider the self-similarity
 $$g_l^\phi=s_{2k+1}^{+}(\phi)\textrm{ with multiplier }
 l=\hbox{\rm det}(\phi)^{2k+1}.$$ 
 For the moment, assume that $V$ and $I$ are trivial over 
 $X.$ Fix a
  global basis $\{e_1,e_2,e_3\}$ for $V$ 
  and set $e'_i=g_l^\phi(e_i).$ 
  It follows from Kneser's  definition 
 of the half-discriminant $d_0$---see formula (3.1.4), Chap.IV, 
\cite{Knus}---that 
$$d_0(q, \{e_i\})
=d_0(q, \{e'_i\})\thinspace\hbox{\rm det}^2(g_l^\phi).$$ 
 Since we have 
 $$g_l^\phi.q=l^{-1}q,$$ a simple computation shows that 
$$d_0(q, \{e'_i\})=
l^{3}d_0(q,\{e_i\}).$$
 The hypothesis that $q\otimes\kappa(x)$ is 
 semiregular means that the image of the 
 element 
$$d_0(q,\{e_i\})\in\Gamma(X,{\cal O}_X)$$ in $\kappa(x)$ is nonzero. 
 Since $X$ is integral, we therefore deduce that 
$$\hbox{\rm det}^2(g_l^\phi)=l^{-3}.$$ On the other hand, we know that 
 $$\hbox{\rm det}^2(g_l^\phi)l^{-3}=\hbox{\rm det}(\phi).$$ It follows 
 that $$\hbox{\rm det}^{12k+7}(\phi)
=1\thinspace\forall\thinspace k\in\mathbb Z,$$
 which forces $\hbox{\rm det}(\phi)=1.$ In general, even if $V$ and 
 $I$ are not necessarily trivial, since this equality holds 
 over a covering of $X$ which trivialises both $V$ and $I$, it 
 also holds over all of $X.$
{\bf Q.E.D., Theorem \ref{lifting-of-automorphisms}.}

\section{Surjectivity of Theorem~\ref{bijectivity}: Bilinear Forms as Specialisations}
\label{sec5}
In this section we reduce the proof of Theorem \ref{bilinear-forms-as-specialisations} to Theorem \ref{multiplication-table-wrt-bilinear}. We 
 prove the latter and using it  
  alongwith Theorem \ref{mult-table-wrt-Theta},  
 deduce Theorem \ref{Thetaisom-Upsilonisom}. The surjectivity of 
 Theorem \ref{bijectivity} is also established. 
\paragraph*{\sc Reduction of Proof of
 Theorem \ref{bilinear-forms-as-specialisations} to the case $I={\cal O}_X.$}
We adopt the notations introduced just before 
  Theorem \ref{bilinear-forms-as-specialisations}.
 Let $T$ be an $X$-scheme. Given a bilinear form 
$$b\in\hbox{\rm Bil}_{(V,I)}(T)
=\Gamma(T, \hbox{\rm Bil}_{(V_T, I_T)}),$$
 consider the linear isomorphism 
$$\psi_{b}:C_0(V_T, q_b, I_T)
\cong {\cal O}_T\oplus\Lambda^2(V_T)\otimes {(I_T)}^{-1}=W_T$$ of 
 (2d), Theorem \ref{Bourbaki}. Let $A_b$ denote the algebra bundle structure
 on $W_T$ with unit $w_T=1$ induced via $\psi_b$ 
 from the even Clifford algebra $C_0(V_T, q_b, I_T)$. By definition, 
 $A_b\in\aaWw(T)$ and we get a map of $T$-valued points 
 $$\Upsilon'(T):\hbox{\rm Bil}_{(V,I)}(T)
\longrightarrow \aaWw(T): b\mapsto A_b.$$
 This is functorial in $T$ because of (3), Theorem \ref{Bourbaki}, and 
 hence defines an $X$-morphism 
$$\Upsilon':\hbox{\rm Bil}_{(V,I)}
\longrightarrow 
 \aaWw.$$ The morphism $\Upsilon'$ is $\hbox{\rm GL}_V$-equivariant 
 due to (2j), Theorem \ref{Bourbaki}. Notice that the schemes 
 $\hbox{\rm Bil}_{(V,I)}, \hbox{\rm Bil}^{sr}_{(V,I)}$ and $\aaWw$ are 
 well-behaved relative to $X$ with respect to base-change. In fact, 
 so are $\azuWw$ and $\spazuWw$, as may be recalled from Theorems 
 \ref{repbility-smoothness-azu} and \ref{smoothness-spazu} of page 
 \pageref{repbility-smoothness-azu}.
In view of these observations, by taking a trivialisation for $I$ over 
 $X$, we may reduce to the case when $I$ is trivial. We treat this
case next.
\paragraph*{\sc Reduction of Proof of
 Theorem \ref{bilinear-forms-as-specialisations} for $I={\cal O}_X$  
to Theorem \ref{multiplication-table-wrt-bilinear}.}
 We first recall the following crucial fact 
(see (1), Prop.3.2.4, Chap.IV \cite{Knus}): 
{\em The even Clifford algebra of a semiregular quadratic form is an 
 Azumaya algebra.} 
 Using this fact and the definition of $\Upsilon'$, we see that
 the morphism $\Upsilon'$ restricted to $\hbox{\rm Bil}_V^{sr}$  
  factors through $\azuWw$ by a morphism $\Upsilon^{sr}$ such 
 that   the following  diagram is commutative  
$$\begin{CD}
\hbox{\rm Bil}_V @>{\Upsilon'}>> \aaWw \\
@AAA @AAA \\
\hbox{\rm Bil}_V^{sr} @>>{\Upsilon^{sr}}> \azuWw 
\end{CD}$$
 where the vertical arrows are the canonical open immersions. 
The above diagram base changes well in view of 
(2), Theorem \ref{repbility-smoothness-azu}, Prop.\ref{bilvsr-def} and 
 (3), Theorem \ref{Bourbaki}.  
Notice that since the structure morphism $\hbox{\rm Bil}_V\longrightarrow X$
 is an affine morphism, and since the same is true of 
$$\aaWw\longrightarrow X,$$
 it is also true of $\Upsilon'.$ 
 In particular, $\Upsilon'$ is quasi-compact and separated, and 
 therefore has a schematic image by case 
 (1) of Prop.\ref{sch-img-existence}. The same is true of
 each of the two vertical arrows and of 
 $\Upsilon^{sr}$ 
 in view of  
 Prop.\ref{bilvsr-def} and  (1) of Theorem 
\ref{repbility-smoothness-azu}.
 Further,  as noted in Prop.\ref{bilvsr-def},  
 $$\hbox{\rm Bil}_V^{sr}\hookrightarrow\hbox{\rm Bil}_V$$ is schematically 
 dominant and therefore by (5), Prop.\ref{sch-img-props},
 the limiting scheme of the former in the latter 
 is the latter itself. 
  So using the commutative diagram above, the transitivity 
 of the schematic image (assertion (3), Prop.\ref{sch-img-props}),
 and the definition of $\spazuWw$ (assertion (1), 
Theorem \ref{smoothness-spazu}), we see 
 that in order to prove (1), Theorem \ref{bilinear-forms-as-specialisations},
 it is enough to show that
\begin{description}
\item[(*)] $\Upsilon^{sr}$ is schematically dominant and surjective, and 
 $\Upsilon'$ is a closed immersion. 
\end{description}
We now claim that the above properties are equivalent to 
\begin{description}
\item[(**)] $\Upsilon^{sr}$ is proper and  $\Upsilon'$ is a closed immersion. 
\end{description}
 Suppose {\bf (**)} holds. To show {\bf (*)}, we only need 
 to show that $\Upsilon^{sr}$ is surjective and schematically dominant.
 From {\bf (**)} it follows 
 that  
$$\Upsilon_K^{sr}:=\Upsilon^{sr}\otimes_X K$$ is functorially injective 
 and proper for each algebraically closed field $K$ with 
 an $X$-morphism $$\hbox{\rm Spec}(K)\longrightarrow X.$$
  That both the $K$-schemes 
$$\hbox{\rm Bil}_V^{sr}\otimes_X K\textrm{ and }
\azuWw\otimes_X K$$ are integral and 
 smooth of the same dimension follows from 
 the smoothness of relative dimension 9
  and geometric irreducibility $/X$ of 
 $\hbox{\rm Bil}_V^{sr}$ (which is obvious), and of 
  $\azuWw$  from (3), Theorem \ref{repbility-smoothness-azu}. 
 Since $\Upsilon_K^{sr}$ is differentially injective at each closed 
 point, it has 
 to be a smooth morphism by Theorem 17.11.1 of EGA IV \cite{ega4}
 and 
 thus has to be an open map. 

But by {\bf (**)} 
 it is also proper and hence a closed map. 
 Thus $\Upsilon_K^{sr}$ 
 is bijective etale, and hence an isomorphism. This also gives 
 that $\Upsilon^{sr}$ is surjective.  
 Now from Cor.11.3.11 of EGA IV and from the flatness of
 $\hbox{\rm Bil}_V^{sr}$ 
 over $X$,  
 it follows that $\Upsilon^{sr}$ is itself flat, and hence schematically 
 dominant since it is faithfully flat (being already surjective). 
 Therefore {\bf (**)}$\implies${\bf(*)}. 
\paragraph*{}
The property of a morphism being proper is
 local on the target (see for e.g., (f), Cor.4.8, Chap.IV, \cite{Hartshorne}) 
  and the same is true of the property of being a closed immersion. 
 Therefore, in verifying (**), we may assume that $V$ is free over $X$ 
 so that $$W=\Lambda^{even}(V)$$ is also free over $X$ and 
  and $w$ is part of a global basis. We are now in the situation of 
 Theorem \ref{multiplication-table-wrt-bilinear}. 
Granting it, we see immediately 
 from the multiplication table that {\bf (**)} holds.
 For the table shows that the composition of the following $X$-morphisms
 $$\hbox{\rm Bil}_V\stackrel{\Upsilon'}{\longrightarrow}\aaWw\stackrel{\hbox{\rm\tiny CLOSED}}{\hookrightarrow}\algW$$ is a closed immersion, which 
 implies that $\Upsilon'$ is also a closed immersion. Further, 
  the multiplication 
 table also shows that both $\Upsilon'$ and $\Upsilon^{sr}$ satisfy the 
 valuative criterion for properness, and are therefore
 proper. Thus 
the conditions {\bf (**)} are verified.  So we have reduced  
 the proof of (1), Theorem \ref{bilinear-forms-as-specialisations}
 to Theorem \ref{multiplication-table-wrt-bilinear}.
\paragraph*{} As for statement (2) 
 of Theorem \ref{bilinear-forms-as-specialisations}, firstly, 
  the  involution 
 $\Sigma$ of  $\aaWw$ 
 defines a unique involution (also denoted by $\Sigma$) on 
 the scheme of specialisations $\spazuWw$
 (leaving the open subscheme $\azuWw$
 invariant) because of   
  the defining property of the 
 schematic image involved; for we may verify  that an automorphism of 
 a scheme $T$ which leaves an open subscheme $U$ stable will also leave 
 stable the limiting scheme of $U$ in $T$ (of course we assume here that 
 the canonical open immersion $U\hookrightarrow T$ is a quasi-compact 
 open immersion, which ensures the existence of the limiting scheme). Secondly,
 a glance at the multiplication table of Theorem \ref{multiplication-table-wrt-bilinear}
 keeping in view the definition of opposite algebra shows that indeed 
 the induced $\Sigma\in\hbox{\rm Aut}_X\left(\hbox{\rm Bil}_V\right)$ takes the $T$-valued point 
$B=(b_{ij})$ to $\hbox{\rm transpose}(-B)=(-b_{ji}).$ Finally, assertion (3)
 of Theorem \ref{bilinear-forms-as-specialisations} is a consequence of 
 (1) taking into account (3), Theorem \ref{smoothness-spazu}.
\paragraph*{\sc Proof of Theorem \ref{multiplication-table-wrt-bilinear}.}
Given $B=(b_{ij})\in\hbox{\rm Bil}_V(T)$, by our definition above, 
 $(\Upsilon'(T))(B)=A_B$ is the algebra structure induced from 
 the linear isomorphism 
$$\psi_B:C_0(V_T, q_B)\cong \Lambda^{even}(V_T)$$ 
 of (2d), Theorem \ref{Bourbaki}. 
The stated multiplication table for $A=A_B$
 is a consequence of straightforward calculation, keeping in 
 mind (2f), Theorem \ref{Bourbaki} and the standard properties 
 of the multiplication in the even Clifford algebra $C_0(V_T,q_B).$
{\bf Q.E.D., Theorems \ref{multiplication-table-wrt-bilinear} \& \ref{bilinear-forms-as-specialisations}.}
\paragraph*{\sc Proof of Theorem \ref{Thetaisom-Upsilonisom}.}
The proof follows by comparing the multiplication table relative to 
 $\Theta$ as computed in Theorem \ref{mult-table-wrt-Theta} with 
 the multiplication table relative to $\Upsilon$ of Theorem  
\ref{multiplication-table-wrt-bilinear} computed above. 
{\bf Q.E.D., Theorem \ref{Thetaisom-Upsilonisom}.}
\paragraph*{\sc Proofs of assertions in
  ({\rm a}), Theorem \ref{structure-of-specialisation}
 and Surjectivity part of Theorem \ref{bijectivity}.} 
Let $W$ be the rank 4 vector bundle underlying 
 the specialised algebra $A$ and $w\in\Gamma(X, W)$
 be the global section corresponding to 
 $1_A.$ We choose an affine open covering 
${\{U_i\}}_{i\in \mathcal{I}}$ of $X$ such that $W|U_i$ is trivial and 
 $w|U_i$ is part of a global basis $\forall\thinspace i.$  
 Therefore on the one hand, 
for each $i\in \mathcal{I}$, we can find a linear isomorphism
 $$\zeta_i:\Lambda^{even}\left({\cal O}_X^{\oplus 3}|U_i\right)\cong 
 W|U_i$$ taking $1_{\Lambda^{even}}$ onto $w|U_i.$ The $(w|U_i)$-unital
 algebra structure  $A|U_i$ induces via $\zeta_i$ an algebra structure $A_i$ on
 $\Lambda^{even}\left({\cal O}_X^{\oplus 3}|U_i\right)$ (so that $\zeta_i$
 becomes an algebra isomorphism). Recall that $A_i$ is also 
 a specialised algebra structure by (3),
 Theorem~\ref{smoothness-spazu}.
 Hence by  Theorem \ref{bilinear-forms-as-specialisations} 
 applied to $X=U_i$, $V={\cal O}_{U_i}^{\oplus 3}$ and $I={\cal O}_{U_i}$, 
  we can also
 find an ${\cal O}_{U_i}$-valued quadratic form 
 $q_i$ on ${\cal O}_X^{\oplus 3}|U_i$ induced from a bilinear form 
 $b_i$ so that the algebra structure $A_i$ is precisely the one 
 induced by the linear isomorphism 
 $$\psi_{b_i}:C_0\left({\cal O}_X^{\oplus 3}|U_i, q_i\right)
\cong\Lambda^{even}\left({\cal O}_X^{\oplus 3}|U_i\right)$$ 
given by (2d) of Theorem \ref{Bourbaki}.  
 For each pair of indices $(i,j)\in \mathcal{I}\times \mathcal{I}$,
 let $\zeta_{ij}$ and 
 $\phi_{ij}$ be defined so that the following diagram commutes:
$$\begin{CD}
C_0({\cal O}_X^{\oplus 3}|U_{ij}, q_i|U_{ij}) 
@>{\psi_{b_i}|U_{ij}}>{\cong}> \Lambda^{ev}({\cal O}_X^{\oplus 3}|U_{ij})
@>{\zeta_i|U_{ij}}>{\cong}> A|U_{ij}\\
@V{\phi_{ij}}V{\cong}V @V{\zeta_{ij}}V{\cong}V @V{=}VV\\
C_0({\cal O}_X^{\oplus 3}|U_{ij}, q_j|U_{ij}) 
@>{\cong}>{\psi_{b_j}|U_{ij}}> \Lambda^{ev}({\cal O}_X^{\oplus 3}|U_{ij})
@>{\cong}>{\zeta_j|U_{ij}}> A|U_{ij}
\end{CD}$$
The above diagram means that the algebras $A_i$ glue along 
 $U_{ij}:=U_i\cap U_j$ via $\zeta_{ij}$ to give (an algebra bundle 
 isomorphic to) $A$, and in the same vein, the even Clifford algebras 
 $C_0({\cal O}_X^{\oplus 3}|U_i, q_i)$ glue along the $U_{ij}$ via 
 $\phi_{ij}$ to give $A$ as well.  
Now consider the similarity 
$$g_{l_{ij}}^{\phi_{ij}}=s_{-1}^{+}(\phi_{ij}):({\cal O}_X^{\oplus
 3}|U_{ij}, q_i|U_{ij})\thinspace\cong_{l_{ij}}\thinspace({\cal
 O}_X^{\oplus 3}|U_{ij}, q_j|U_{ij})$$ with multiplier
 $$l_{ij}:=\hbox{\rm det}(\phi_{ij})^{-1}$$
 given by (c), 
 Theorem \ref{lifting-of-isomorphisms}. Since $s_{-1}^{+}$
 is  multiplicative, and since $\phi_{ij}$ satisfy 
 the cocycle condition, it follows that $s_{-1}^{+}(\phi_{ij})$ also 
 satisfy the cocycle condition and therefore glue the
 ${\cal O}_X^{\oplus 3}|U_i$
 along the $U_{ij}$ to give a rank 3 vector bundle $V$ on $X.$ While 
 the $q_i$ do not glue to give an ${\cal O}_X$-valued quadratic form on 
 $V$, the facts that the multipliers $\{l_{ij}\}$ form a cocycle for 
 $$I:=\hbox{\rm det}^{-1}(A)$$ and that 
 $s^+_{(-1)}$ is a section together imply, taking into account the 
 uniqueness in (1), Prop.\ref{simil-induces-iso-of-even-cliff},
  that actually the $q_i$ glue to give an 
 $I$-valued quadratic form $q$ on $V$ and 
 that $\covqi\cong A.$ 
We shall now 
 revert to the notations of Section \ref{sec4}.
 By Theorem \ref{b-co-phi-lambda-indep-of-q}, we have 
$$h_s(g_{l_{ij}}^{\phi_{ij}},l_{ij},q_i|U_{ij},q_j|U_{ij})=\begin{pmatrix}1 & {\bf 0}\\{\bf 0}& l_{ij}^{-1}\Lambda^2(g_{l_{ij}}^{\phi_{ij}})\end{pmatrix}$$
 which means that 
$${(\phi_{ij})}_{\Lambda^2}=\hbox{\rm det}(\phi_{ij})\times
 \Lambda^2(g_{l_{ij}}^{\phi_{ij}}).$$ 
 This immediately implies part (1) of
 assertion (a) of Theorem \ref{structure-of-specialisation}, from which 
 parts (2)---(4) can be deduced using the standard properties of the
  determinant  
 and the perfect 
 pairings between suitable exterior powers of a bundle.
\paragraph*{\sc Proofs of assertions in ({\rm b}),
 Theorem \ref{structure-of-specialisation}.}  We first prove (b1). 
Let $A$ be a given specialisation, and let $A\cong\covqi$ as in part 
 (a) of Theorem \ref{structure-of-specialisation} with
 $I=\hbox{\rm det}^{-1}(A).$ By the injectivity part of 
 Theorem \ref{bijectivity}, we have 
$$\covqi\cong A\cong C_0(V',q',{\cal O}_X)$$ iff 
there exists a twisted discriminant bundle $(L,h,J)$ and an isomorphism 
 $$\vqi\cong (V',q',{\cal O}_X)\otimes (L,h,J).$$ The latter implies that 
 $I\cong J\cong L^2$ and hence $\hbox{\rm det}(A)\in 2.\hbox{\rm Pic}(X).$ 
 On the other hand, if this last condition holds, we could take 
 for $L$ a square root of $J:=I^{-1}$, alongwith an isomorphism
 $h:L^2\cong J$ and we would have by Prop.\ref{isom-covq-covlqh}
 an algebra isomorphism 
$$\gamma_{(L,h,J)}:C_0(V\otimes L, q\otimes h, {\cal O}_X)\equiv
C_0\left(\vqi\otimes(L,h,J)\right)\cong C_0(V,q,I)\cong A.$$
For the proof of (b2), 
suppose that the line subbundle ${\cal O}_X.1_A\hookrightarrow A$ is a
 direct summand of $A.$  We
 may choose a splitting 
$$A\cong {\cal O}_X.1_A\oplus 
 (A/{\cal O}_X.1_A).$$
 Using assertion (1) of (a),
 Theorem \ref{structure-of-specialisation}, we see that there exists a
 rank 3 vector bundle $V$ on $X$ such that
 \begin{equation}\begin{split}
A&\cong {\cal O}_X.1_A\oplus (A/{\cal O}_X.1_A)\notag\\
&\cong {\cal O}_X.1_A\oplus 
 (\Lambda^2(V)\otimes I^{-1})\cong {\cal O}_X.1\oplus \Lambda^2(V)\otimes 
 I^{-1}=:W\notag\end{split}\end{equation}
 where $I:=\hbox{\rm det}^{-1}(A)$ and   
 the last isomorphism is chosen so as 
 to map ${\cal O}_X.1_A$ isomorphically onto ${\cal O}_X.1.$
 Therefore if 
$$(W,w):=({\cal O}_X.1\oplus\Lambda^2(V)\otimes I^{-1},1),$$
 then  by the above identification $A$ induces an element of $\spazuWw(X)$, 
 and since 
$$\Upsilon:\hbox{\rm Bil}_{(V,I)}
\cong \spazuWw$$ is an $X$-isomorphism 
 by (1), Theorem \ref{bilinear-forms-as-specialisations}, it follows that 
 there exists an $I$-valued 
  global quadratic form $q=q_b$ induced from an $I$-valued global
  bilinear form 
  $b$ on $V$ such that the algebra structure $\Upsilon(b)\cong A.$ 
  (We recall that 
 $\Upsilon(b)$ is the algebra structure induced from the 
 linear isomorphism  
$$\psi_b:C_0(V,q=q_b, I)\cong {\cal O}_X.1\oplus
 \Lambda^2(V)\otimes I^{-1}=W$$ of 
 (2d), Theorem \ref{Bourbaki}, which preserves 1 by (2a) of the same 
 Theorem). The proof of (b3) follows from a combination of 
 those of (b1) and (b2). {\bf Q.E.D., 
Theorem \ref{structure-of-specialisation} and 
 surjectivity of Theorem \ref{bijectivity}.}

\section{Specialised Algebras on Self-Dual 
  Bundles}
\label{sec6}
 In this section, we investigate the 
 specialised algebras when the underlying bundle is self-dual and 
 prove Theorem~\ref{surjectivity-affine-limited}.
We first have the following general 
result, which yields part of the assertions in 
 (a), Theorem \ref{surjectivity-affine-limited} that concern only 
 rank 4 vector bundles. 
\begin{proposition} Let $X$ be any scheme, $n$ an integer $\geq 2$
  and $W$ a rank $n^2$ vector
  bundle over $X$ with nowhere-vanishing global section $w.$  Then  
 $X$ is irreducible (resp. reduced) iff $\azuWw$ is irreducible
  (resp. reduced) iff 
 $\spazuWw$ is irreducible (resp. reduced). \label{partial-claim} 
\end{proposition}
\paragraph*{\sc Proof:} 
We first make certain observations when 
  $X$ is any reduced scheme and $W$ is a rank $n^2$ vector bundle 
 over $X$ with nowhere-vanishing global section $w.$ It is not hard to 
 see that the structure morphism 
$$\azuWw\longrightarrow X$$ is in fact a 
 morphism of finite presentation. Hence, in view of Prop.17.5.7, EGA IV 
 \cite{ega4}, assertion  (3) of Theorem 
\ref{repbility-smoothness-azu} implies that $\azuWw$ is reduced. 
 Since $\spazuWw$ is the schematic image of $\azuWw$, it follows from 
 (1), Prop.\ref{sch-img-props} that $\spazuWw$ is reduced as 
 well.  
When  $X$ is integral, we next verify that 
 $\azuWw$ is also integral. Consider any affine open subscheme 
 $$U=\hbox{\rm Spec}(R)\hookrightarrow X$$ 
such that $W|U$ is trivial and $w|U$ is part of a global 
 basis. There are $(w|U)$-unital Azumaya algebra structures on 
 $W|U$ which are isomorphic to the $(n\times n)$-matrix algebra over 
 $R.$ Consider the orbit morphism 
$$\hbox{\rm Stab}_{w|U}\longrightarrow 
 \hbox{\rm Azu}_{W|U, w|U}$$ corresponding to 
 one such algebra structure. Assertions (2) and   
 (3) of Theorem \ref{repbility-smoothness-azu} show that the topological 
 image of this morphism is dense.
 Further, $\hbox{\rm Stab}_{w|U}$ is integral since it is 
$\cong {\mathbb A}^{12}_U$ and since $U$ is integral. 
 Thus 
 $$\hbox{\rm Azu}_{W|U, w|U}$$ is integral.
Recall (e.g., Prop.2.1.6 \& Cor.2.1.7, Chap.0.2, EGA I \cite{ega1}) 
 that a nonempty topological space, whose set of irreducible components 
is locally finite, is locally irreducible iff each 
 irreducible component is open; further  
 it is irreducible iff it is locally irreducible and 
 connected. Now 
 since $\azuWw$ can be covered by irreducible open subschemes which 
 pairwise intersect (as $X$ is irreducible), 
 it follows that $\azuWw$ is integral as well. 

Putting the above facts together with 
 (1), Prop.\ref{sch-img-props} shows that 
 $\spazuWw$ is also integral if $X$ is integral. 
 Using again the facts that $\azuWw$ behaves well under base change 
 (Theorem \ref{repbility-smoothness-azu}) and
 that the schematic image of an irreducible scheme is irreducible 
 (by (1), Prop.\ref{sch-img-props}), we may infer from the foregoing 
 that $X$ is irreducible iff $\azuWw$ is irreducible iff 
 $\spazuWw$ is irreducible. {\bf Q.E.D., Prop.\ref{partial-claim}.}
\paragraph*{\sc Proof of assertions ({\rm a})---({\rm d}) of 
 Theorem \ref{surjectivity-affine-limited}.}
The proofs  follow 
 essentially from the properties 
 of 
$$({\spadesuit})\hspace*{2cm} 
\spazuWw\longrightarrow X\textrm{ and }\azuWw\longrightarrow X$$
 as mentioned in Theorems \ref{smoothness-spazu} and
 \ref{repbility-smoothness-azu}. Part of the assertions in (a) are
 valid more generally and were proved in 
 Prop.\ref{partial-claim}. 
  We indicate  
 proofs for the not-so-obvious assertions,
 especially (c)  and for 
 the implication 
$$X\textrm{ irreducible }
 \Rightarrow D_X\textrm{ irreducible }$$ in statement (a). The 
 converse implication would follow from the fact that 
$D_X\longrightarrow X$ 
 is topologically surjective. Since the structure morphisms 
$(\spadesuit)$
 are smooth, 
 they are faithfully flat (hence surjective).  
 Therefore in view of (b), which is actually the result of good base-change 
 properties of $\spazuWw$ and $\azuWw$ relative to $X$, 
 the  irreducible components of $D_X$ and of $\spazuWw$ 
 are induced from those of $X$ by pullback---provided we check 
 the particular case when $X$ is irreducible and reduced. 
 
So let $X$ be integral and first assume that $W$ is trivial and $w$ is
 part of a global basis. Without loss of generality,  
 we may take 
$$(W,w)=(\Lambda^{even}(V), 1)$$ for $V$ a rank 3 
trivial vector bundle over $X.$ After fixing a suitable 
 basis for $V$, we may define the morphism $\Theta$, which 
 by Theorem \ref{Thetaisom} is an isomorphism 
 that maps the closed subset 
 $$(\hbox{\rm Quad}_V\times_X \hbox{\sf L}_w) \backslash(\hbox{\rm Quad}_V^{sr}\times_X\hbox{\sf L}_w)$$ onto
 $$D_0:=\spazuWw\backslash\azuWw.$$ Therefore the irreducibility of 
 $D_0$ is equivalent to that of 
  the closed subset
 $$\hbox{\rm Quad}_V\backslash\hbox{\rm Quad}_V^{sr}.$$ Recall 
 from the discussion on semiregular forms (page \pageref{semireg-bil-forms},
 Section \ref{sec2}) that the open subscheme $\hbox{\rm Quad}_V^{sr}$ 
 corresponds to localisation by the polynomial $P_3.$ This 
 polynomial is irreducible as an element of 
$$R[\zeta_i, \zeta_{ij}]\cong 
 R[\hbox{\rm Quad}_V]$$
 when $X=\hbox{\rm Spec}(R)$ is affine, though it is not 
 clear if it is a prime element (unless we assume something more 
 e.g., $R$ a UFD). The
 closed subset
 $$\hbox{\rm Quad}_V\backslash\hbox{\rm Quad}_V^{sr}$$ may be 
 given the canonical closed subscheme structure $Z(P_3)$ corresponding to 
 the vanishing of $P_3.$ Let $q^{(2)}$ be the global quadratic form 
 on $V$ given by 
$$(x_1, x_2, x_3)\mapsto x_1x_2.$$ It can be checked 
 that $q^{(2)}$ is not semiregular, but that its restriction to 
 the rank two (direct summand) vector subbundle generated by $\{e_1, e_2\}$
 is regular. Therefore the $X$-valued point corresponding to $q^{(2)}$ 
 lands topologically inside the closed subset underlying $Z(P_3).$
 Consider the orbit morphism 
$$O(q^{(2)}):\hbox{\rm GL}_V\longrightarrow 
 \hbox{\rm Quad}_V$$
 corresponding to this $X$-valued point which also  
 lands topologically inside $Z(P_3).$ It will follow from assertion 
 (2), Theorem \ref{stratification}, that the topological image of 
 $O(q^{(2)})$ is dense in $Z(P_3).$ On the other hand this topological 
 image is irreducible, since $\hbox{\rm GL}_V\cong {\mathbb A}^9_X.$ 
  It follows that $\hbox{\rm Quad}_V\backslash
\hbox{\rm Quad}_V^{sr}$, and hence $D_0$,
  is irreducible in the case when $W$ is trivial and 
 $w$ is part of a global basis over $X.$ Since 
 the reduced closed subscheme structure on $Z(P_3)$ is given by 
 the radical of the ideal $(P_3)$, it follows that $\textrm{rad}(P_3)$
 is the minimal prime divisor of $(P_3)$ and 
 by Krull's {\em Hauptidealsatz} this prime has height
 1. Therefore, the codimension of $D_0$ is also 1 in the present 
 case. 
 Now consider the case of a general $(W,w)$, choose an affine open covering 
  $${\{U_i=\hbox{\rm Spec}(R_i)\}}_{i\in I}$$ 
 of $X$ such that $W_i:=W|U_i$ is trivial and $w_i:=w|U_i$ is part 
 of a global basis $\forall\thinspace i.$  The subset  
  $$D_i:=\spazuWiwi\backslash\azuWiwi$$ is irreducible for each $i$ by 
 the preceding discussion. 
 On the other hand, by Theorems \ref{repbility-smoothness-azu}
 and  \ref{smoothness-spazu}, the subsets $D_i$ form 
 an open cover of 
$$D:=\spazuWw\backslash\azuWw.$$ 
 Since $X$ is irreducible, we
 have $D_i\cap D_j\neq\emptyset$  when $i\neq j.$
  Thus $D$ is 
 locally irreducible and connected, and hence irreducible (for e.g., by 
 Cor.2.1.7, Chap.0, EGA I \cite{ega1}).   
 Since $X$ is integral and noetherian, 
 assertion (2) of Theorem \ref{smoothness-spazu} implies that $\spazuWw$ is 
 integral and noetherian as well.
  Therefore the codimension of $D$ in 
 $\spazuWw$ is atleast 1. On the other hand (for e.g., 
 by Prop.14.2.3, Chap.0,  EGA IV, \cite{ega4})  
 this codimension is bounded above by 
$$1=\hbox{\rm Codim}\thinspace
(D_i, \spazuWiwi)$$
 for any $i$, since 
 $$\spazuWiwi\hookrightarrow\spazuWw$$
 is an open subset whose intersection with $D$ is precisely 
 $D_i.$ 
\paragraph*{\sc Proof of assertion ({\rm e}) of 
 Theorem \ref{surjectivity-affine-limited}.}
Let $X$ be an integral separated quasi-compact scheme and let 
$\pi:\spazuWw\rightarrow X$ be the structure morphism. We need the
following result. 
\begin{lemma}\label{ho-gm-isom}
The canonical map $H^0(X, {\cal O}^*)\rightarrow H^0(\spazuWw, {\cal
  O}^*)$ induced by $\pi$ is an isomorphism. 
\end{lemma}
For the proof, we first consider the case when $X=\textrm{Spec}(R)$ is
affine, $W$ is trivial and further $w$ is part of a global basis for
$W.$ In this case, by Theorem~\ref{Thetaisom}, $\spazuWw$ may be
identified with $${\mathbb A}^9_X=\textrm{Spec}(R[X_1,\ldots, X_9])$$ 
and the assertion of the Lemma is trivial since
$R[X_1,\ldots,X_9]^*=R^*$, $R$ being an integral domain.

 Since $X$ is quasi-compact, we
can always find a finite affine open covering which trivializes $W$ and
restricted to each member of which $w$ becomes part of a global basis.
We now proceed by induction on
the number $n$ of open sets in such an affine covering ${\cal V}$ of
$X.$ The case $n=1$ was dealt with in the previous paragraph. Consider
$n=2$ i.e., 
let ${\cal V}=\{V_0, V_1\}$ be a covering by two members.
 We have the following natural commutative diagram with
exact rows and with the vertical arrows given by the natural maps
induced by $\pi$, where ${\cal U}:=\pi^{-1}{\cal V}$ is the Zariski-open
affine covering of $\spazuWw$ induced from ${\cal V}$ via $\pi.$ 
We write $H^0(-)$ for $H^0(-, {\cal O}^*)$ for typographical reasons. 
$$\begin{CD}
0@>>>H^0(\spazuWw)@>>>H^0(U_0)\times H^0(U_1)@>>>H^0(U_{01})\\
& & @AAA @AAA @AAA & & \\
0@>>>H^0(X)@>>>H^0(V_0)\times H^0(V_1)@>>>H^0(V_{01})\\
\end{CD}$$
Here $V_{01}:=V_0\cap V_1$ and $U_{01}:=\pi^{-1}(V_{01}).$ Note that
both $V_{01}$ and $U_{01}$ are affine since $X$ is separated and since
$\pi$ is an affine morphism. Also note that $U_i$ (resp.~$U_{ij}$) may
be identified with $\textrm{SpAzu}_{W_i, w_i}$ (resp.~with
$\textrm{SpAzu}_{W_{ij},w_{ij}}$) due to the good base-change property
(3) of Theorem~\ref{smoothness-spazu}, where $W_i:=W|V_i$ etc. By the
previous paragraph, the middle and right vertical arrows are
isomorphisms, hence so is the first.
Suppose therefore that Lemma~\ref{ho-gm-isom} is true
for any integral separated quasi-compact $X$
 with an affine open covering of cardinality $\leq n$ which
 trivializes $W$ and to each member of which $w$ restricts to part of
 a basis.  Let
$${\cal V}=\{V'_i|0\leq i\leq n\}$$ be such a covering with $(n+1)$ members
and let ${\cal U}:=\pi^{-1}{\cal V}$ be the covering of $\spazuWw$
induced by $\pi$ from ${\cal V}.$ Let $V_0$ denote the union of the
first $n$ open sets $V'_i$, $U_0:=\pi^{-1}(V_0)$, $V_1:=V'_n$ and
$U_1:=\pi^{-1}(V_1).$ Then writing out a commutative diagram as done
above for the case of a covering with 2 members, and using the
induction hypotheses (for $n$ for $V_0$, for $n=1$ for $V_1$ and again
for $n$ for $V_{01}$), we conclude the proof of
Lemma~\ref{ho-gm-isom}.\\[1mm]

We now proceed with the proof of assertion (e) of
Theorem~\ref{surjectivity-affine-limited}. So let $X$ be
locally-factorial. We want to show that $\textrm{Pic}(\pi)$ is an
isomorphism. Recall that we have canonical isomorphisms
$\textrm{Pic}(T)\equiv H^1(T, {\cal O}_T^*)$ (Ex.4.5, Ch.III,
  \cite{Hartshorne}) and $\textrm{Pic}(T)\equiv \textrm{CL}(T)$
 for $T$ locally-factorial 
  (Cor.6.16, Ch.II, \cite{Hartshorne}). Here $\textrm{CL}$ denotes
  the Weil divisor class group. We shall use these identifications in
  what follows without particular mention. 

 Consider first the case when $X=\textrm{Spec}(R)$ is
affine, $W$ is trivial and further $w$ is part of a global basis for
$W.$ In this case, by Theorem~\ref{Thetaisom}, $\spazuWw$ may be
identified with ${\mathbb A}^9_X$ 
and so the assertion (e) is a consequence of Prop.6.6, Ch.II,
\cite{Hartshorne}.
 
Next, let ${\cal V}=\{V_0, V_1\}$ be a Zariski-open affine covering of
$X$ which trivializes $W$ and makes $w|V_i$ part of a global basis for
$W|V_i$ ($i=0,1$). We have the following commutative diagram with
exact rows given by Mayer-Vietoris sequences (Ex.2.24, Ch.III,
\cite{Milne})  
and with the vertical arrows given by the natural maps
induced by $\pi$, where ${\cal U}:=\pi^{-1}{\cal V}$ is the Zariski-open
affine covering of $\spazuWw$ induced from ${\cal V}$ via $\pi.$ 
We again write $H^i(-)$ for $H^i(-,{\cal O}^*)$ for typographical
reasons. 
$$\begin{CD}
H^0(U_0)\times H^0(U_1)
@>>>H^0(U_{01})
@>>>H^1(\spazuWw)
@>(=)>>
\\
  @AAA @AAA @AAA &  \\
H^0(V_0)\times H^0(V_1)
@>>>H^0(V_{01})
@>>>H^1(X)
@>(=)>>
\end{CD}$$
$$\begin{CD}
@>(=)>>H^1(\spazuWw)
@>>>H^1(U_0)\times H^1(U_1)
@>>>H^1(U_{01})
\\
& & @AAA @AAA @AAA &  \\
@>(=)>>H^1(X)
@>>>H^1(V_0)\times H^1(V_1)
@>>>H^1(V_{01})
\end{CD}$$
Here $V_{01}:=V_0\cap V_1$ and $U_{01}:=\pi^{-1}(V_{01}).$
By Lemma~\ref{ho-gm-isom}, the first and second vertical arrows are
isomorhisms. By the previous paragraph, the fourth and fifth vertical
arrows are also isomorphisms. Therefore the 5-Lemma implies that the
central vertical arrow is an isomorphism too. We now conclude the
proof of assertion (e) by induction on the number of members in an
open covering of $X$ (just as in the proof of
Lemma~\ref{ho-gm-isom}). 

\paragraph*{\sc Proof of assertion ({\rm f}) of 
 Theorem \ref{surjectivity-affine-limited}.} Let us remind the reader 
 that the pullbacks 
$$W':=W\otimes_X\spazuWw\textrm{ and }W\otimes_X\azuWw$$
 are naturally 
 endowed with $(w\otimes_X\spazuWw)$-unital associative
 algebra structures $A'$ and $A'|\azuWw$ with which they become respectively  
 the universal specialisation and universal Azumaya algebra relative 
 to the pair $(W,w)$---for details see the proof of Theorem 3.4, Part I,
 \cite{tevb-paper1}.  
By (c) and (d), the natural map 
$$\hbox{\rm Pic}(\spazuWw)\longrightarrow \hbox{\rm Pic}(\azuWw)$$ is 
 surjective and its kernel is generated by the image of ${\mathbb Z}.(D_X).$
 Given an isomorphism 
$$\phi:W^\vee\stackrel{\cong}{\longrightarrow}W,$$ we consider its pullback 
 $\phi'$ to $\spazuWw$, 
 and let 
$$(X', A', \phi'):=(\spazuWw, W\otimes_X\spazuWw,
 \phi\otimes_X\spazuWw).$$
 Note that $D_X=X'\backslash U(X',A').$
 Consider the composition of 
 the following morphisms of vector bundles: 
$$W'\otimes_{X'} W'\equiv A'\otimes_{X'} {A'}^{op}
\stackrel{(a\otimes b)\mapsto 
 (x\mapsto axb)}{\longrightarrow} \hbox{\rm End}(A')
\equiv {W'}^\vee\otimes_{X'} W'
\stackrel{\phi'\times Id}{\longrightarrow} W'\otimes_{X'} W'.$$ 
The above composite gives an endomorphism of the vector bundle $W'\otimes W'$ 
 which is an isomorphism precisely at the local rings of the points of 
 $U(X',A').$ Therefore, the induced element 
$$s(\phi')\in\Gamma(X', {\cal O}_{X'})
\equiv \hbox{\rm H}^0(X', \hbox{\rm End}(\hbox{\rm det}(W'\otimes W')))$$ 
 goes into the maximal ideal of the local ring at each point of $D(X',A')$ and 
 to a unit of the local ring at each point of $U(X',A').$ It follows 
 that $n.(D_X)$ is principal for some $n\geq 1.$ 
Now 
$A'|\azuWw$ is an Azumaya algebra bundle,
 and as seen in Theorem \ref{Max-Knus-Theorem} (cf. also 
the discussion following Prop.~\ref{spln-is-quaternion})  
 it is isomorphic to the even Clifford algebra of a canonically 
 obtained rank 3 quadratic bundle  on $\azuWw$ with values in the 
 structure sheaf. Therefore it follows 
 from (b1), Theorem \ref{structure-of-specialisation} that 
 $$\hbox{\rm det}(W\otimes_X\spazuWw)$$ maps to an element of 
 $2.\hbox{\rm Pic}(\azuWw).$ The last assertion in (f) is now 
 again a consequence of (b1), Theorem \ref{structure-of-specialisation}.   
\paragraph*{\sc Proof of assertion ({\rm g}) of 
 Theorem \ref{surjectivity-affine-limited}.}
Suppose that $s$ is a section to $\pi':\azuWw\rightarrow X$
corresponding to an Azumaya algebra structure on $W$ with identity
$w.$ By assertions (e) and (f) (notice that $W$ is self-dual!) of
Theorem~\ref{surjectivity-affine-limited}, the homomorphism 
$$\textrm{Pic}(\pi'):\textrm{Pic}(X)\rightarrow\textrm{Pic}(\azuWw)$$
is surjective. But $\pi'\circ s=\textrm{Identity}$ implies that
$\textrm{Pic}(\pi')$ is injective as well. The rest of the assertions
in (g) now follow from (f). 
{\bf Q.E.D., Theorem \ref{surjectivity-affine-limited}.}

\section{Stratification of the Variety of Specialisations}
\label{sec7}
In this section, we prove Prop.\ref{semilocal-unique-semireg} and 
 Theorem \ref{stratification}. 
\paragraph*{\sc Proof of Prop.\ref{semilocal-unique-semireg}.} 
Fix an $S$-basis  $\{e_1, e_2, e_3\}$ for 
 $V$, and with respect to this basis, let $q^1$ denote the
 quadratic form given by 
$$(x_1e_1+x_2e_2+x_3e_3)\longmapsto x_1x_2+x_3^2.$$  
 It is easy to see that this quadratic form is semiregular. 
 We show that any 
 semiregular quadratic form $q$ can be moved to $q^1$ i.e., that $\exists\  
 g\in \hbox{\rm GL}(V)\hbox{ such that } g\cdot q=q^1.$ 
   By Prop.3.17, Chap.IV, \cite{Knus}, there exists 
 a basis  $\{e'_1, e'_2, e'_3\}$ for 
 $V$ such that $q$ restricted to the submodule
 generated by $e'_1$ and $e'_2$ is {\em regular} and further such that
 $$q(e'_3)\in S^{*},\ 
b_q(e'_1, e'_2)=1 
\textrm{  and  }
 b_q(e'_1, e'_3)=0=b_q(e'_2, e'_3).$$
  Let $g'\in \hbox{\rm GL}(V)$
 be the automorphism  that maps $e'_i$ onto the 
 $e_i$ for each $i$ and consider the quadratic form $q':=g'\cdot q.$ 
 Then by definition of the $\hbox{\rm GL}(V)$-action on the 
 set $\hbox{\rm Quad}(V)$ of quadratic $S$-forms on $V$ we have  
$$q'(e_3)\in S^{*},\ 
b_{q'}(e_1, e_2)=1 
\textrm{  and  }
 b_{q'}(e_1, e_3)=0=b_{q'}(e_2, e_3).$$
 So if we assume that
 $$q'(e_i)=\lambda_i\ (\Leftrightarrow q(e'_i)=\lambda_i),$$ 
then we would have 
$$q'(x_1e_1+x_2e_2+x_3e_3)=\lambda_1 x_1^2+\lambda_2x_2^2+
\lambda_3x_3^2+x_1x_2\ \forall\ x_i\in S.$$
 Thus it is enough to show that $q'$ can be moved to $q^1.$ 
 We look for  an invertible matrix 
$$g''=(u_{ij})\in
 \hbox{\rm GL}(V)\textrm{ such that }g''\cdot q'=q^1.$$
 Writing this condition equivalently as
 $q'={(g'')}^{-1}\cdot q^1$ and comparing the polynomials in the $x_i$ gives 
  6 equations in terms of the $u_{ij}$ and the $\lambda_i$ which are to 
 be satisfied. We choose the $u_{ij}$ as follows. 
First set 
$$u_{11}=u_{22}=0\textrm{ and  let }u_{12}\in S^{*}$$ be a free parameter.  
Since every element of  $S$  has square roots in $S$, it makes sense to 
 choose 
 $$u_{31}=\pm\sqrt{\lambda_1}\textrm{ and }u_{32}=\pm\sqrt{\lambda_2}.$$
 We let 
 $$\alpha=1+2u_{31}u_{32},\ \beta=1-2u_{31}u_{32}\textrm{ and 
  }u_{21}=\beta/u_{12}.$$ 
 Since $q'=g\cdot q$ and since $q$ is semiregular,  
  $q'$ is also semiregular. Its half-discriminant relative to the present 
 basis of $V$ is  (remembering that $\lambda_3\in S^{*}$)   
 $$d_{q'}(e_1, e_2, e_3)=\lambda_3.(4\lambda_1\lambda_2-1)\in S^{*}.$$
 This implies that 
 $$\alpha\beta=1-4\lambda_1\lambda_2\in S^{*}\Longrightarrow \alpha,\beta
\in S^{*}.$$
 Therefore it makes sense to define 
$$u_{33}=\pm\sqrt{(\beta\lambda_3)/\alpha}, 
 u_{13}=-2u_{31}u_{33}u_{12}/\beta\textrm{   and   } 
 u_{23}=-2u_{32}u_{33}/u_{12}.$$
Note that $u_{33}\in S^{*}.$ 
A computation shows 
 that the determinant of the matrix $g''=(u_{ij})$ defined above
 is $-u_{33}\alpha\in S^{*}$ and hence $g''$ is invertible. 
It is also easily  checked that $g''\cdot q'=q^1.$  
{\bf Q.E.D., Prop.\ref{semilocal-unique-semireg}.}
\paragraph*{\sc Proof of assertion (2) of Theorem \ref{stratification}.}
(We shall not prove assertion (1) since it is well-known).
 Recall 
 from the discussion on semiregular forms (page \pageref{semireg-bil-forms},
 Section \ref{sec4}) that the open subscheme $\hbox{\rm Quad}_V^{sr}$ 
 corresponds to localisation by the polynomial $P_3$ and that this 
 polynomial is prime as an element of 
$$k[\zeta_i, \zeta_{ij}]\cong 
 k[\hbox{\rm Quad}_V].$$
 Here a quadratic form $q$ corresponding to 
 the point $(\lambda_i,\lambda_{ij})\in{\mathbb A}^6_k$ is 
 given by 
$$(x_1,x_2,x_3)\mapsto \Sigma_i \lambda_ix_i^2
+\Sigma_{i<j}\lambda_{ij}x_ix_j.$$
 For 
 ease of readability (and typesetting !) let us denote the closure
 $\overline{T}$ of a subset $T$  
 (given the reduced closed subscheme structure)   
 by $\closure{T}$ in what follows.  
Since $\hbox{\rm Quad}_V^{(1)}$ is the same as the variety underlying 
the open subscheme $\hbox{\rm Quad}_V^{sr}$ of semiregular quadratic 
 forms, that 
$$\closure{\hbox{\rm Quad}_V^{(1)}}={\hbox{\rm Quad}_V}$$ follows 
 from the fact that ${\hbox{\rm Quad}_V}$ is irreducible. 
By assertion (1) of the present Theorem,
 ${\hbox{\rm Quad}_V}$ is the disjoint union of the $\hbox{\rm Quad}_V^{(i)}$; 
 therefore 
$$\closure{\hbox{\rm Quad}_V^{(1)}}\backslash\hbox{\rm Quad}_V^{(1)}$$
 is the disjoint union of $$\{\hbox{\rm Quad}_V^{(i)}|2\leq i\leq 4\}$$
 and also equals the closed subset $Z(P_3)$ defined by 
 the vanishing of $P_3.$ An explicit computation shows that 
the dimension of the  stabilizer of $q^{(2)}$ in $\hbox{\rm GL}(V)$ is 4.
Since $\hbox{\rm Quad}_V^{(2)}$ is an open dense subvariety of
 $$\closure{\hbox{\rm Quad}_V^{(2)}}
\subset V(P_3),$$ its closure is thus 5-dimensional. But since $P_3$ is an 
 irreducible polynomial, $Z(P_3)$ is also an irreducible 5-dimensional 
 subvariety. It  follows that 
$$\closure{\hbox{\rm Quad}_V^{(2)}}
=\closure{\hbox{\rm Quad}_V^{(1)}}\backslash\hbox{\rm Quad}_V^{(1)}.$$
Since $\hbox{\rm Quad}_V^{(2)}$ is smooth in its closure
 ($=Z(P_3)$ as seen above),  
 the singularities of its closure are contained in
 $$\hbox{\rm Quad}_V^{(3)}\cup \hbox{\rm Quad}_V^{(4)}$$
 which consists of quadratic forms that are perfect squares i.e., squares 
 of linear forms. 
 These singularities may be identified with 
  points 
$$(\lambda_i, \lambda_{ij})\in {\mathbb A}^6_k
\cong\hbox{\rm Quad}_V$$
 at which 
 all the partial derivatives of $P_3$ vanish. A simple computation shows that 
 this set is the symmetric determinantal variety given by the vanishing of 
 the $(2\times 2)$-minors of the matrix of the symmetric bilinear form 
 associated to the generic quadratic form 
 given by 
$$(x_1,x_2,x_3)\mapsto \Sigma_i \zeta_ix_i^2
+\Sigma_{i<j}\zeta_{ij}x_ix_j;$$ but it  can also be 
 shown that this set precisely corresponds to the  
 perfect squares. Therefore
 $$\hbox{ Sing }(\closure{\hbox{\rm Quad}_V^{(2)}})
=\closure{\hbox{\rm Quad}_V^{(2)}}\backslash\hbox{\rm Quad}_V^{(2)}.$$
That  
$$\closure{\hbox{\rm Quad}_V^{(i+1)}}
=\closure{\hbox{\rm Quad}_V^{(i)}}\backslash
\hbox{\rm Quad}_V^{(i)}\textrm{ for }i=2$$
follows from the above and the obvious fact that any quadratic form 
 can be specialised to the zero quadratic form. The case $i=3$ is trivial.
To see that 
$$\closure{\hbox{\rm Quad}_V^{(3)}}$$ is smooth 
 if char($k$)=2, we first note from the above and assertion (1) of the 
 present Theorem that the closure of 
$\hbox{\rm Quad}_V^{(3)}$ consists
 of perfect squares. Since char($k$)=2,
 under the identification 
$${\hbox{\rm Quad}_V}\cong{\mathbb A}^6_k,$$ 
 the perfect squares are seen to correspond to the copy of 
 ${\mathbb A}_k^3$ in ${\mathbb A}_k^6$
 given by the vanishing of the last three coordinates 
 $\lambda_{ij},\ 1\leq i<j\leq 3.$ 
 To see that the zero quadratic form is a 
 singularity of
 $$\closure{\hbox{\rm Quad}_V^{(3)}}$$ if char($k$) $\neq 2$, we note 
 from the above that 
 $\closure{\hbox{\rm Quad}_V^{(3)}}$ is defined by the same equations 
 that define the singularities of $Z(P_3)$ and
 is a certain symmetric determinantal variety.
 It is a well-known fact---an application of Standard Monomial 
 Theory (for e.g., see \cite{laxmi-css} or \cite{deconcini-procesi})---that the
 ideal defining these equations 
 is itself {\em reduced}, i.e., it is the {ideal of the variety}.  
 Checking the Jacobian criterion now shows that the zero quadratic 
 form is indeed a singular  point. 
\paragraph*{\sc Proof of assertion (3) of Theorem \ref{stratification}.} The 
 proof will follow from a series of lemmas. Recall that on page 
\pageref{thetadef} we had defined  
the $X$-morphism 
$$\theta:\hbox{\rm Quad}_{V}
\longrightarrow \hbox{\rm Id-}w\hbox{\rm -Sp-Azu}_W$$ 
 which was used to define the 
 isomorphism $\Theta.$  
\begin{lemma}\label{Claim1}
There are exactly 4 $H$-orbits in \hbox{\rm SpAzu}.
\end{lemma}
By Theorem \ref{Thetaisom},
 any point of \hbox{\rm SpAzu}
 is of the form $\underline{t}\cdot \theta(q).$ Its $H$-orbit is 
 $H\cdot \theta(q).$ There exists a unique $i$, $(1\leq i\leq 4)$,
 and some $g\in\hbox{\rm GL}(V)$ 
  such that $q=g\cdot q^{(i)}.$ Consider the algebra isomorphism
 $C_0(g,1)$ of (1), 
 Prop.\ref{simil-induces-iso-of-even-cliff} and 
 the induced isomorphism $h(g,1,q^{(i)},q)$ of
 Theorem \ref{b-co-phi-lambda-indep-of-q}. By definition,  
\begin{equation}\begin{split}
\theta(q)&=\theta(g\cdot q^{(i)})=h(g,1, q^{(i)},q)\cdot
\theta(q^{(i)})
\notag \\
&\Rightarrow H\cdot \theta(q)= H\cdot \theta(q^{(i)})\notag\end{split}
\end{equation} so that 
 there are atmost 4 orbits. To see that there are atleast 4, we use  
 the following result. 
\begin{lemma}\label{Claim2} 
For $q\in \hbox{\rm Quad}_V$, we have
$$\hbox{\sf L}_w\cdot\Theta((GL(V)\cdot q)\times\{I_4\})=\Theta((GL(V)\cdot q) \times \hbox{\sf L}_w)
=H\cdot \theta(q).$$
\end{lemma}
On the one hand, we have 
\begin{equation}
\begin{split}
\Theta((GL(V)\cdot q) \times \hbox{\sf L}_w) &= \left\{ \underline{t}\cdot \theta( g\cdot q)
 \left|\right. \underline{t}\in \hbox{\sf L}_w \hbox{ and } g\in GL(V)\right\} \notag \\
 &= \left\{ \underline{t}\cdot(h(g,1,q,g\cdot q)\cdot \theta(q))
 \left|\right. \underline{t}\in \hbox{\sf L}_w \hbox{ and } g\in GL(V)\right\} \notag \\
 &\subset H\cdot \theta(q).\notag 
\end{split}
\end{equation}
Conversely take any  $h\cdot \theta(q)\in H\cdot \theta(q).$
 By Theorem~\ref{Thetaisom}, there exists a unique $\underline{t'}
\in\hbox{\sf L}_w$ and
 a unique $q'\in\hbox{\rm Quad}_V$ such that
 $$h\cdot \theta(q)=\underline{t'}\cdot \theta(q').$$
Therefore 
$$({(\underline{t'})}^{-1}.h)\cdot \theta(q)=\theta(q').$$
 Since $k$ is 
 algebraically- and hence quadratically closed, by (b), 
 Theorem \ref{lifting-of-isomorphisms},  
 there exists $g\in GL(V)$ such that 
 $q'=g\cdot q$ and 
$${(\underline{t'})}^{-1}.h=h(g,1,q,q').$$ Hence
$$h\cdot \theta(q)=\underline{t'}\cdot \theta(g\cdot q)
=\Theta(g\cdot q, \underline{t'})
\Rightarrow H\cdot \theta(q)\subset\Theta((GL(V)\cdot q)
\times \hbox{\sf L}_w).$$
This settles Lemma \ref{Claim2}.  As for Lemma \ref{Claim1}, 
if 
$q,q'\in\hbox{\rm Quad}_V$ are such that
 their $GL(V)$-orbits are distinct, then because 
 $\Theta$ is an isomorphism, Lemma \ref{Claim2} shows that 
 $H\cdot \theta(q)$ is distinct from $H\cdot \theta(q').$ 
\begin{lemma}
\label{Claim3} 
For each $i$, with $1\leq i\leq 4$,  
$$\closure{\hbox{\rm Quad}_V^{(i)}\times \hbox{\sf L}_w} 
= \closure{\hbox{\rm Quad}_V^{(i)}}\times \hbox{\sf L}_w.$$
\end{lemma}
If $f:X\longrightarrow Y$ is a {\em  smooth 
 morphism} and $U\hookrightarrow Y$ is an open subset, then 
 $$f^{-1}(\closure{U})=\closure{f^{-1}(U)}.$$
 Since 
$\hbox{\sf L}_w\longrightarrow 
\hbox{Spec}(k)$ is smooth, so is the induced morphism 
$$\closure{\hbox{\rm Quad}_V^{(i)}}\times\hbox{\sf L}_w\longrightarrow \closure{\hbox{\rm Quad}_V^{(i)}}.$$
Taking $f$ to be this morphism and 
$U=\hbox{\rm Quad}_V^{(i)}$ gives Lemma \ref{Claim3}. 
\begin{lemma}\label{Claim4}
  The $GL(V)$-stratification of
 $\hbox{\rm Quad}_V$ 
 induces a $GL(V)$-stratification of
 $$\hbox{\rm Quad}_V\times\hbox{\sf L}_w$$ (the 
 $GL(V)$-action on $\hbox{\sf L}_w$ taken to be trivial) with 
 strata given by 
$$(\hbox{\rm Quad}_V\times\hbox{\sf L}_w)^{(i)}:= \hbox{\rm
   Quad}_V^{(i)}\times\hbox{\sf L}_w,\ (1\leq i\leq 4).$$
\end{lemma}
To prove Lemma \ref{Claim4}, 
the only thing that needs to be checked is that 
 $$\closure{(\hbox{\rm Quad}_V\times\hbox{\sf L}_w)^{(i+1)}}
=\closure{(\hbox{\rm Quad}_V\times\hbox{\sf L}_w)^{(i)}}\backslash(\hbox{\rm Quad}_V\times\hbox{\sf L}_w)^{(i)}.$$
This follows by applying Lemma \ref{Claim3} twice:
\begin{equation}
\begin{split}
\closure{(\hbox{\rm Quad}_V\times\hbox{\sf L}_w)^{(i+1)}} &=
\closure{\hbox{\rm Quad}_V^{(i+1)}\times\hbox{\sf L}_w}\notag \\
&=  \closure{\hbox{\rm Quad}_V^{(i+1)}}\times\hbox{\sf L}_w \notag \\
&=  (\closure{\hbox{\rm Quad}_V^{(i)}}\backslash\hbox{\rm Quad}_V^{(i)})\times\hbox{\sf L}_w \notag \\
 &=  \closure{\hbox{\rm Quad}_V^{(i)}}\times\hbox{\sf
  L}_w\backslash\hbox{\rm Quad}_V^{(i)}\times\hbox{\sf L}_w \notag \\
&=  \closure{(\hbox{\rm Quad}_V\times\hbox{\sf L}_w)^{(i)}}\backslash(\hbox{\rm Quad}_V\times\hbox{\sf L}_w)^{(i)}. \notag 
\end{split}
\end{equation}
Now according to Lemma \ref{Claim2}, we have  
$$\hbox{\rm SpAzu}^{(i)}=\Theta(GL(V)\cdot q^{(i)} \times \hbox{\sf
 L}_w)=\Theta((\hbox{\rm Quad}_V\times\hbox{\sf L}_w)^{(i)}).$$
 This combined with 
Lemma \ref{Claim4}  and the fact that $\Theta$
 is an isomorphism (Theorem \ref{Thetaisom})
 completes  the proof of assertion (3) of Theorem 
 \ref{stratification}. {\bf Q.E.D., Theorem \ref{stratification}.}

\paragraph{Acknowledgements}\addcontentsline{toc}{section}{Acknowledgements}
 The author is grateful for the Postdoctoral Fellowship 
  of the {\em Graduierten\-kolleg Gruppen
 und Geometrie} 
 under support from the {\em Deutschen Forschungsgemeinschaft}  and the 
 State of Niedersachsen at the Mathematisches Institut G\"ottingen 
 where this paper was written.

 It is a pleasure to thank the 
 National Board for Higher Mathematics, Department of Atomic Energy, 
 Government of India, for its Postdoctoral Fellowship (August 2002--August 2003) at the 
 Chennai Mathematical Institute, Chennai, India, during the latter half of 
 which certain special cases of the results of this work were obtained.



\end{document}